\documentclass{article}
\usepackage[utf8]{inputenc}

\usepackage{graphicx}

\usepackage{mathtools}
\usepackage{amssymb}
\usepackage{amsmath}
\usepackage{amsthm}

\usepackage{caption}
\usepackage{subcaption}

\usepackage{booktabs}
\usepackage{comment}

\usepackage{cases}
\usepackage{hyperref}

\usepackage{float}
\usepackage{enumitem}

\usepackage{algorithm}
\usepackage{algorithmic}
\usepackage{verbatim}

\usepackage[margin=1.3cm]{geometry}

\newtheorem{theorem}{Theorem}

\newtheorem{lemma}[theorem]{Lemma}
\newtheorem{proposition}[theorem]{Proposition}

\newtheorem{assumption}[theorem]{Assumption}
\newtheorem{example}[theorem]{Example}

\theoremstyle{definition}

\theoremstyle{remark}
\newtheorem{remark}{Remark}

\newtheorem{innercustomgeneric}{\customgenericname}
\providecommand{\customgenericname}{}
\newcommand{\newcustomtheorem}[2]{%
  \newenvironment{#1}[1]
  {%
   \renewcommand\customgenericname{#2}%
   \renewcommand\theinnercustomgeneric{##1}%
   \innercustomgeneric
  }
  {\endinnercustomgeneric}
}

\newcustomtheorem{customthm}{Theorem}
\newcustomtheorem{customlemma}{Lemma}
\newcustomtheorem{customassump}{Assumption}
\newcustomtheorem{customnot}{Notation}

\def\NN{{\mathbb N}}
\def\ZZ{{\mathbb Z}}
\def\RR{{\mathbb R}}
\def\CC{{\mathbb C}}

\DeclareMathOperator{\sign}{sign}

\usepackage{tikz}

\usepackage{tocloft}

\newcommand{\nocontentsline}[3]{}
\let\origcontentsline\addcontentsline
\newcommand\stoptoc{\let\addcontentsline\nocontentsline}
\newcommand\resumetoc{\let\addcontentsline\origcontentsline}

\title{Computation of whispering gallery modes for spherical symmetric, heterogeneous Helmholtz problems with piecewise smooth refractive index}
\author{Bouchra Bensiali\footnote{(bouchra.bensiali@centrale-casablanca.ma), \'Ecole Centrale Casablanca, Ville Verte, 27182 Bouskoura, Morocco} \and Stefan Sauter\footnote{(stas@math.uzh.ch), Institut für Mathematik, Universität Zürich, Winterthurerstr 190, CH-8057 Zürich, Switzerland}}
\date{\today}

\begin{document}

\maketitle

\begin{abstract}
In this paper, we develop a numerical method for the computation of (quasi-)resonances in spherical symmetric, heterogeneous Helmholtz problems with piecewise smooth refractive index. Our focus lies in resonances very close to the real axis, which characterize the so-called whispering gallery modes. Our method involves a modal equation incorporating fundamental solutions to decoupled problems, extending the known modal equation to the case of piecewise smooth coefficients. We first establish the well-posedness of the fundamental system, then we formulate the problem of resonances as a nonlinear eigenvalue problem, whose determinant will be the modal equation in the piecewise smooth case. In combination with the numerical approximation of the fundamental solutions using a spectral method, we propose a Newton method to solve the nonlinear modal equation with a proper scaling. We prove the local convergence of the algorithm in the piecewise constant case by showing the simplicity of the roots. We confirm our approach through a series of numerical experiments in the piecewise constant as well as in the variable case.


\end{abstract}


\addtocontents{toc}{\protect\setcounter{tocdepth}{-1}}

\section{Introduction}

Whispering gallery modes (WGM) are one type of interference phenomena of waves characterized by their localization along the jump interface of the
wave speed. These modes occur in general at high frequencies and are associated with complex scattering resonances very close to the real axis.
These waves have many applications such as in whispering gallery modes
resonators in optics (with important applications in medicine) 
where the goal is to confine light waves in 
a local region of an ambient medium~\cite{ilchenko2006optical,balac2021asymptotics}. The determination of 
geometric and material configurations for such interference phenomena 
is an important task in order to improve the performance of the system. The
study of whispering gallery modes thus requires the study of high-frequency
scattering problems, that are usually modelled by the Helmholtz equation. In
this paper, we propose a numerical method to detect the critical
states associated with scattering resonances for heterogeneous Helmholtz
problems with piecewise smooth refractive index. The developed software allows for an efficient
computation of many (quasi-)resonances for spherical symmetric
problems with varying coefficients.

Different numerical methods have been suggested in the literature for
computing resonances in Helmholtz problems. Some of them rely on the
perfectly matched layer technique~\cite{moitier2019etude,kim2009computation}%
, while others are based on a Dirichlet-to-Neumann map~\cite%
{araujo2018efficient,araujo2017spurious}. Further approaches are boundary integral methods~\cite{steinbach2017combined,poisson1995etude,heider2010computation} and Hardy space methods~\cite{hohage2009hardy,nannen2013exact}.  All these methods have in common that the problem for resonances is
formulated as a 
non-selfadjoint 
or nonlinear eigenvalue problem.

After discretization, the nonlinear eigenvalue problem can be reduced to the general form: given an holomorphic matrix function $T\colon \Omega\to \CC^{n\times n}$ on an open set $\Omega\subset \CC$, find $(u,k)\in \CC^n\times \Omega$ with $u\ne 0$ such that
\begin{equation}
T(k)u=0.
\end{equation}
This type of nonlinear matrix eigenvalue problems can be solved for instance using a Newton method~\cite{kressner2009block} or a contour integration based method~\cite{asakura2009numerical,beyn2012integral,steinbach2017combined}, among other methods. While Newton's method consists in solving a nonlinear equation (for instance $\det(T(k))=0$), the contour integral method consists in finding all eigenvalues and eigenvectors inside a given contour in the complex plane by reducing the original nonlinear eigenvalue problem to a linear eigenvalue problem that has identical eigenvalues in the domain. It is well known that Newton-type methods are sensitive to the initial starting point and can thus miss some eigenvalues, moreover, different formulations of Newton-type methods may exhibit non-robust convergence behavior ~\cite{jarlebring2012convergence}. On the other hand, the contour integration method is free from fixed point iterations required in Newton's method; however, it relies on a considerable number of control parameters: the choice of the contour itself (the radius and midpoint if the contour is a circle), quadrature points for the approximation of the contour integral and thresholds for the singular value decomposition filtering. A comparison between the two methods thus requires an extensive analysis of numerical experiments to assess the robustness, accuracy and speed of convergence of each one of them; see the review papers of Voß~\cite{voss2013nonlinear}, and Güttel and Tisseur~\cite{guttel2017nonlinear}, which contain a presentation of various Newton-type methods and contour integration methods for certain types of nonlinear eigenvalue problems.

In the physical community, whispering gallery modes are also studied numerically, e.g., for the applications mentioned earlier. Numerical codes are used for this purpose, based on Finite Difference Time Domain~\cite{hagness1997fdtd} or Finite Elements~\cite{degtyarev20163d}. Another widely used physical approach is based on the coupled mode theory (CMT)~\cite{hiremath2006modeling,franchimon2013interaction}, but this method is limited to specific geometries and refractive indices where the modes can be explicitly computed. For complex systems or those with significant variations in material properties, alternative modeling approaches may be necessary to capture the full dynamics. For the computation of resonances associated with whispering gallery modes, it is also common in the physical literature to truncate the domain with non transparent boundary conditions~\cite{oxborrow2007traceable}. 
Other recent works provide asymptotic expansion of WGM resonances at high polar
frequency for cavities with radially varying optical index~\cite%
{balac2021asymptotics}.

Next, we will sketch our new approach for spherical symmetric, heterogeneous
Helmholtz problems with piecewise smooth refractive index. As the starting
point we write the Helmholtz equation in spherical coordinates so that the
equation for the radially dependent part becomes a Bessel-type ordinary
differential equation with a transmission condition at the interface point.
The solution is a linear combination of two linear independent solutions
(fundamental system) in each of the two sub-intervals and the linearity of
the problem allows us to express explicitly two of the four coefficients in
their linear combinations by the boundary conditions. The remaining two
coefficients are determined by the transmission conditions and are the
solution of a $2\times2$ system whose entries depend on the fundamental
system and on the wavenumber. Since also the fundamental system depends on
the wavenumber, the problem is highly non-linear. For varying refractive index
the fundamental system as well as its dependence on the wavenumber is not
known explicitly. To determine the resonances, i.e., the wavenumbers such that
the determinant of the $2\times2$ system is zero we propose a Newton method
which takes into account the fact that, for changing wavenumber in the Newton
iteration, also the corresponding fundamental systems have to be updated
(numerically) along their derivatives with respect to the wavenumber. The
resulting algorithm allows us to compute the WGM for varying refractive index
to high accuracy for a large range of modes.

The paper is organized as follows. After presenting the problem setting, the formulation of the problem of resonances as a nonlinear eigenvalue problem is presented and analyzed using a system of fundamental solutions in Section~\ref{sec:problemformulation}. In Section~\ref{sec:newton}, the Newton method is presented and the local convergence result under suitable assumptions. We investigate in Section~\ref{sec:simplicity} the simplicity of the roots, and in Section~\ref{sec:scaling} the impact of scaling on the Newton algorithm. Finally, we present our final approach in Section~\ref{sec:pwscase} through targeted numerical experiments, before concluding remarks. By concentrating on Newton's approach, this work provides a focused analysis that can serve as a basis for more extensive methodological benchmarks in future research.

\section{Problem formulation}\label{sec:problemformulation}

\subsection{Polar coordinates and fundamental systems}

We consider the following heterogeneous Helmholtz problem in a spherical symmetric setting in dimension $d=2$:
\begin{equation}\label{eq:helmholtz}
\begin{dcases}
-\Delta u - (k n)^2 u=0 & \text{in $\Omega:=B_1^2:=\{x\in \RR^2, |x|<1\}$},\\
\frac{\partial u}{\partial \nu}-T_{k \hat{n}(1)} u=g & \text{on $\gamma_1:=\partial\Omega$},
\end{dcases}
\end{equation}
where $k\in \CC^*$ (set of nonzero complex numbers) denotes the wavenumber, and $n^2$ the refractive index. In the following we restrict to piecewise Lipschitz continuous coefficient $n$
and formulate these equations as a transmission problem. We will need the sets $\CC^*_{\ge 0}:=\{\zeta\in\CC^* \mid \operatorname{Im} \zeta \ge 0\}$ and $\CC_{< 0}:=\{\zeta\in\CC^* \mid \operatorname{Im} \zeta < 0\}$. Here $\partial/\partial\nu$ is the (outward) normal derivative at
$\partial\Omega$. Since we are interested in spherical symmetric problems we
assume that the refractive index only varies in radial direction: $n\left(
x\right)  =\hat{n}\left(  \left\vert x\right\vert \right)  $ for some
univariate function $\hat{n}\in L^{\infty}\left(  \tau\right)  $,
$\tau:=\left]  0,1\right[  $, which is piecewise smooth and positive. More
precisely, we assume that there exists a \textit{jump point} $\xi\in\left]
0,1\right[  $ which divides $\tau$ into the subintervals $\tau_{1}:=\left]
0,\xi\right[  $, $\tau_{2}:=\left]  \xi,1\right[  $ and the domain into
$\Omega_{j}:=\left\{  x\in\mathbb{R}^{2}\mid\left\vert x\right\vert \in
\tau_{j}\right\}  $, $j\in\left\{  1,2\right\}  $. The interface is denoted by $\gamma_\xi
:=\overline{\Omega_{1}}\cap\overline{\Omega_{2}}$. The exterior complement of
$\Omega$ is $\Omega^{+}:=\mathbb{R}^{2}\backslash\overline{\Omega}$. We assume
that there are positive constants $n_{\min}$, $n_{\max}$ such that the
restrictions of $\hat{n}$ to the subintervals $\tau_{1}
$, $\tau_{2}$ satisfy%
\begin{align*}
	\hat{n}_{j}  & :=  \hat{n}\vert_{\tau_{j}}\in C^{\infty}\left(\tau_{j}\right)  \\
	0  & <n_{\min}:=\inf_{r\in \,{]0},1[  }\hat{n}\left(  r\right) \leq\sup_{r\in \,{]0},1[  }\hat{n}\left(r\right)=:n_{\max}<\infty.
\end{align*}
Finally we set $n_{j}\left(  x\right)  :=\hat{n}_{j}\left(  \left\vert
x\right\vert \right)  $. To express functions in polar coordinates
$x=r\binom{\cos\theta}{\sin\theta}$ with $r:=\left\vert x\right\vert $ and
$\theta\in\left[  -\pi,\pi\right[  $ we use the \textquotedblleft hat
notation\textquotedblright: $\hat{u}\left(  r,\theta\right)  :=u\left(
r\binom{\cos\theta}{\sin\theta}\right)  $. If $\hat{u}$ is independent of
$\theta$ we write short $\hat{u}\left(  r\right)  $ for $\hat{u}\left(
r,\theta\right)  $.  We consider Dirichlet-to-Neumann boundary conditions, which are equivalent to an outgoing radiation condition (for $g=0$). 
We employ the definition of the Dirichlet-to-Neumann operator $T_{\kappa}$ on
the sphere $\gamma_{1}$ by a series expansion. For $\varphi\in H^{1/2}\left(
\gamma_{1}\right)  $ with Fourier expansion%
\[
\varphi\left(  x,y\right)  =\sum_{m\in\mathbb{Z}}\varphi_{m}\operatorname*{e}%
\nolimits^{\operatorname*{i}m\theta}\quad\text{for }\left(  x,y\right)
=\left(  \cos\theta,\sin\theta\right)  \in\gamma_{1}%
\]
the operator $T_{\kappa}$ is given by%
\begin{equation}
\left(  T_{\kappa}\varphi\right)  \left(  x,y\right)  =\sum_{m\in\mathbb{Z}%
}\kappa\frac{H_{m}^{\prime}\left(  \kappa\right)  }{H_{m}\left(
\kappa\right)  }\varphi_{m}\operatorname*{e}\nolimits^{\operatorname*{i}%
m\theta},\label{DefTk}%
\end{equation}
where $H_{m}\left(  r\right)  :=H_{m}^{\left(  1\right)  }\left(  r\right)  $
is the Hankel function of first kind and
order $m$~\cite{nedelec2001acoustic}. From
\cite[\S 15.7]{watson1922treatise} it follows that the total number of zeros of $H_{m}$
does not exceed $m+1$ and all of them are located in $\mathbb{C}_{<0}$. In
this way, $\mathcal{Z}:=\left\{  z\in\mathbb{C}\mid\exists m\in\mathbb{Z}\quad
H_{m}\left(  z\right)  =0\right\}  \subset\mathbb{C}_{<0}$ is countable and we
assume that $k\hat{n}\left(  1\right)  \notin\mathcal{Z}$ in~\eqref{eq:helmholtz}. It
is shown, e.g., in \cite{GrSa_DtN} that the series in (\ref{DefTk}) defines a
continuous mapping from $H^{1/2}\left(  \gamma_{1}\right)  \rightarrow
H^{-1/2}\left(  \gamma_{1}\right)  $ for all $\kappa\in\mathbb{C}^*_{\geq0}$.

We consider the nonlinear PDE eigenvalue problem associated with~\eqref{eq:helmholtz}
\begin{equation}\label{eq:pdenep}
\begin{dcases}
\text{Find $(k,u)\in\CC^*\times (H^1(\Omega)\setminus{\{0\}})$ such that}\\
-\Delta u = (k n)^2 u & \text{in $\Omega_1 \cup\Omega_2$},\\
[u]_{\gamma_\xi}=[\nabla u\cdot \nu]_{\gamma_\xi}=0, &\text{interface condition,}\\
\frac{\partial u}{\partial \nu}=T_{k \hat{n}(1)} u& \text{on $\gamma_1
	$}.
\end{dcases}
\end{equation}

\begin{proposition}
	\label{PropWellPosed}All solutions (resonances) of problem~\eqref{eq:pdenep}
	satisfy $k\in\mathbb{C}_{<0}$.
\end{proposition}

\proof

Assume by contradiction that there exists a non-trivial solution $\left(
k,u\right)  \in\mathbb{C}_{\geq0}^{\ast}\times H^{1}\left(  \Omega\right)
\backslash\left\{  0\right\}  $. Then, the weak form of~\eqref{eq:pdenep} reads
\[
a\left(  u,v\right)  :=\left(  \nabla u,\nabla v\right)  _{L^{2}\left(
	\Omega\right)  }+k^{2}\left(  n^2u,v\right)  _{L^{2}\left(  \Omega\right)
}-\left(  T_{k\hat{n}\left(  1\right)  }u,v\right)  _{L^{2}\left(
\partial\Omega\right)  }=0\quad\forall v\in H^{1}\left(  \Omega\right)  .
\]
We choose $v=\mu u$ with $\mu=-\operatorname*{i}k/\left\vert k\right\vert $
and obtain for the real part%
\begin{align*}
	\operatorname{Re}a\left(  u,\mu u\right)  =  & \frac{\operatorname{Im}%
		k}{\left\vert k\right\vert }\left(  \left\Vert \nabla u\right\Vert
	_{L^{2}\left(  \Omega\right)  }^{2}+\left\vert k\right\vert ^{2}\left\Vert
	nu\right\Vert _{L^{2}\left(  \Omega\right)  }^{2}\right)  \\
	& -\operatorname{Re}\left(  \overline{\mu}\left(  T_{k\hat{n}\left(  1\right)
	}u,u\right)  _{L^{2}\left(  \partial\Omega\right)  }\right)  =0\quad\forall
	v\in H^{1}\left(  \Omega\right)  .
\end{align*}
Since $\operatorname{Im}k\geq0$ we get%
\[
0=\operatorname{Re}a\left(  u,\mu u\right)  \geq-\operatorname{Re}\left(
\overline{\mu}\left(  T_{k\hat{n}\left(  1\right)  }u,u\right)  _{L^{2}\left(
	\partial\Omega\right)  }\right)  .
\]
From \cite[(3.4b,c)]{MelenkSauterMathComp} and \cite{GrSa_DtN} we conclude that%
\begin{equation}
-\operatorname{Re}\left(  \overline{\mu}\left(  T_{k\hat{n}\left(  1\right)
}u,u\right)  _{L^{2}\left(  \partial\Omega\right)  }\right)  \geq C\left(
\operatorname{Im}k\right)  \left\Vert u\right\Vert _{L^{2}\left(
	\partial\Omega\right)  }^{2}\quad\forall u\in H^{1/2}\left(  \partial
\Omega\right)  \label{dissipative}%
\end{equation}
so that $\left.  u\right\vert _{\partial\Omega}=0$. The homogeneous
$\operatorname*{DtN}$ boundary conditions on $\partial\Omega$ in~\eqref{eq:pdenep}
then imply $\partial u/\partial \nu=0$ on $\partial\Omega$. For
$\operatorname{Im}k=0$ we apply the unique continuation principle for
piecewise Lipschitz continuous refractive index (see \cite{KenigUCP},
\cite{Alessandrini_ucp} and, e.g., \cite[Thm. 2.4]{GrahamSauter_Helm}) and
obtain $u=0$ which is a contradiction. For $\operatorname{Im}k>0$ we use
(\ref{dissipative}) and conclude directly from%
\[
0=\operatorname{Re}a\left(  u,\mu u\right)  \geq\frac{\operatorname{Im}%
	k}{\left\vert k\right\vert }\left(  \left\Vert \nabla u\right\Vert
_{L^{2}\left(  \Omega\right)  }^{2}+\left\vert k\right\vert ^{2}\left\Vert
nu\right\Vert _{L^{2}\left(  \Omega\right)  }^{2}\right)
\]
that $u=0$ which again is a contradiction.%
\endproof

Using polar coordinates, the solution and the boundary data can be expanded as a Fourier series via the ansatz
\begin{subequations}
	\label{fourrierexp}
\end{subequations}%
\begin{equation}
\begin{dcases}
u(x,y)=\hat{u}(r,\theta)=\sum_{m\in\ZZ} \hat{u}_m(r) e^{im\theta}\\
g(x,y)=\hat{g}(\theta)=\sum_{m\in\ZZ} \hat{g}_m e^{im\theta}
\end{dcases} \tag{%
	\ref{fourrierexp}%
	a}\label{fourrierexp1}%
\end{equation}
which leads to a system of ODEs for the Fourier coefficients
\begin{equation}
\hat{u}_{j,m}:=\left.  \hat{u}_{m}\right\vert _{\tau_{j}}\quad j\in\left\{
1,2\right\}  . \tag{%
	\ref{fourrierexp}%
	b}\label{fourrierexp2}%
\end{equation}
Let the differential operator $L^m_{j,k}$ be given by%
\[
L^m_{j,k}v:=-\frac{1}{r}\left(  rv^{\prime}\left(  r\right)  \right)  ^{\prime
}+\left(  \frac{m^{2}}{r^{2}}-\left(  k\hat{n}_{j}\left(  r\right)  \right)
^{2}\right)  v\left(  r\right)  \quad\text{in }\tau_{j}%
\]
and the boundary operators 
by%
\[
B^m_{1,k}(0)v:=\begin{dcases}
v\left(  0\right)   & \text{for odd }m\text{,}\\
v^{\prime}\left(  0\right)   & \text{for even }m\text{,}%
\end{dcases}  \quad\text{and\quad}B^m_{2,k}(1)v:=v^{\prime}\left(  1\right)
-k\hat{n}_{2}\left(  1\right)  \frac{H_{m}^{\prime}\left(  k\hat{n}_{2}\left(  1\right)
	\right)  }{H_{m}\left(  k\hat{n}_{2}\left(  1\right)  \right)  }v\left(  1\right)  ,
\]
where we again assume that $k\hat{n}_{2}\left(
1\right)  $ does not belong to the set of zeros $\mathcal{Z}$ (see~\cite{sauter2021heterogeneous} for the distinction between odd and even $m$ for the boundary condition at $0$). In this paper we investigate resonances which are caused by the jump of the
refractive coefficient across the interface at $\xi$. 
Then, we consider the problem%
\begin{equation}%
\begin{dcases}
\begin{array}
[c]{cl}%
L^m_{j,k}\hat{u}_{j,m}=0 & \text{in }\tau_{j},\\
\left[  \hat{u}_m\right]  _{\xi}=\left[  \hat{u}_m^{\prime}\right]  _{\xi}=0 &
\text{transmission condition,}\\
B^m_{2,k}(1)\hat{u}_{2,m}=\hat{g}_m & \operatorname*{DtN}\text{ boundary conditions
	at }1,\\
B^m_{1,k}(0)\hat{u}_{1,m}=0 & \text{boundary condition at }0.
\end{array}
\end{dcases}
\label{eq:helmholtzpolarcordinates}%
\end{equation}

\begin{remark}
	The restriction to two-dimensional domains is only made to simplify the
	exposition. For general dimension, spherical coordinates can be used to
	transform~\eqref{eq:helmholtz} to a Bessel-type ordinary differential equation in
	analogy to~\eqref{eq:helmholtzpolarcordinates}, see, e.g., \cite[\S 2.2]{sauter2021heterogeneous}.
\end{remark}

In the next step we derive equations for a fundamental system of $L^m_{j,k}$. We define the functions $f^m_{j,\ell,k}$, $\left(  j,\ell\right)  \in\left\{
\left(  1,1\right)  ,\left(  2,1\right)  ,\left(  2,2\right)  \right\}  $ as the solutions of%
\begin{align}
	&
	\begin{dcases}
		\begin{array}
			[c]{lc}%
			L^m_{1,k}f^m_{1,1,k}=0 & \text{in }\tau_{1}\\
			B^m_{1,k}(0)f^m_{1,1,k}=0 & \\
			\left(  f_{1,1,k}^{m}\right)  ^{\prime}\left(  \xi\right)  -\beta_{1,k}%
			^{m}\left(  \xi\right)  f_{1,1,k}^{m}\left(  \xi\right)     =h_{1,1,k}%
			^{m}\neq0  &
		\end{array}
		\label{eq:f11k}
	\end{dcases}
\end{align}
\noindent\begin{minipage}{.5\linewidth}
	\begin{align}
		\begin{dcases}
			\begin{array}
				[c]{lc}%
				L^m_{2,k}f^m_{2,1,k}=0 & \text{in }\tau_{2}\\
				-\left(  f_{2,1,k}^{m}\right)  ^{\prime}\left(  \xi\right)  -\beta_{2,k}%
				^{m}\left(  \xi\right)  f_{2,1,k}^{m}\left(  \xi\right)     =0\\
				\left(  f_{2,1,k}^{m}\right)  ^{\prime}\left(  1\right)  -\beta_{2,k}%
				^{m}\left(  1\right)  f_{2,1,k}^{m}\left(  1\right)     =h_{2,1,k}^{m}\neq0
			\end{array}
			\label{eq:f21k}
		\end{dcases}
	\end{align}
\end{minipage}
\begin{minipage}{.5\linewidth}
	\begin{align}
		\begin{dcases}
			\begin{array}
				[c]{lc}%
				L^m_{2,k}f^m_{2,2,k}=0 & \text{in }\tau_{2}\\
				-\left(  f_{2,2,k}^{m}\right)  ^{\prime}\left(  \xi\right)  -\beta_{2,k}%
				^{m}\left(  \xi\right)  f_{2,2,k}^{m}\left(  \xi\right)     =h_{2,2,k}%
				^{m}\neq0\\
				\left(  f_{2,2,k}^{m}\right)  ^{\prime}\left(  1\right)  -\beta_{2,k}%
				^{m}\left(  1\right)  f_{2,2,k}^{m}\left(  1\right)     =0
			\end{array}
			\label{eq:f22k}%
		\end{dcases}
	\end{align}
\end{minipage}

\medskip
\noindent for multipliers $\beta_{j,k}^{m}\left(  \xi\right)  \in\mathbb{C}$ which are
fixed later. They should be chosen such that a) problems~\eqref{eq:f11k}-\eqref{eq:f22k} are well-posed and b) the dependence on the parameter $k$ is
smooth either on the whole $\mathbb{C}^{\ast}$, or on the half plane
$\mathbb{C}_{\geq0}^{\ast}$, or, more relevant in view of Proposition
\ref{PropWellPosed}, on the half plane $\mathbb{C}_{<0}$.

Let%
\[
\mathcal{J}:=\left\{  \left(  1,\xi\right)  ,\left(  2,\xi\right)  ,\left(
2,1\right)  \right\}
\]
and introduce for $\left(  j,x_{0}\right)  \in\mathcal{J}$, the operators
\begin{equation}
\mathfrak{R}_{j,k}\left(  x_{0}\right)  \left(  \hat{u}_{j}\right)
=\sum_{m\in\mathbb{Z}}\beta_{j,k}^{m}\left(  x_{0}\right)  \hat{u}%
_{j,m}\left(  x_0\right)  \operatorname*{e}\nolimits^{\operatorname*{i}%
	m\theta}.\label{DefRfrak}%
\end{equation}
The coefficients $h_{j,\ell,k}^{m}$ induce anti-linear forms via%
\begin{equation}
\mathfrak{F}_{j,\ell,k}\left(  v\right)  :=\left\{
\begin{array}
[c]{cc}%
\int_{\gamma_{\xi}}\mathfrak{h}_{1,1,k}\overline{v} & \text{for problem~\eqref{eq:f11k}, i.e.,
}\left(  j,\ell\right)  =\left(  1,1\right) \\
\int_{\gamma_{1}}\mathfrak{h}_{2,1,k}\overline{v} & \text{for problem~\eqref{eq:f21k}, i.e., }\left(
j,\ell\right)  =\left(  2,1\right)\\
\int_{\gamma_{\xi}}\mathfrak{h}_{2,2,k}\overline{v} & \text{for problem~\eqref{eq:f22k}, i.e.,
}\left(  j,\ell\right)  =\left(  2,2\right)%
\end{array}
\right.  \quad\text{for\quad}\mathfrak{h}_{j,\ell,k}:=\sum_{m\in\mathbb{Z}%
}h_{j,\ell,k}^{m}\operatorname*{e}\nolimits^{\operatorname*{i}m\theta
}.\label{DefFfrak}%
\end{equation}

\begin{assumption}
	\label{AssumptionCont}The coefficients $\beta_{j,k}^{m}\left(  \xi\right)  $,
	$\beta_{2,k}^{m}\left(  1\right)  $ are such that the resulting operators
	$\mathfrak{R}_{j,k}\left(  x_{0}\right)  :H^{1/2}\left(  \gamma_{x_{0}%
	}\right)  \rightarrow H^{-1/2}\left(  \gamma_{x_{0}}\right)  $ in
	(\ref{DefRfrak}) are bounded linear operators.
	
	The coefficients $h_{j,\ell,k}^{m}$ are such that the resulting anti-linear
	forms $\mathfrak{F}_{j,\ell,k}:H^{1/2}\left(  \gamma_{x_{0}}\right)
	\rightarrow\mathbb{C}$ in (\ref{DefFfrak}) are bounded.
\end{assumption}

\begin{example}
	\label{ExCasesI}\quad
	
	\begin{enumerate}
		\item[a.] The choice $\beta_{j,k}^{m}\left(  x_{0}\right)  =\operatorname*{i}%
		k \hat{n}_j(x_0)$ 
		for all $m\in\mathbb{Z}$ corresponds to Robin boundary conditions
		$\mathfrak{R}_{j,k}\left(  x_{0}\right)  u=\operatorname*{i}k\hat{n}_j(x_0)\left.
		u\right\vert _{\gamma_{x_{0}}}$ and
		
		\item[b.] the choice $\beta_{j,k}^{m}\left(  x_{0}\right)  =-\operatorname*{i}%
		k\hat{n}_j(x_0)$ to its flipped version $\mathfrak{R}_{j,k}\left(  x_{0}\right)
		u=-\operatorname*{i}k\hat{n}_j(x_0)\left.  u\right\vert _{\gamma_{x_{0}}}$.
		
		\item[c.] The choice $\beta_{2,k}^{m}\left(  1\right)  =k\hat{n}\left(
		1\right)  \frac{\left(  H_{m}^{\left(  1\right)  }\right)  ^{\prime}\left(
			k\hat{n}\left(  1\right)  \right)  }{H_{m}^{\left(  1\right)  }\left(  k\hat{n}\left(
			1\right)  \right)  }$ corresponds to the $\operatorname*{DtN}$ operator:
		$\mathfrak{R}_{2,k}\left(  1\right)  u=T_{k\hat{n}\left(  1\right)  }\left(
		\left.  u\right\vert _{\gamma_{1}}\right)  $.
	\end{enumerate}
\end{example}

\subsection{Well-posedness of ODEs for the fundamental system}\label{sec:wellposedness}
Next, we will investigate the well-posedness of the problems~\eqref{eq:f11k}-\eqref{eq:f22k}. The idea is to switch back from the ODE in polar
coordinates to an elliptic PDE on $\Omega_{j}$. We define the sesquilinear
form $a_{j}:H^{1}\left(  \Omega_{j}\right)  \times H^{1}\left(  \Omega
_{j}\right)  \rightarrow\mathbb{C}$ by%
\[
a_{j}\left(  u,v\right)  :=\left(  \nabla u,\nabla v\right)  _{\Omega_{j}%
}-k^{2}\left(  n^2u,v\right)  _{\Omega_{j}}-\left(  \mathfrak{R}_{j,k}\left(
\xi\right)  u,v\right)  _{\gamma_{\xi}}-\left\{
\begin{array}
[c]{ll}%
0 & \text{for }j=1,\\
\left(  \mathfrak{R}_{j,k}\left(  1\right)  u,v\right)  _{\gamma_{1}} &
\text{for }j=2,
\end{array}
\right.
\]
where $\left(  \cdot,\cdot\right)  _{\Omega_{j}}$ and $\left(  \cdot
,\cdot\right)  _{\gamma_{x}}$ denote the $L^{2}\left(  \Omega_{j}\right)  $
and $L^{2}\left(  \gamma_{x}\right)  $ scalar products. The problems~\eqref{eq:f11k}-\eqref{eq:f22k} are then equivalent to%
\begin{align*}
	\text{find }f_{1,1,k}  & \in H^{1}\left(  \Omega_{1}\right)  \text{ s. t.
	}a_{1}\left(  f_{1,1,k},v\right)  =\mathfrak{F}_{1,1,k}\left(  v\right)  \quad\forall
	v\in H^{1}\left(  \Omega_{1}\right)  ,\\
	\text{find }f_{2,1,k}  & \in H^{1}\left(  \Omega_{2}\right)  \text{ s. t.
	}a_{2}\left(  f_{2,1,k},v\right)  =\mathfrak{F}_{2,1,k}\left(  v\right)  \quad\forall
	v\in H^{1}\left(  \Omega_{2}\right)  ,\\
	\text{find }f_{2,2,k}  & \in H^{1}\left(  \Omega_{2}\right)  \text{ s. t.
	}a_{2}\left(  f_{2,2,k},v\right)  =\mathfrak{F}_{2,2,k}\left(  v\right)  \quad\forall
	v\in H^{1}\left(  \Omega_{2}\right)  .
\end{align*}
By choosing $v=\mu u$ for $\mu
:=-\operatorname*{i}k/\left\vert k\right\vert $ as a test function and
considering the real part we obtain for $j\in\left\{  1,2\right\}:$
\begin{equation}
0=\operatorname{Re}a_{j}\left(  u,\mu u\right)  =\frac{\operatorname{Im}%
	k}{\left\vert k\right\vert }\left(  \left\Vert \nabla u\right\Vert
_{\Omega_{j}}^{2}+\left\vert k\right\vert ^{2}\left\Vert nu\right\Vert
_{\Omega_{j}}^{2}\right)  -\sum_{x:\left(  j,x\right)  \in\mathcal{J}%
}\operatorname{Re}\left(  \bar{\mu}\left(  \mathfrak{R}_{j,k}\left(
x\right)  u,u\right)  _{\gamma_{x}}\right)  .\label{0ReCond}%
\end{equation}
Note that%
\[
-\operatorname{Re}\left(  \bar{\mu}\left(  \mathfrak{R}_{j,k}\left(  x\right)
u,u\right)  _{\gamma_{x}}\right)  =\frac{\operatorname{Re}k\operatorname{Im}%
	\left(  \mathfrak{R}_{j,k}\left(  x\right)  u,u\right)  _{\gamma_{x}%
	}-\operatorname{Im}k\operatorname{Re}\left(  \mathfrak{R}_{j,k}\left(
	x\right)  u,u\right)  _{\gamma_{x}}}{\left\vert k\right\vert }%
\]
for $\left(  j,x\right)  \in\mathcal{J}$.

\begin{lemma}
	\label{LemWellPosed}Let Assumption \ref{AssumptionCont} be satisfied.
	
	\begin{enumerate}
		\item Assume that there exists a non-negative function $\mu_{\operatorname*{R}%
		}:\mathbb{C}^*_{\ge0}\rightarrow\mathbb{R}_{\ge0}$ with $\mu_{\operatorname*{R}%
	}:\mathbb{C}_{>0}\rightarrow\mathbb{R}_{>0}$ such that for all $\left(
	j,x\right)  \in\mathcal{J}$, $k\in\mathbb{C}_{\geq0}^{\ast}$, and $u\in
	H^{1/2}\left(  \gamma_{x}\right)  \backslash\left\{  0\right\}  :$%
	\begin{equation}%
	\begin{array}
	[c]{rll}%
	-\operatorname{Re}\left(  \mathfrak{R}_{j,k}\left(  x\right)  u,u\right)
	_{\gamma_{x}}\geq & \mu_{\operatorname*{R}}\left(  k\right)  \left\Vert
	u\right\Vert _{\gamma_{x}}^{2}, & \text{for }\operatorname{Im}k>0,\\
	\quad &  & \\
	\left(  \operatorname*{sign}\operatorname{Re}k\right)  \operatorname{Im}%
	\left(  \mathfrak{R}_{j,k}\left(  x\right)  u,u\right)  _{\gamma_{x}}> & 0 &
	\text{for }\operatorname{Re}k\neq0.
	\end{array}
	\label{CondImkpos}%
	\end{equation}
	Then, problems~\eqref{eq:f11k}-\eqref{eq:f22k} are well-posed for any $k\in
	\mathbb{C}_{\geq0}^{\ast}$.
	
	\item Assume that there exists a non-negative function $\mu_{\operatorname*{R}%
	}:\mathbb{C}^*_{\le0}\rightarrow\mathbb{R}_{\ge0}$ with $\mu_{\operatorname*{R}%
}:\mathbb{C}_{<0}\rightarrow\mathbb{R}_{>0}$ such that for all $\left(
j,x\right)  \in\mathcal{J}$, $k\in\mathbb{C}_{<0}$ and $u\in H^{1/2}\left(
\gamma_{x}\right)  \backslash\left\{  0\right\}  $%
\begin{equation}%
\begin{array}
[c]{rll}%
-\operatorname{Re}\left(  \mathfrak{R}_{j,k}\left(  x\right)  u,u\right)
_{\gamma_{x}}\geq & \mu_{\operatorname*{R}}\left(  k\right)  \left\Vert
u\right\Vert _{\gamma_{x}}^{2} & \text{for }\operatorname{Im}k<0,\\
\quad &  & \\
-\left(  \operatorname*{sign}\operatorname{Re}k\right)  \operatorname{Im}%
\left(  \mathfrak{R}_{j,k}\left(  x\right)  u,u\right)  _{\gamma_{x}}> & 0 &
\text{for }\operatorname{Re}k\neq0.
\end{array}
\label{CondImkneg}%
\end{equation}
Then, problems~\eqref{eq:f11k}-\eqref{eq:f22k} are well-posed for any $k\in
\mathbb{C}_{<0}$.
\end{enumerate}
\end{lemma}

\begin{proof}
	It follows, e.g., as in \cite[Prop. 2.5]{doerfler_sauter} that the
	sesquilinear form $a_{j}\left(  \cdot,\cdot\right)  $ is continuous in
	$H^{1}\left(  \Omega_{j}\right)  $ and satisfies a G\aa rding inequality. The
	anti-linear forms $\mathfrak{F}_{j,\ell,k}$ are continuous so that
	well-posedness can be concluded from Fredholm's alternative via uniqueness,
	i.e., from the implication%
	\begin{equation}
	\left(  a_{j}\left(  u,v\right)  =0\quad\forall v\in H^{1}\left(  \Omega
	_{j}\right)  \right)  \implies\left(  u=0\right)  .\label{uniqueness}%
	\end{equation}
	We start with some consideration of the boundary terms in $a_{j}\left(
	\cdot,\cdot\right)  $.
	
	\textbf{First case: }Let $k\in\mathbb{C}_{\geq0}^{\ast}$ and the corresponding
	conditions in the lemma be satisfied.
	
	Then, for $\left(  j,x\right)  \in\mathcal{J}$ it holds%
	\[
	-\operatorname{Re}\left(  \bar{\mu}\left(  \mathfrak{R}_{j,k}\left(  x\right)
	u,u\right)  _{\gamma_{x}}\right)  \geq\frac{\operatorname{Re}k}{\left\vert
		k\right\vert }\operatorname{Im}\left(  \mathfrak{R}_{j,k}\left(  x\right)
	u,u\right)  _{\gamma_{x}}+\mu_{\operatorname*{R}}\left(  k\right)
	\frac{\operatorname{Im}k}{\left\vert k\right\vert }\left\Vert u\right\Vert
	_{\gamma_{x}}^{2},
	\]
	where we recall $\bar{\mu}=\operatorname*{i}\bar{k}/\left\vert k\right\vert $.
	It follows from (\ref{CondImkpos}) that both summands are non-negative and we
	split the analysis into two cases:
	
	\begin{enumerate}
		\item $\operatorname{Im}k>0$. Then%
		\begin{equation}
		-\operatorname{Re}\left(  \bar{\mu}\left(  \mathfrak{R}_{j,k}\left(  x\right)
		u,u\right)  _{\gamma_{x}}\right)  \geq\left(  \mu_{\operatorname*{R}}\left(
		k\right)  \frac{\operatorname{Im}k}{\left\vert k\right\vert }\right)
		\left\Vert u\right\Vert _{\gamma_{x}}^{2}.\label{Case1}%
		\end{equation}
		Since $\operatorname{Im}k>0$ we have $\mu_{\operatorname*{R}}\left(  k\right)
		>0$ so that the prefactor $\left(  \mu_{\operatorname*{R}}\left(
		k\right)  \frac{\operatorname{Im}k}{\left\vert k\right\vert }\right)  $ in the right-hand side of
		(\ref{Case1}) is positive.
		
		\item $\operatorname{Im}k=0$. Since $k\neq0$, we have $\operatorname{Re}%
		k\neq0$ and%
		\[
		-\operatorname{Re}\left(  \bar{\mu}\left(  \mathfrak{R}_{j,k}\left(  x\right)
		u,u\right)  _{\gamma_{x}}\right)  \geq\frac{\left\vert \operatorname{Re}%
			k\right\vert }{\left\vert k\right\vert }\left(  \operatorname*{sign}%
		\operatorname{Re}k\right)  \operatorname{Im}\left(  \mathfrak{R}_{j,k}\left(
		x\right)  u,u\right)  _{\gamma_{x}}>0
		\]
		for all $H^{1/2}\left(  \gamma_{x}\right)  \backslash\left\{  0\right\}  $.
		Hence, in both cases the implication%
		\begin{equation}
		\left(  \operatorname{Re}\left(  \bar{\mu}\left(  \mathfrak{R}_{j,k}\left(
		x\right)  u,u\right)  _{\gamma_{x}}\right)  =0\right)  \implies\left(  \left.
		u\right\vert _{\gamma_{x}}=0\right)  \label{Reunique1}%
		\end{equation}
		follows.
	\end{enumerate}
	
	\textbf{Second case: }Let $k\in\mathbb{C}_{<0}$ and the corresponding
	conditions in the Lemma be satisfied.
	
	In a similar way we obtain for $k\in\mathbb{C}_{<0}$ the estimate%
	\begin{align*}
		-\operatorname{Re}\left(  \bar{\mu}\left(  \mathfrak{R}_{j,k}\left(  x\right)
		u,u\right)  _{\gamma_x}\right)   &  \leq\frac{\operatorname{Re}k}{\left\vert
			k\right\vert }\operatorname{Im}\left(  \mathfrak{R}_{j,k}\left(  x\right)
		u,u\right)  _{\gamma_{x}}+\mu_{\operatorname*{R}}\left(  k\right)
		\frac{\operatorname{Im}k}{\left\vert k\right\vert }\left\Vert u\right\Vert
		_{\gamma_{x}}^{2}\\
		&  \leq\left(  \mu_{\operatorname*{R}}\left(  k\right)  \frac
		{\operatorname{Im}k}{\left\vert k\right\vert }\right)  \left\Vert u\right\Vert
		_{\gamma_{x}}^{2}.
	\end{align*}
	This time, the pre-factor $\left(  \mu_{\operatorname*{R}}\left(  k\right)  \frac
	{\operatorname{Im}k}{\left\vert k\right\vert }\right)  $ is negative yielding also
	in this case the implication%
	\begin{equation}
	\left(  \operatorname{Re}\left(  \bar{\mu}\left(  \mathfrak{R}_{j,k}\left(
	x\right)  u,u\right)  _{\gamma_x}\right)  =0\right)  \implies\left(  \left.
	u\right\vert _{\gamma_{x}}=0\right)  .\label{Reunique2}%
	\end{equation}

	We insert this into the left-hand side in~\eqref{uniqueness} and use (\ref{0ReCond})
	to get
	\begin{align*}
		0 &  =\left\vert \operatorname{Re}a_{j}\left(  u,\mu u\right)  \right\vert
		\geq\frac{\left\vert \operatorname{Im}k\right\vert }{\left\vert k\right\vert
		}\left(  \left\Vert \nabla u\right\Vert _{\Omega_{j}}^{2}+\left\vert
		k\right\vert ^{2}\left\Vert nu\right\Vert _{\Omega_{j}}^{2}\right)
		+\sum_{x:\left(  j,x\right)  \in\mathcal{J}%
		}\left\vert \operatorname{Re}\left(  \bar{\mu}\left(  \mathfrak{R}%
		_{j,k}\left(  x\right)  u,u\right)  _{\gamma_x}\right)  \right\vert \\
		&  \geq \sum_{x:\left(  j,x\right)  \in\mathcal{J}%
		}\left\vert \operatorname{Re}\left(  \bar{\mu}\left(  \mathfrak{R}%
		_{j,k}\left(  x\right)  u,u\right)  _{\gamma_x}\right)  \right\vert .
	\end{align*}
	Conditions (\ref{Reunique1}) and (\ref{Reunique2}) imply $\left.  u\right\vert
	_{\partial\Omega_{j}}=0$. The (homogeneous)  
	boundary conditions $\partial u/\partial \nu- \mathfrak{R}_{j,k}(\xi)(u)=0$ at
	$\gamma_\xi$ for the strong formulation of the left-hand side in~\eqref{uniqueness} yield
	$\left.  \partial u/\partial \nu\right\vert _{\gamma_\xi}=0$. This, in combination with the unique continuation principle (see~\cite{Alessandrini_ucp}) leads to $u=0$.
\end{proof}

\begin{example}
	\label{ExCasesII}\quad
	
	\begin{enumerate}
		\item[a.] If $\operatorname{Im}k\geq0$, the boundary conditions $\mathfrak{R}%
		_{j,k}\left(  x\right)  u=\operatorname*{i}k\hat{n}_j(x)u$ 
		satisfy (\ref{CondImkpos}) with
		$\mu_{\operatorname*{R}}=\left\vert \operatorname{Im}k\right\vert \hat{n}_j(x) $ and \\
		$\left(  \operatorname*{sign}\operatorname{Re}k\right)  \operatorname{Im}%
		\left(  \mathfrak{R}_{j,k}\left(  x\right)  u,u\right)  _{\gamma_{x}%
		}=\left\vert \operatorname{Re}k\right\vert \hat{n}_j(x) \left\Vert u\right\Vert
		_{\gamma_{x}}^{2}$ so that all conditions in (\ref{CondImkpos}) are satisfied.
		
		\item[b.] If $\operatorname{Im}k<0$, the boundary conditions $\mathfrak{R}%
		_{j,k}\left(  x\right)  u=-\operatorname*{i}k\hat{n}_j(x)u$ satisfy (\ref{CondImkneg})
		with $\mu_{\operatorname*{R}}=\left\vert \operatorname{Im}k\right\vert \hat{n}_j(x) $ and \\
		$-\left(  \operatorname*{sign}\operatorname{Re}k\right)  \operatorname{Im}%
		\left(  \mathfrak{R}_{j,k}\left(  x\right)  u,u\right)  _{\gamma_{x}%
		}=\left\vert \operatorname{Re}k\right\vert \hat{n}_j(x) \left\Vert u\right\Vert
		_{\gamma_{x}}^{2}$.
		
		\item[c.] If $k\in\mathbb{C}_{\geq0}^{\ast}$, the boundary condition
		$\mathfrak{R}_{2,k}\left(  1\right)  u=T_{k\hat{n}\left(  1\right)  }u$
		satisfy (\ref{CondImkpos}) with $\mu_{\operatorname*{R}}\left(  k\right)
		=c>0$ independent of $k$. This is proved in \cite{GrSa_DtN}.
	\end{enumerate}
\end{example}

The choices of boundary conditions as in Examples \ref{ExCasesI}.b and
\ref{ExCasesII}.b imply that the corresponding equations~\eqref{eq:f11k}-\eqref{eq:f22k} have unique solutions $f_{1,1,k}^{m}$, $f_{2,1,k}^{m}$,
$f_{2,2,k}^{m}$ for $k\in\mathbb{C}_{<0}$. According to Proposition
\ref{PropWellPosed} all possible resonances of problem~\eqref{eq:pdenep} belong to
$\mathbb{C}_{<0}$ and this choice of boundary conditions is justified. Any
solution $u$ in~\eqref{eq:pdenep} with Fourier coefficients $\hat{u}_{j,m}$ as in
(\ref{fourrierexp}) then can be written in the form
\begin{equation}%
\begin{array}
[c]{rl}%
\hat{u}_{1,m}= & A^m_{1,1}f_{1,1,k}^{m}\\
\hat{u}_{2,m}= & A^m_{2,1}f_{2,1,k}^{m}%
+A^m_{2,2}f_{2,2,k}^{m}%
\end{array}
\label{urep}%
\end{equation}
and the conditions for the coefficients $A^m_{1,1}$, $A^m_{2,1}$, $A^m_{2,2}$ are
such that $\hat{u}_{j,m}$ satisfy the two transmission conditions as well as
the $\operatorname*{DtN}$ boundary condition at $\gamma_{1}$:%
\begin{equation}
\left[
\begin{array}
[c]{ccc}%
f_{1,1,k}^{m}\left(  \xi\right)   & -f_{2,1,k}^{m}\left(  \xi\right)   &
-f_{2,2,k}^{m}\left(  \xi\right)  \\
\left(  f_{1,1,k}^{m}\right)  ^{\prime}\left(  \xi\right)   & -\left(
f_{2,1,k}^{m}\right)  ^{\prime}\left(  \xi\right)   & -\left(  f_{2,2,k}%
^{m}\right)  ^{\prime}\left(  \xi\right)  \\
0 & \left(  f_{2,1,k}^{m}\right)  ^{\prime}\left(  1\right)  -\beta_{2,k}%
^{m}\left(  1\right)  f_{2,1,k}^{m}\left(  1\right)   & \left(  f_{2,2,k}%
^{m}\right)  ^{\prime}\left(  1\right)  -\beta_{2,k}^{m}\left(  1\right)
f_{2,2,k}^{m}\left(  1\right)
\end{array}
\right]  \left(
\begin{array}
[c]{c}%
A^m_{1,1}\\
A^m_{2,1}\\
A^m_{2,2}%
\end{array}
\right)  =0\label{NLEVP}%
\end{equation}
with the Fourier multiplier corresponding to $\operatorname*{DtN}$ boundary
conditions%
\begin{equation}
\beta_{2,k}^{m}\left(  1\right)  =k\hat{n}\left(  1\right)  \frac{\left(
	H_{m}^{\left(  1\right)  }\right)  ^{\prime}\left(  k\hat{n}\left(  1\right)
	\right)  }{H_{m}^{\left(  1\right)  }\left(  k\hat{n}\left(  1\right)  \right)
}.\label{defbetaskm}%
\end{equation}
We denote the $3\times3$ matrix in (\ref{NLEVP}) by $T_{+}^{m}\left(
k\right)  $ and problem~\eqref{eq:pdenep} is equivalent to find resonances
$k\in\mathbb{C}_{<0}$ such that $\det T_{+}^{m}\left(  k\right)  =0$. In the
following we will simplify the problem by imposing the following assumption, to ensure that the resonances of the non-linear eigenvalue
problem~\eqref{eq:pdenep} are not simultaneously critical frequencies for the
auxiliary problems in (\ref{finaleqs}).

\begin{assumption}\label{assumption}
	Let $\left(  k,u\right)  \in\mathbb{C}_{<0}\times H_{0}^{1}\left(
	\Omega\right)\setminus\{0\}  $ be a resonance of problem~\eqref{eq:pdenep} and $m\in\mathbb{Z}$.
	Then, there exists a complex neighborhood $\omega_{k}\subset\mathbb{C}_{<0}$
	of $k$ such that the problems%
	\begin{equation}
	\left\{
	\begin{array}
	[c]{c}%
	L_{1,k}^{m}f_{1,1,k}^{m}=0\quad\text{in }\tau_{1}\\
	B_{1,k}^{m}\left(  0\right)  f_{1,1,k}^{m}=0\\
	\left(  f_{1,1,k}^{m}\right)  ^{\prime}\left(  \xi\right)  -\operatorname*{i}%
	k\hat{n}_1(\xi)f_{1,1,k}^{m}\left(  \xi\right)  =h_{1,1,k}^{m}\neq0
	\end{array}
	\right\}  \text{,\quad}\left\{
	\begin{array}
	[c]{c}%
	L_{2,k}^{m}f_{2,2,k}^{m}=0\quad\text{in }\tau_{2}\\
	-\left(  f_{2,2,k}^{m}\right)  ^{\prime}\left(  \xi\right)  -\operatorname*{i}%
	k\hat{n}_2(\xi)f_{2,2,k}^{m}\left(  \xi\right)  =h_{2,2,k}^m\neq0,\\
	\left(  f_{2,2,k}^{m}\right)  ^{\prime}\left(  1\right)  -\beta_{2,k}%
	^{m}\left(  1\right)  f_{2,2,k}^{m}\left(  1\right)  =0
	\end{array}
	\right\}  \label{finaleqs}%
	\end{equation}
	for $\beta_{2,k}^{m}(1)$ as in (\ref{defbetaskm}) have unique solutions.
\end{assumption}

This allows to reduce the non-linear
eigenvalue problem (\ref{NLEVP}) by using the ansatz
\[
  \hat{u}_{j,m}=A^m_{j}f_{j,j,k}^{m}\quad\forall
j\in\left\{  1,2\right\}
\]
and obtaining from the jump relations the condition%
\begin{equation}
\label{eq:nep}
T^m(k)A=
\left[
\begin{array}
[c]{cc}%
f_{1,1,k}^{m}\left(  \xi\right)   & -f_{2,2,k}^{m}\left(  \xi\right)  \\
\left(  f_{1,1,k}^{m}\right)  ^{\prime}\left(  \xi\right)   & -\left(
f_{2,2,k}^{m}\right)  ^{\prime}\left(  \xi\right)
\end{array}
\right]  \left(
\begin{array}
[c]{c}%
A^m_{1}\\
A^m_{2}%
\end{array}
\right)  =0,
\end{equation}
thus the link to an algebraic nonlinear eigenvalue problem (\ref{appendix:nep}~\footnote{See Supplementary Material.}). Going back to the system of ODEs~\eqref{eq:helmholtzpolarcordinates}, this is equivalent to finding $k$ such that there exists a nontrivial solution $\hat{u}_m$ of~\eqref{eq:helmholtzpolarcordinates} for $\hat{g}_m=0$, which is a nonlinear eigenvalue problem since the DtN boundary conditions depend in a nonlinear way on $k$.   In turn, this is equivalent to finding scattering resonances $k$ such that there exists a nontrivial solution $u$ of~\eqref{eq:helmholtz} with $g=0$;~\cite{araujo2018efficient}. In the spherical symmetric setting, we thus have reduced the nonlinear PDE eigenvalue problem~\eqref{eq:pdenep} to an algebraic nonlinear eigenvalue problem~\eqref{eq:nep}.

Thus the condition for critical complex $k$ (for which the equation~\eqref{eq:helmholtzpolarcordinates} is not well-posed) is given by
\begin{equation}\label{eq:determinant}
\mathrm{det}(k):=\det(T(k))=-[f_{1,1,k}(\xi) f'_{2,2,k}(\xi) - f'_{1,1,k}(\xi) f_{2,2,k}(\xi)]=0,
\end{equation}
where we skipped the superscript $m$.

\subsection{Quasi-resonances}

We will also be interested in quasi-resonances $\underline{k}=\operatorname{Re}\left( k\right)\ne0$ where $k$ is a resonance of problem~\eqref{eq:pdenep}. Since $\underline{k}$ is real, problem~\eqref{eq:helmholtz} is well posed.  Using similar reasoning as in the proof of Lemma~\ref{LemWellPosed} (see Examples \ref{ExCasesI}.a-
\ref{ExCasesII}.a and Examples \ref{ExCasesI}.c-\ref{ExCasesII}.c), we can show that the problems
\begin{align}
	&
	\begin{dcases}
		\begin{array}
			[c]{lc}%
			L^m_{1,\underline{k}}f^m_{1,1,\underline{k}}=0 & \text{in }\tau_{1}\\
			B^m_{1,\underline{k}}(0)f^m_{1,1,\underline{k}}=0 & \\
			\left(  f_{1,1,\underline{k}}^{m}\right)  ^{\prime}\left(  \xi\right)  -i\underline{k} \hat{n}_1(\xi) f_{1,1,\underline{k}}^{m}\left(  \xi\right)     =h_{1,1,\underline{k}}%
			^{m}\neq0  &
		\end{array}
		\label{eq:f11kquasi}
	\end{dcases}\\
	&
	\begin{dcases}
		\begin{array}
			[c]{lc}%
			L^m_{2,\underline{k}}f^m_{2,1,\underline{k}}=0 & \text{in }\tau_{2}\\
			-\left(  f_{2,1,\underline{k}}^{m}\right)  ^{\prime}\left(  \xi\right)  -i\underline{k}  \hat{n}_2(\xi)f_{2,1,\underline{k}}^{m}\left(  \xi\right)     =0\\
			\left(  f_{2,1,\underline{k}}^{m}\right)  ^{\prime}\left(  1\right)  -\beta_{2,\underline{k}}%
			^{m}\left(  1\right)  f_{2,1,\underline{k}}^{m}\left(  1\right)     =\hat{g}_m\phantom{\neq0}
		\end{array}
		\label{eq:f21kquasi}
	\end{dcases}
	\\
	&
	\begin{dcases}
		\begin{array}
			[c]{lc}%
			L^m_{2,\underline{k}}f^m_{2,2,\underline{k}}=0 & \text{in }\tau_{2}\\
			-\left(  f_{2,2,\underline{k}}^{m}\right)  ^{\prime}\left(  \xi\right)  -i\underline{k} \hat{n}_2(\xi) f_{2,2,\underline{k}}^{m}\left(  \xi\right)     =h_{2,2,\underline{k}}%
			^{m}\neq0\\
			\left(  f_{2,2,\underline{k}}^{m}\right)  ^{\prime}\left(  1\right)  -\beta_{2,\underline{k}}%
			^{m}\left(  1\right)  f_{2,2,\underline{k}}^{m}\left(  1\right)     =0
		\end{array}
		\label{eq:f22kquasi}%
	\end{dcases}
\end{align}
are well-posed for $\beta_{2,\underline{k}}^{m}(1)$ as in (\ref{defbetaskm}).

We thus can express the
solution of~\eqref{eq:helmholtzpolarcordinates} for a quasi-resonance $\underline{k}$ in terms of local homogeneous solutions via the ansatz%
\begin{equation}
\hat{u}_{1,m}=A^m_{1,1}f^m_{1,1,\underline{k}}\quad\text{and\quad}\hat{u}_{2,m}=f^m_{2,1,\underline{k}}%
+A^m_{2,2}f^m_{2,2,\underline{k}},
\label{eq:ansatzquasi}%
\end{equation}
and the transmission conditions lead to a linear system for the coefficients:
\begin{equation}\label{eq:linearsystemquasi}
\begin{pmatrix}
f^m_{1,1,\underline{k}}(\xi) & -f^m_{2,2,\underline{k}}(\xi) \\
(f^m_{1,1,\underline{k}})'(\xi) & -(f^m_{2,2,\underline{k}})'(\xi)
\end{pmatrix}
\begin{pmatrix}
A^m_{1,1}\\
A^m_{2,2}
\end{pmatrix}=
\begin{pmatrix}
f^m_{2,1,\underline{k}}(\xi)\\
(f^m_{2,1,\underline{k}})'(\xi)
\end{pmatrix}.
\end{equation}
In this way, the function~\eqref{eq:ansatzquasi} with $A^m_{1,1}$ and $A^m_{2,2}$ chosen as
in~\eqref{eq:linearsystemquasi} solves the boundary value problem~\eqref{eq:helmholtzpolarcordinates} for $k\leftarrow\underline{k}$.

For a resonance $k$ very close to the real axis (i.e. with small imaginary part), we expect the matrix of the system~\eqref{eq:linearsystemquasi} to be close to being singular and we refer to this case as a quasi-mode (see Section~\ref{quasi}).

\section{A Newton method for determining critical $k$'s}\label{sec:newton}

For a compact formulation, we suppress the superscript $m$ in the operators $L^m_{j,k}$, the functions $f^m_{j,\ell,k}$, and the data $\beta^m_{2,k}(1)$ and $h^m_{j,\ell,k}$. We want to determine the zeros of $\mathrm{det}(k)$~\eqref{eq:determinant} 
in $\CC^*$ knowing only approximations of $f_{1,1,k}$ and $f_{2,2,k}$ given by an ODE solver (see Section~\ref{sec:odesolver} for details on the ODE solver). For this purpose, we formulate a Newton method. Since Newton's method requires the derivative of $\mathrm{det}(k)$ with respect to $k$, we have to evaluate $\partial_k f_{1,1,k}(\xi)$ and $\partial_k f_{2,2,k}(\xi)$ along $\partial_k f'_{1,1,k}(\xi)$ and $\partial_k f'_{2,2,k}(\xi)$. In order to avoid finite-difference approximation, we are going to establish the equations satisfied by $\partial_k f_{1,1,k}$ and $\partial_k f_{2,2,k}$ (see Proposition~\ref{prop:derivatives}). Also, the choice of boundary conditions at $r=\xi$ on $f_{1,1,k}$~\eqref{eq:f11k} and $f_{2,2,k}$~\eqref{eq:f22k} allows us to write
\begin{equation}
\begin{dcases}
f'_{1,1,k}(\xi)=ikn_1(\xi) f_{1,1,k}(\xi)+h_{1,1,k}\\
f'_{2,2,k}(\xi)=-ikn_2(\xi) f_{2,2,k}(\xi)- h_{2,2,k} 
\end{dcases}
\end{equation}
without additional approximation. For the sake of simplicity, and where there is no ambiguity, we will use the notations $n_1$ and $n_2$ instead of $\hat{n}_1$ and $\hat{n}_2$.

The expressions of the determinant and its derivative are then given by
\begin{align}\label{eq:determinant2}
	\mathrm{det}(k)&=-[f_{1,1,k}(\xi) (-i k n_2(\xi) f_{2,2,k}(\xi)-h_{2,2,k}) - (i k n_1(\xi)f_{1,1,k}(\xi)+h_{1,1,k}) f_{2,2,k}(\xi)]\notag\\
	&=i k (n_2(\xi)+n_1(\xi)) f_{1,1,k}(\xi) f_{2,2,k}(\xi) + h_{2,2,k}f_{1,1,k}(\xi)+h_{1,1,k}f_{2,2,k}(\xi)
\end{align}
and
\begin{align}\label{eq:derdeterminant2}
	\partial_k\mathrm{det}(k)=&i(n_2(\xi)+n_1(\xi))(f_{1,1,k}(\xi)f_{2,2,k}(\xi)+k\partial_k f_{1,1,k}(\xi) f_{2,2,k}(\xi)+k f_{1,1,k}(\xi)\partial_k f_{2,2,k}(\xi)) \notag\\ 
	&+h_{2,2,k}\partial_k f_{1,1,k}(\xi)+ \partial_k h_{2,2,k} f_{1,1,k}(\xi)+h_{1,1,k}\partial_k f_{2,2,k}(\xi)+\partial_kh_{1,1,k} f_{2,2,k}(\xi).
\end{align}

\begin{proposition}\label{prop:derivatives}
	Let $(\tilde{k},\tilde{u})\in\CC_{<0}\times H^1(\Omega)\setminus\{0\}$  be a resonance pair of~\eqref{eq:pdenep} and Assumption~\ref{assumption} be satisfied for some ${\omega_{\tilde{k}}} \subset \CC_{<0}$. The equations satisfied by $\partial_k f_{1,1,k}$ and $\partial_k f_{2,2,k}$ are the following:
	\begin{equation}\label{eq:dkf11k}
	\begin{dcases}
	L_{1,k} (\partial_kf_{1,1,k})=2 k n_1^2 f_{1,1,k} & \text{in $\tau_1$}\\
	B_{1,k}(0) (\partial_kf_{1,1,k})=0  & \text{Dirichlet or Neumann boundary conditions} \\
	(\partial_kf_{1,1,k})'(\xi)-ik n_1(\xi) \partial_kf_{1,1,k}(\xi)= \partial_k h_{1,1,k} +i n_1(\xi) f_{1,1,k}(\xi) & \text{Robin boundary conditions}
	\end{dcases}
	\end{equation}
	\begin{equation}\label{eq:dkf22k}
	\begin{dcases}
	L_{2,k} (\partial_k f_{2,2,k})=2 k n_2^2 f_{2,2,k} & \text{in $\tau_2$}\\
	-(\partial_k f_{2,2,k})'(\xi)-ik n_2(\xi)\partial_k f_{2,2,k}(\xi) = \partial_k h_{2,2,k}+i n_2(\xi) f_{2,2,k}(\xi) & \text{Robin boundary conditions}\\
	(\partial_k f_{2,2,k})'(1)-\beta_{2,k}(1) \partial_k f_{2,2,k}(1)=(\partial_k \beta_{2,k}(1)) f_{2,2,k}(1)& \text{DtN boundary conditions}  \phantom{............................}
	\end{dcases}
	\end{equation}
	where $\beta_{2,k}(1)=kn_2(1)\frac{H'_{m}(kn_2(1))}{H_{m}(kn_2(1))}$. 
	
	These equations are well-posed in $V=H^1(\tau_1)$ (even $m$) or $V=\{v\in H^1(\tau_1), v(0)=0\}$ (odd $m$), and $H^1(\tau_2)$ respectively, for all $k\in \omega_{\tilde{k}}$.
\end{proposition}

\begin{proof}
	
	These results can be established by switching back from spherical to Cartesian coordinates, and passing to the limit $\varepsilon\to 0$ ($\varepsilon\in \CC^*$) in the equation satisfied by $\frac{\tilde{f}_{k+\varepsilon}-\tilde{f}_k}{\varepsilon}$, where for instance $\tilde{f}_k(x)=f_{1,1,k}(r)e^{im\theta}$, see~\cite[\S 3]{BBSS}.

\end{proof}

\subsection{Newton's algorithm}\label{sec:newtonalgo}

We use Newton's method for complex differentiable functions:

\begin{algorithm}[H]
	\caption{Newton's algorithm}
	\begin{algorithmic}
		\STATE $\mathrm{Newton} (k_0,\varepsilon,l_{\max})$
		\FOR{$l=0\ldots$ until stopping criterion is reached}
		\STATE $\bullet$ Start with $k=k_l;$
		\STATE $\bullet$ Solve two two-point BVPs~\eqref{finaleqs} exactly or numerically for $k=k_l$: result: approximations $f_{1,1,k}$ and $f_{2,2,k}$ to~\eqref{finaleqs};
		\STATE $\bullet$ Compute the derivatives with respect to $k$, $\partial_k f_{1,1,k}$ and $\partial_k f_{2,2,k}$, as the solutions of two two-point BVPs~\eqref{eq:dkf11k}-\eqref{eq:dkf22k}; 
		\STATE $\bullet$ Compute $\mathrm{det}(k_l)$ by using~\eqref{eq:determinant2} and $\partial_k\mathrm{det}(k_l)$ by~\eqref{eq:derdeterminant2};
		\STATE $\bullet$ Compute
		\begin{equation}\label{eq:newton}
		k_{l+1}=k_l-\frac{\mathrm{det}(k_l)}{\partial_k\mathrm{det}(k_l)};
		\end{equation}
		\STATE $\bullet$ $l\leftarrow l+1;$
		\STATE The stopping criterion is given by $\frac{|\mathrm{det}(k_l)|}{\|T(k_l)\|_F}\le \varepsilon$ or when a maximal number of iterations $l_{\max}$ is reached.
		\ENDFOR
	\end{algorithmic}
\end{algorithm}
We use a stopping criterion based on the relative residual $\frac{|\mathrm{det}(k_l)|}{\|T(k_l)\|_F}$, where $\|\cdot\|_F$ is the Frobenius norm, as in~\cite{guttel2017nonlinear}.

\begin{remark}
	We used the direct computation of the determinant since $T$ is a $2\times2$ matrix, but for larger matrices (as in the case of multiple jump points), it would be more pertinent to use the following formula, holding for $k$ such that $\det(k)\ne0$,
	\begin{equation}\label{eq:traceformula}
	\partial_k\mathrm{det}(k)=\mathrm{det}(k)\, \mathrm{trace}({T^{-1}(k)\partial_k T(k)}),
	\end{equation}
	leading to the \emph{Newton-trace} iteration (see \cite[Eq.~(4.2)--(4.3)]{guttel2017nonlinear})
	\begin{equation}
	k_{l+1}=k_l-\frac{1}{\mathrm{trace}({T^{-1}(k_l)\partial_k T(k_l)})}.
	\end{equation}
\end{remark}

\subsection{Starting values}\label{sec:startingvalues}

Since we are interested in whispering gallery modes (WGM) which are associated with complex scattering resonances very close to the real axis, we will consider these different choices for the starting values of $k$ in the Newton method:

\begin{itemize}
	\item Trust region strategies~\cite{voss2013nonlinear,yang2007trust}:
	\begin{itemize}
		\item Start from the real axis
		
		\item Start from asymptotic expansions: in~\cite{balac2021asymptotics, moitier2019etude}, asymptotic expansion of WGM resonances when $m\to\infty$ for cavities with radially varying optical index are proposed in some configurations. Even if our problem may differ from the one of~\cite{balac2021asymptotics}, their asymptotics could represent a good first guess of the critical states in our case. For instance, considering the leading term in the asymptotics (\ref{appendix:dauge}), a possible starting point is
		\begin{equation}\label{eq:startingvalue}
		k=\frac{m}{\xi n_0} \quad\text{where}\quad  n_0=\lim\limits_{r\nearrow\xi} n(r).
		\end{equation}
	\end{itemize}
	\item Homotopy methods~\cite{daya2001numerical}:
	\begin{itemize}
		\item If Newton's method does not converge for a given refractive index and a given starting value $k_0$, it could be useful to first solve the problem for a piecewise constant refractive index case or an intermediate case close to $n$ for which the algorithm converges, and then take the resulting $k$ as a starting value for the initial problem. 
	\end{itemize}
\end{itemize}

\subsection{Local convergence}

We state here an important result of this paper, which deals with the local convergence of the suggested complex Newton method in Section~\ref{sec:newtonalgo}, under suitable conditions. We assume here that $\det(k)$~\eqref{eq:determinant2} and $\partial_k \det(k)$~\eqref{eq:derdeterminant2} are evaluated using the exact 
fundamental system for some $k$.

\begin{theorem}\label{th:newton}
	Let $n\colon (0,1)\to \RR_+^*$ be piecewise smooth (on $(0,\xi)$ and $(\xi,1))$ and bounded from above and below by positive numbers. Let $k_\infty\in\CC_{<0}$  be a resonance of~\eqref{eq:pdenep} and Assumption~\ref{assumption} be satisfied for some ${\omega_{k_\infty}} \subset \CC_{<0}$. Then $\det\colon {\omega_{k_\infty}}\to\CC$ is analytic. If $k_\infty\neq 0$ is a simple root of $\det$ then there exists $\delta>0$ such that for any starting point $k_0\in B(k_\infty,\delta)$, the sequence
	\begin{align*}
		k_{l+1}=k_l-\frac{\det(k_l)}{\partial_k \det(k_l)}
	\end{align*}
	converges to $k_\infty$ and the convergence is quadratic.
\end{theorem}

\begin{proof}
	\textbf{Step 1:} Analyticity of $\det$
	
	Under the assumptions on $n$, we know from Assumption~\ref{assumption} that the solutions $f_{1,1,k}$ and $f_{2,2,k}$ exist and are unique in $H^1(0,\xi)$ and $H^1(\xi,1)$ respectively for all $k\in{\omega_{k_\infty}}$. Thus $\det$ in~\eqref{eq:determinant2} is well defined as a function from ${\omega_{k_\infty}}$ to $\CC$.
	
	We also established in Proposition~\ref{prop:derivatives} that $\partial_k f_{1,1,k}$ and $\partial_k f_{2,2,k}$ exist and are unique in $H^1(0,\xi)$ and $H^1(\xi,1)$ respectively for all $k\in{\omega_{k_\infty}}$. Thus $\det$ is complex-differentiable, and hence analytic in ${\omega_{k_\infty}}$.
	
	\textbf{Step 2:} Convergence of the complex Newton method
	
	Since $\det$ is analytic, the local quadratic convergence of the complex Newton method to a simple root starting sufficiently close to the root follows classically~\cite{ostrowski1973solution,alexander2012early,wilkinson1988algebraic}. 
\end{proof}

\section{Simplicity of the roots}\label{sec:simplicity}

In order to apply Theorem~\ref{th:newton} providing the local convergence of the Newton algorithm, it remains to show that the roots of $\det$~\eqref{eq:determinant2} are simple.

\subsection{Piecewise constant case}\label{sec:pwccase}

In the piecewise constant case, we have explicit expressions of $f_{1,1,k}$ and $f_{2,2,k}$ in terms of Bessel and Hankel functions, leading to an explicit expression of the determinant~\eqref{eq:determinant}:
\begin{equation}\label{eq:determinantpwconstant}
\mathrm{det}(k)/k=-[J_m(kn_1\xi) n_2 H_m'(k n_2 \xi)- n_1 J_m'(k n_1\xi) H_m (kn_2 \xi)]=D(k),    
\end{equation}
which is exactly Eq.(1.6) in~\cite{balac2021asymptotics} with $n_2=1$, $n_1=n_0$ and $\xi=R$ (see also~\cite{moitier2019etude,moiola2019acoustic}). 
By forming $K=kn_2\xi$ and $N=\frac{n_1}{n_2}$ in~\eqref{eq:determinantpwconstant}, 
\begin{equation}\label{eq:wlog}
D(k)=Nn_2J_m'(NK)H_m(K)-n_2 J_m(NK) H_m'(K)=n_2\tilde{D}_N(K),
\end{equation}
where
\begin{equation}\label{eq:Dtilde}
\tilde{D}_N(K)=NJ_m'(NK)H_m(K)- J_m(NK) H_m'(K).
\end{equation}
Thus if $n_1$ and $n_2$ are given such that $N>1$ then we can recover the zeros $k_0$ of $D$ for a given $\xi$ from the zeros $K_0$ of $\tilde{D}_N$ using $k_0=\frac{K_0}{n_2\xi}$. Thus we are interested in the zeros of the following function of $z$ (by replacing $N\leftarrow n$ and $K\leftarrow z$ in~\eqref{eq:Dtilde}):
\begin{equation}\label{eq:D(k)}
D(z)=n J_m'(nz) H_m(z) - J_m(n z) H_m'(z)
\end{equation}
whose derivative with respect to $z$ is given by
\begin{equation}\label{eq:DerD(k)}
D'(z)=n^2 J_m''(n z ) H_m(z) - J_m(n z) H_m''(z). 
\end{equation}

\begin{theorem}\label{th:simplicity}
	Let $D$ be as in~\eqref{eq:D(k)} and $n>1$. For $m\ge 0$, all roots of $D$ except, possibly, $z=0$ are simple.
\end{theorem}

\begin{proof}
	We start by simplifying the expression of the derivative $D'$~\eqref{eq:DerD(k)} to remove the second derivatives. We know that $J_m$ and $H_m$, for $m\in\RR$, are solutions of the Bessel differential equation
	\begin{equation}\label{eq:besselEq}
	z^2 y''(z)+ zy'(z)+(z^2-m^2) y(z)=0,
	\end{equation}
	i.e., for $z\ne0$,
	\begin{equation}\label{eq:besselEq2}
	y''(z)= -\frac{y'(z)}z-\Bigl(1-\frac{m^2}{z^2}\Bigr) y(z).
	\end{equation}
	Thus
	\begin{align}\label{eq:derDEq}
		D'(z)&=n^2\Bigl[-\frac{J_m'(nz)}{nz}-\Bigl(1-\frac{m^2}{(nz)^2}\Bigr) J_m(nz)\Bigr] H_m(z)-J_m(nz) \Bigl[-\frac{H_m'(z)}{z}-\Bigl(1-\frac{m^2}{z^2}\Bigr) H_m(z)\Bigr]\notag\\
		&=-n\frac{J_m'(nz)}{z}H_m(z)-\Bigl(n^2-\frac{m^2}{z^2}\Bigr) J_m(nz) H_m(z)+J_m(nz)\frac{H_m'(z)}{z}+\Bigl(1-\frac{m^2}{z^2}\Bigr) J_m(nz)H_m(z)\notag\\
		&=-\frac{1}{z}\Bigl[n J_m'(nz)H_m(z)-J_m(nz) H_m'(z) \Bigr]-(n^2-1) J_m(nz) H_m(z)\notag\\
		&=-\frac{D(z)}z-(n^2-1)J_m(nz) H_m(z).
	\end{align}
	
	If $z_0\ne 0$ is a zero of $D(z)$ with multiplicity $\ge2$, then we have $D(z_0)=D'(z_0)=0$. Since $n>1$ and $z_0\ne 0$, it follows from the differential equation~\eqref{eq:derDEq} that $z_0$ is also a zero of $J_m(nz) H_m(z)$. Since, for $m\ge -1$, the Bessel function $J_m(z)$ has only real zeros~\cite{abramowitz1968handbook}\cite[\S 15.25]{watson1922treatise}, and $D(z)\ne 0$ for all $z\in\RR\setminus\{0\}$ (recall that resonances $z\ne 0$ must satisfy $\Im(z)<0$), we obtain that $z_0$ is a zero of $H_m(z)$. Using $D(z_0)=H_m(z_0)=0$ in~\eqref{eq:D(k)}, we get also $H_m'(z_0)=0$ (since $J_m(nz_0)\ne0$ for the same reasons stated earlier). Thus $z_0$ is a multiple zero of $H_m$. But all zeros of $H_m$, except $z=0$, are simple~\cite[\S 15.21]{watson1922treatise}~\cite{sturm2009memoire}. This is a contradiction.    
\end{proof}

\subsection{Piecewise smooth case}

We prove here the geometric simplicity of eigenvalues (see definition in~\cite{voss2013nonlinear} and \ref{appendix:nep}) that is necessary for the roots to be simple. The study of the algebraic simplicity of eigenvalues (which is equivalent to the simplicity of the roots) is reported in~\cite[\S5.1.2]{BBSS}.

\begin{proposition}\label{prop:geometricsimplicity} 
	Let $\lambda\in\CC^*$ be an eigenvalue of the nonlinear eigenvalue problem~\eqref{eq:nep} (under Assumption~\ref{assumption}). Then $\lambda$ is geometrically simple.
\end{proposition}

\begin{proof}
	
	Let $\begin{pmatrix}x_1\\x_2\end{pmatrix}\in\CC^2\setminus\{0\}$ be an eigenvector corresponding to the eigenvalue $\lambda$, i.e., 
	\begin{equation}\label{eq:detsimplicity}
	f'_{1,1,\lambda}(\xi) f_{2,2,\lambda}(\xi)-f_{1,1,\lambda}(\xi) f'_{2,2,\lambda}(\xi)=0
	\end{equation}
	and
	\begin{equation}\label{eq:geomsimplicity}
	\begin{dcases}
	f_{1,1,\lambda}(\xi) x_1 - f_{2,2,\lambda}(\xi) x_2 =0\\
	f'_{1,1,\lambda}(\xi) x_1 - f'_{2,2,\lambda}(\xi) x_2 =0.
	\end{dcases}
	\end{equation}
	From $f'_{1,1,\lambda}(\xi)-i\lambda n_1(\xi)f_{1,1,\lambda}(\xi)=h_{1,1,\lambda}\ne 0$~\eqref{eq:f11k} (resp. $-f'_{2,2,\lambda}(\xi)-i\lambda n_2(\xi)f_{2,2,\lambda}(\xi)=h_{2,2,\lambda}\ne 0$~\eqref{eq:f22k}), we get that $f_{1,1,\lambda}(\xi)$ and $f'_{1,1,\lambda}(\xi)$ cannot vanish simultaneously (resp. $f_{2,2,\lambda}(\xi)$ and $f'_{2,2,\lambda}(\xi)$ cannot vanish simultaneously). Moreover, from~\eqref{eq:detsimplicity}, $f_{1,1,\lambda}(\xi)=0$ if and only if $f_{2,2,\lambda}(\xi)=0$ (resp. $f'_{1,1,\lambda}(\xi)=0$ if and only if $f'_{2,2,\lambda}(\xi)=0$). By using this in~\eqref{eq:geomsimplicity} we conclude that, there exists a nonzero complex constant $C_\lambda\in\Bigl\{\frac{f_{2,2,\lambda}(\xi)}{f_{1,1,\lambda}(\xi)},\frac{f'_{2,2,\lambda}(\xi)}{f'_{1,1,\lambda}(\xi)}\Bigr\}$
	such that
	\begin{equation}
	x_1=C_\lambda x_2,
	\end{equation}
	and any eigenvector corresponding to the same eigenvalue $\lambda$ can be written $x_2\begin{pmatrix}C_\lambda\\1\end{pmatrix}$ hence ${\dim(\ker(T(\lambda)))=1}$, i.e., the eigenvalue $\lambda$ is geometrically simple.
\end{proof}
\begin{remark}
	One could also adapt the proof of the simplicity of eigenvalues for Sturm-Liouville problems~\cite[Thm.~11.2.3]{boyce2012elementary}~\cite[Thm.~2.4]{pryce1993numerical} to our case, and show the simplicity of eigenvalues for the ODE~\eqref{eq:helmholtzpolarcordinates} with $\hat{g}_m=0$ and then deduce the geometric simplicity of the eigenvalues for the corresponding nonlinear eigenvalue problem.
\end{remark}

\section{A proper conditioning}\label{sec:scaling}

The functions $f_{1,1,k}$ and $f_{2,2,k}$ in~\eqref{finaleqs} depend on the choice
of $h_{1,1,k}$, $h_{2,2,k}$ and have impact to the performance of the Newton
method. Our concrete choice of these control parameters is motivated by the
following consideration for piecewise constant refractive index. In this case,
the unique solutions $f_{1,1,k}$ and $f_{2,2,k}$ of~\eqref{finaleqs} (under Assumption~\ref{assumption}) are given by
$r\mapsto J_m(kn_1r)$ and $r\mapsto H_m(k n_2 r)$ respectively, for a specific choice of $h_{1,1,k}$ and $h_{2,2,k}$ namely
\begin{equation}\label{eq:databoundarychoice}
\begin{dcases}
h_{1,1,k}=k n_1(\xi) J'_m(k \xi n_1(\xi))-ik n_1(\xi) J_m(k\xi n_1(\xi))=k n_1(\xi)[J'_m(k \xi n_1(\xi))-iJ_m(k\xi n_1(\xi))]=:\frac1{c_{1,k}}\\
h_{2,2,k}=-k n_2(\xi) H'_m(k \xi n_2(\xi))-ikn_2(\xi) H_m(k \xi n_2(\xi))=-k n_2(\xi) [H'_m(k \xi n_2(\xi))+iH_m(k \xi n_2(\xi))]=:\frac1{c_{2,k}}.
\end{dcases}
\end{equation}

\begin{lemma}\label{lem:positivityc1kc2k}
	For all $k\in\mathbb{C}_{\geq0}^{\ast}$ and all $\xi\in{]}0,1[$, $J_m'(k\xi n_1(\xi))-i J_m(k\xi n_1(\xi))\ne 0$ and $H_m'(k\xi n_2(\xi))+i H_m(k\xi n_2(\xi))\ne 0$.
\end{lemma}
\begin{proof}
	Let $k\in\mathbb{C}_{\geq0}^{\ast}$. Following Lemma~\ref{LemWellPosed} (Examples \ref{ExCasesI}.a-\ref{ExCasesII}.a), the function $r\mapsto J_m(kn_1(\xi)r)$ is the unique solution of the differential equation
	\begin{equation}
	\begin{dcases}
	-y''(r)-\frac1r y'(r)+\Bigl(\frac{m^2}{r^2}-k^2n_1^2(\xi)\Bigr)y(r)=0 & \text{in $(0,\xi)$}\\
	y'(\xi)-ikn_1(\xi) y(\xi)=kn_1(\xi) J'_m(k\xi n_1(\xi))-ik n_1(\xi) J_m(k\xi n_1(\xi)) & \text{Robin boundary conditions}\\
	y(0)=0 & \text{for odd $m$}\\
	y'(0)=0& \text{for even $m$}.
	\end{dcases}
	\end{equation}
	If there exists $k\in\mathbb{C}_{\geq0}^{\ast}$ such that 
	$J'_m(k\xi n_1(\xi))-iJ_m(k\xi n_1(\xi))=0$ then  the solution $r\mapsto J_m(kn_1(\xi)r)$ is identically zero. This is a contradiction.\footnote{Note that this reasoning (and adapted one to $\mathbb{C}_{<0}$) also allows us to obtain another proof of the simplicity of the zeros of $J_m$ that was used in the proof of Theorem~\ref{th:simplicity}.}  
	For the result concerning $H_m$, we use a similar reasoning considering the unique solution of the differential equation
	\begin{equation}
	\begin{dcases}
	-y''(r)-\frac1r y'(r)+\Bigl(\frac{m^2}{r^2}-k^2n_2^2(\xi)\Bigr)y(r)=0 & \text{in $(\xi,1)$}\\
	-y'(\xi)-ikn_2(\xi) y(\xi)=-kn_2(\xi) H'_m(k\xi n_2(\xi))-ikn_2(\xi) H_m(k\xi n_2(\xi)) & \text{Robin boundary conditions}\\
	y'(1)-kn_2(\xi) \frac{H'_{m}(kn_2(\xi))}{H_{m}(kn_2(\xi))} y(1)=0 & \text{DtN boundary conditions}
	\end{dcases}
	\end{equation}
	since $H_m(kn_2(\xi))\ne 0$ for $k\in\mathbb{C}_{\geq0}^{\ast}$ (for integer orders $m$, the zeros of the Hankel function of the first kind 
	always have negative imaginary parts~\cite{abramowitz1968handbook}).%
\end{proof}

\begin{lemma}\label{lem:positivityc1kc2kadditional}
	Let Assumption~\ref{assumption} be satisfied in the piecewise constant case, for a resonance $(\tilde{k},\tilde{u})\in \CC_{<0}\times H^1(\Omega)\setminus\{0\}$ of problem~\ref{eq:pdenep} and neighborhood $\omega_{\tilde{k}}$ as in Assumption~\ref{assumption}. Then, for any $k\in\omega_{\tilde{k}}$, $J_m'(k\xi n_1(\xi))-i J_m(k\xi n_1(\xi))\ne 0$ and $H_m'(k\xi n_2(\xi))+i H_m(k\xi n_2(\xi))\ne 0$.
\end{lemma}
\begin{proof}
	Under Assumption~\ref{assumption}, problems~\eqref{finaleqs} are well posed for all $k\in\omega_{\tilde{k}}$. Using the same reasoning as in the proof of Lemma~\ref{lem:positivityc1kc2k}, it follows that for any $k\in\omega_{\tilde{k}}$, $J_m'(k\xi n_1(\xi))-i J_m(k\xi n_1(\xi))\ne 0$ and $H_m'(k\xi n_2(\xi))+i H_m(k\xi n_2(\xi))\ne 0$, since the contrary leads to a contradiction.
\end{proof}

\begin{remark}\label{rem:zerosc1kc2k}
	It follows that, in the piecewise constant case, the resonances of the auxiliary problems in~\eqref{finaleqs} are respectively the zeros of $\frac1{c_{1,k}}$ and the zeros of $\frac1{c_{2,k}}$~\eqref{eq:databoundarychoice}.
\end{remark}

For the simplest choice $h_{1,1,k}=h_{2,2,k}=1$, it holds
\begin{equation}\label{eq:databoundarychoice2}
\begin{dcases}
f_{1,1,k}(r)=c_{1,k} J_m(k n_1 r)\\
f_{2,2,k}(r)=c_{2,k} H_m(k n_2 r),
\end{dcases}
\end{equation}
when the fundamental system~\eqref{finaleqs} is well posed, 
and the Newton algorithm on this choice (leading to ``$\mathrm{det}_2$'') could have a different behavior than the one (``$\mathrm{det}_1$'') in~\eqref{eq:determinantpwconstant}  as, even if
\begin{equation}\label{eq:relationship}
\mathrm{det}_2(k)=c_{1,k} c_{2,k} \mathrm{det}_1(k)=\alpha(k) \mathrm{det}_1(k) 
\end{equation}
for some function $\alpha$ satisfying $\alpha(k)\ne 0$ for all $k\in\omega_{\tilde{k}}$ (Lemma~\ref{lem:positivityc1kc2kadditional})\footnote{Note that the relation~\eqref{eq:relationship} holds when both $\det_1$ and $\det_2$ are defined, i.e. in the piecewise constant case for $k\in\CC^*_{\ge0}\cup (\CC_{<0}\setminus Z)$ where $Z$ is the set of the zeros of $\frac{1}{c_{1,k}}$ or $\frac{1}{c_{2,k}}$ (Remark~\ref{rem:zerosc1kc2k}).}, one has
\begin{align}
	\frac{\partial_k \mathrm{det}_2(k)}{\mathrm{det}_2(k)}=\frac{\partial_k c_{1,k}}{c_{1,k}}+\frac{\partial_k c_{2,k}}{ c_{2,k}}+ \frac{\partial_k \mathrm{det}_1(k)}{\mathrm{det}_1(k)}
	\neq \frac{\partial_k \mathrm{det}_1(k)}{\mathrm{det}_1(k)},
\end{align}
so that different choices of $h_{1,1,k}$ and $h_{2,2,k}$ lead to different scalings of the fundamental solutions and hence to different scalings of the determinant, which in turn may affect the Newton algorithm.

\begin{figure}[h!]
	\centering
	\begin{subfigure}[b]{0.35\textwidth}
		\centering
		\includegraphics[width=\textwidth]{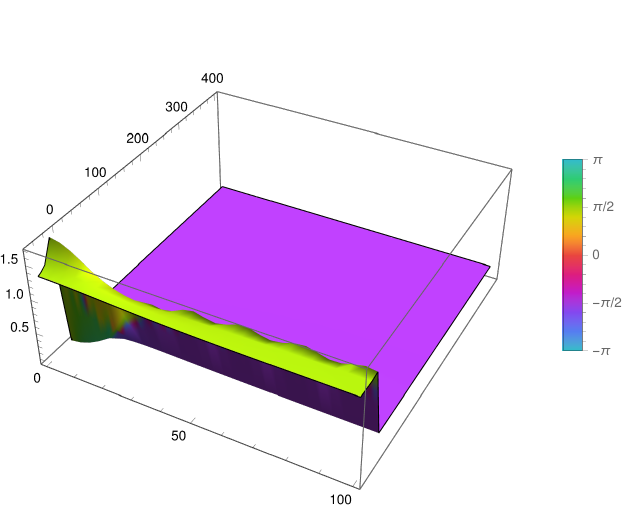}
		\caption{$\frac{\partial_k \mathrm{det}_1(k)}{\mathrm{det}_1(k)}$}
		\label{fig:ratiodet1}
	\end{subfigure}
	\begin{subfigure}[b]{0.35\textwidth}
		\centering
		\includegraphics[width=\textwidth]{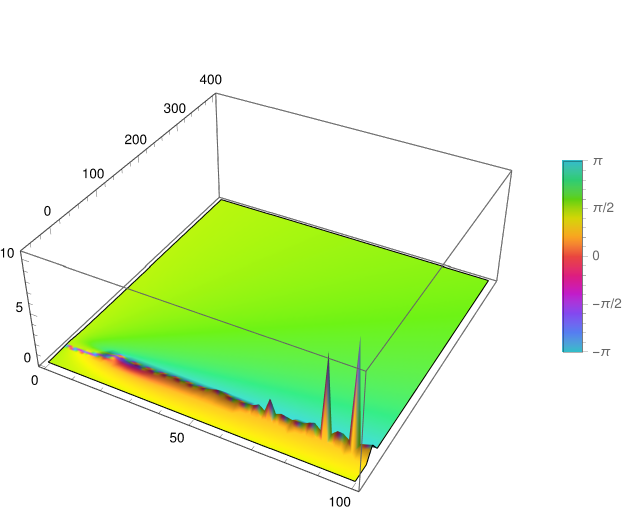}
		\caption{$\frac{\partial_k \mathrm{det}_2(k)}{\mathrm{det}_2(k)}$}
		\label{fig:ratiodet2}
	\end{subfigure}
	
	\caption{Comparison between $\frac{\partial_k \mathrm{det}_2(k)}{\mathrm{det}_2(k)}$ and $\frac{\partial_k \mathrm{det}_1(k)}{\mathrm{det}_1(k)}$ for $\xi=0.5$, $n_1=1.5$, $n_2=1$ and $m=10$. Here $k$ is varying over the complex rectangle with vertices $-50i$ and $100+400i$.}
	\label{fig:c1kc2k}
\end{figure}

Indeed, in our experiments, the use of $\mathrm{det_2}$ may  cause some roundoff problems as the ratio  $\frac{\partial_k \mathrm{det}_2(k)}{\mathrm{det}_2(k)}$ may be close to zero in contrast to $\frac{\partial_k \mathrm{det}_1(k)}{\mathrm{det}_1(k)}$, see Figure~\ref{fig:c1kc2k}\footnote{The comparison should be considered only for $k$ outside the set $Z$ since the plotted $\det_1$ is the exact one~\eqref{eq:determinantpwconstant}, while the constructed $\det_1$ obtained solving $f_{1,1,k}$ and $f_{2,2,k}$ with $h_{1,1,k}=\frac1{c_{1,k}}$ and $h_{2,2,k}=\frac1{c_{2,k}}$ is also undefined for $k\in Z$.}.

\begin{figure}[h!]
	\centering
	\includegraphics[scale=0.35]{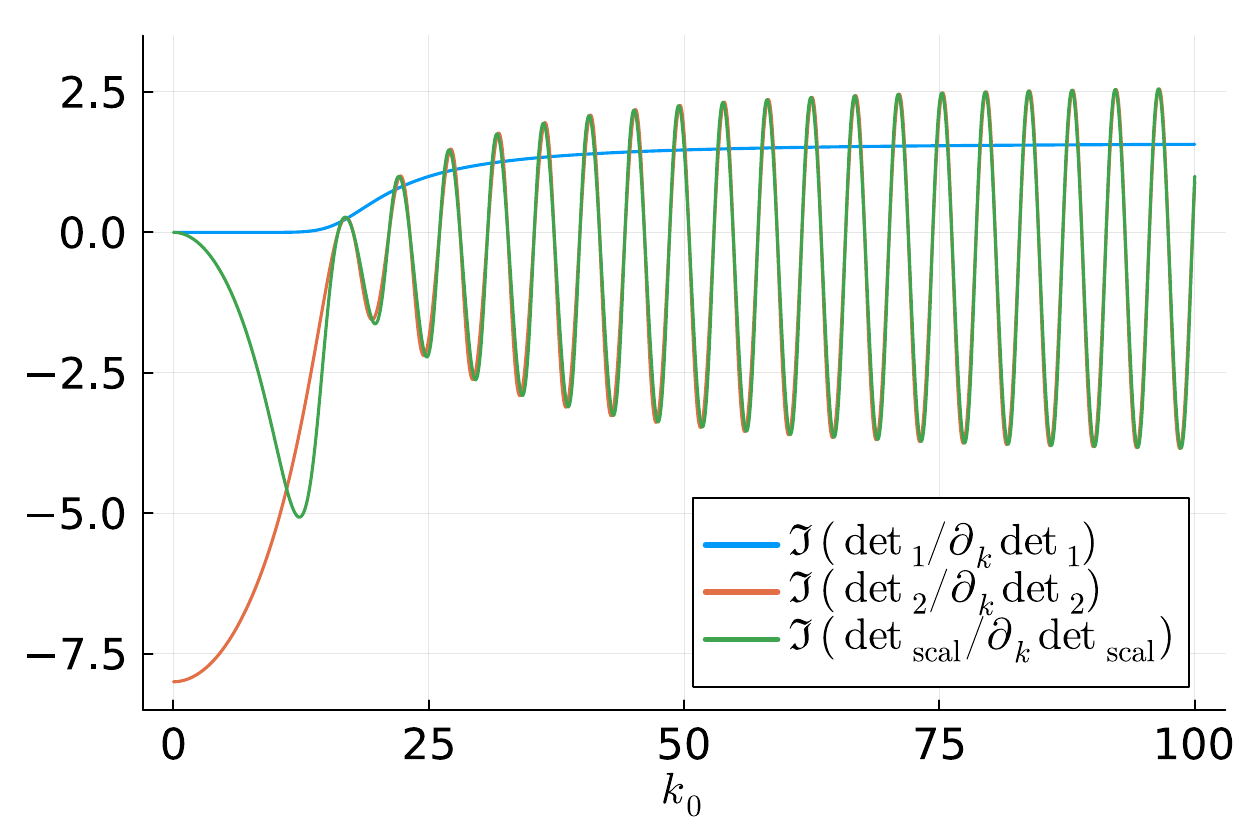}
	\caption{Plot of $\Im\bigl(\frac{\det(k_0)}{\partial_k \mathrm{det}(k_0)}\bigr)$ on the real axis for $\xi=0.5$, $n_1=1.5$, $n_2=1$ and $m=10$.}
	\label{fig:imaginarypart}
\end{figure}

Another observation regarding the use of the Newton method with $\mathrm{det}_2$ starting from the real axis is that, often, the first iterate has in our experiments a positive imaginary part which is not optimal knowing that the roots satisfy $\Im(k)< 0$.  This is related to the sign of $\Im\bigl(\frac{\det(k_0)}{\partial_k \mathrm{det}(k_0)}\bigr)$: if it is positive then the first iterate has a negative imaginary part, otherwise it has a positive imaginary part. We plot this quantity $\Im\bigl(\frac{\det(k_0)}{\partial_k \mathrm{det}(k_0)}\bigr)$ for $k_0\in\RR_{>0}$ in Figure~\ref{fig:imaginarypart} for different choices of determinants: $\mathrm{det}_1$ or $\mathrm{det}_2$, we also considered $\mathrm{det}_\mathrm{scal}$ corresponding to the scaling
\begin{equation}
\mathrm{det}_\mathrm{scal}(k)=k^2 \mathrm{det}_2(k)
\end{equation}
in order to investigate whether a simpler scaling than multiplying $\mathrm{det}_2$ by $\frac1{c_{1,k}c_{2,k}}$ would be suitable in practice.

As shown in Figure~\ref{fig:imaginarypart}, the considered quantity is always positive for $\mathrm{det}_1$ in contrast to $\mathrm{det}_2$ which has an oscillatory behavior. It appears that the simpler scaling by $k^2$ does not resolve the problem of negative imaginary parts on the real axis, and it seems that a scaling that is oscillatory with $k$ is advantageous. For certain values of $k_0$, however, the quantity $\frac{\det_2(k_0)}{\partial_k \mathrm{det}_2(k_0)}$ may be positive, and the Newton algorithm should work better for these initial starting points (tests performed for $k_0=17$ or $k_0=22$ for $\xi=0.5$, $n_1=1.5$, $n_2=1$ and $m=10$). For example, Figure~\ref{fig:newtontrajectories} and Table~\ref{tab:newtontrajectories} in \ref{appendix:figures} show the Newton iterations in the complex plane for two choices of the initial guess on the real axis, where one can observe either a convergence or a gross divergence using $\det_2$. 

\begin{figure}[h]
	\centering
	\centering
	\begin{subfigure}[b]{0.25\textwidth}
		\centering
		\includegraphics[width=\textwidth]{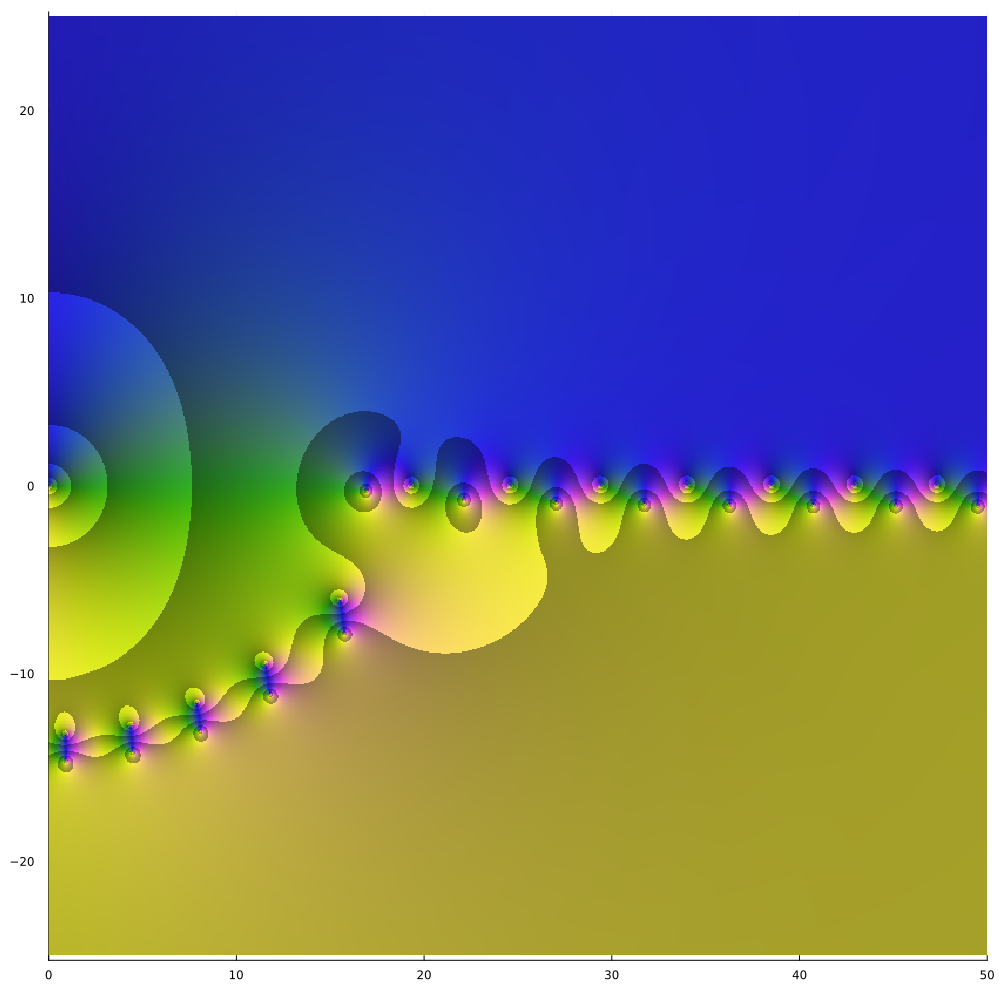}
		\caption{$\frac{\partial_k \det_1(k)}{\det_1(k)}$}
		\label{fig:derdetscal}
	\end{subfigure}\qquad \qquad
	\begin{subfigure}[b]{0.25\textwidth}
		\centering
		\includegraphics[width=\textwidth]{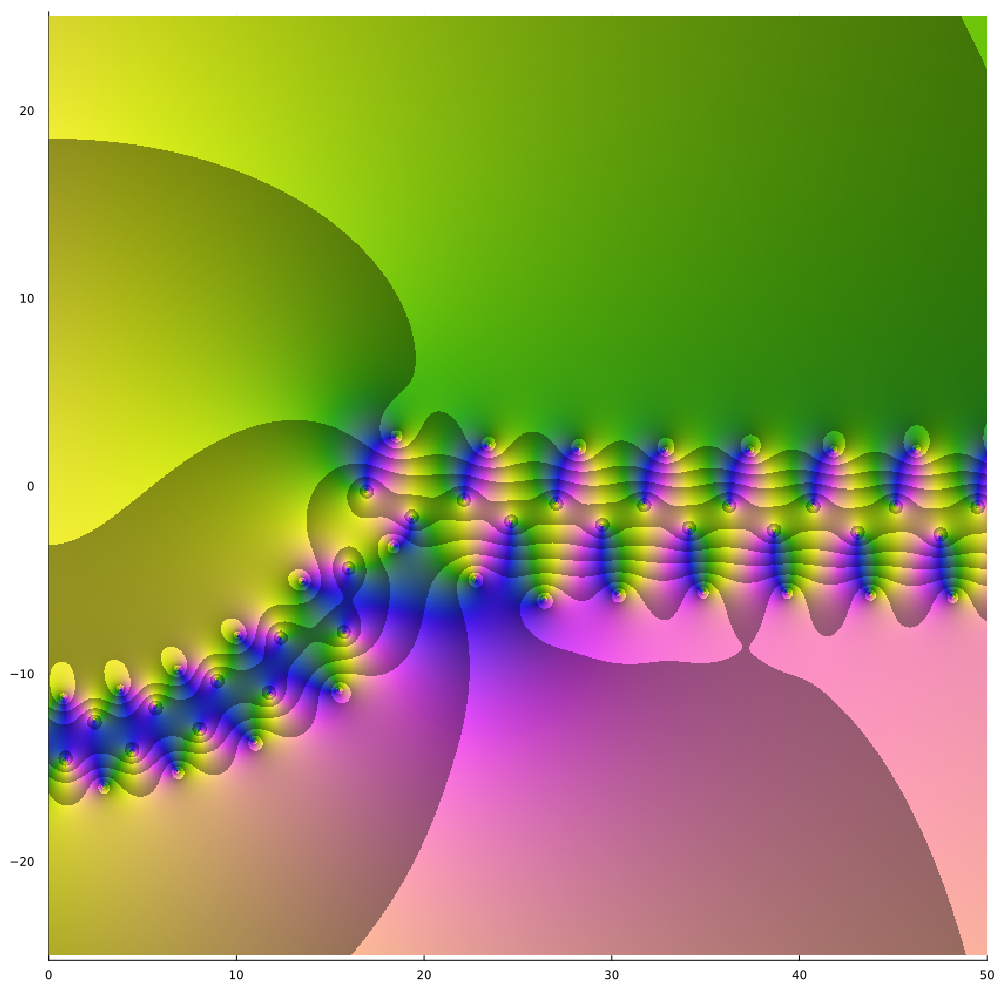}
		\caption{$\frac{\partial_k \det_2(k)}{\det_2(k)}$}
		\label{fig:derdetscal}
	\end{subfigure}\qquad \qquad
	\begin{subfigure}[b]{0.25\textwidth}
		\centering
		\includegraphics[width=\textwidth]{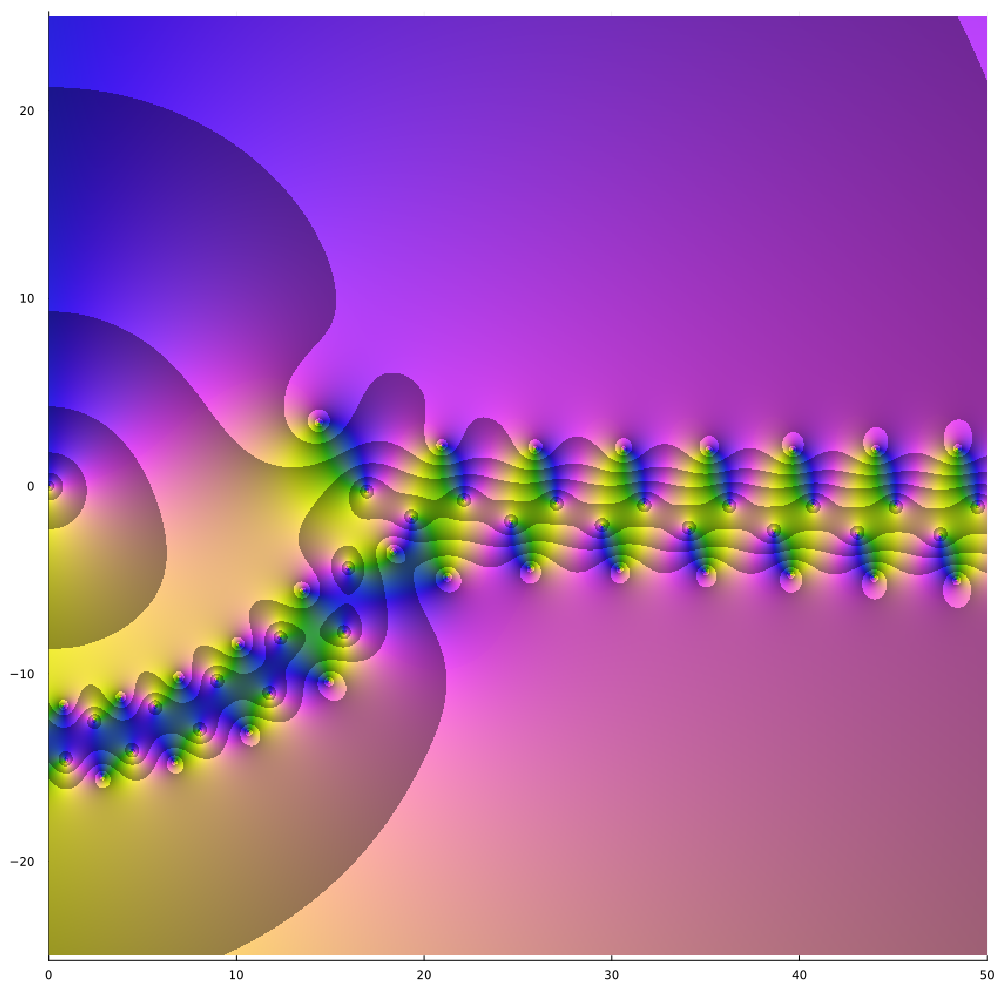}
		\caption{$\frac{\partial_k \det_\mathrm{scal}(k)}{\det_\mathrm{scal}(k)}$}
		\label{fig:derdetscal}
	\end{subfigure}\\
	\begin{subfigure}[b]{0.3\textwidth}
		\centering
		\includegraphics[width=\textwidth]{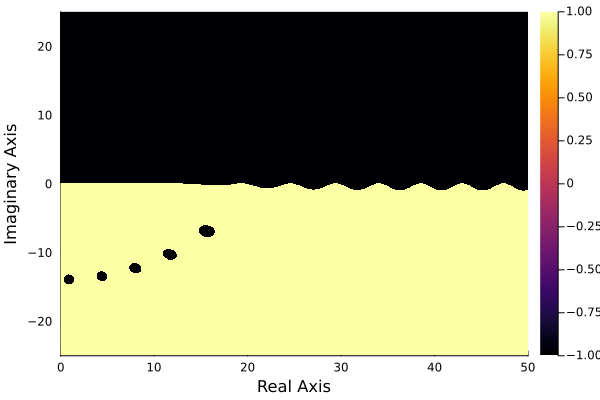}
		\caption{$\sign(\Im(\frac{\partial_k \det_1(k)}{\det_1(k)}))$}
		\label{fig:signderdet1ratio}
	\end{subfigure}
	\quad
	\begin{subfigure}[b]{0.3\textwidth}
		\centering
		\includegraphics[width=\textwidth]{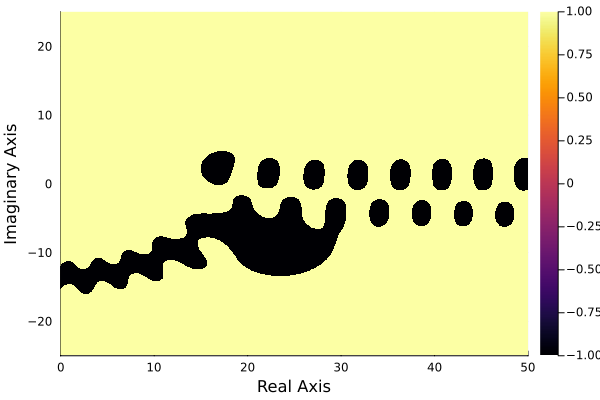}
		\caption{$\sign(\Im(\frac{\partial_k \det_2(k)}{\det_2(k)}))$}
		\label{fig:signderdet2ratio}
	\end{subfigure}
	\quad
	\begin{subfigure}[b]{0.3\textwidth}
		\centering
		\includegraphics[width=\textwidth]{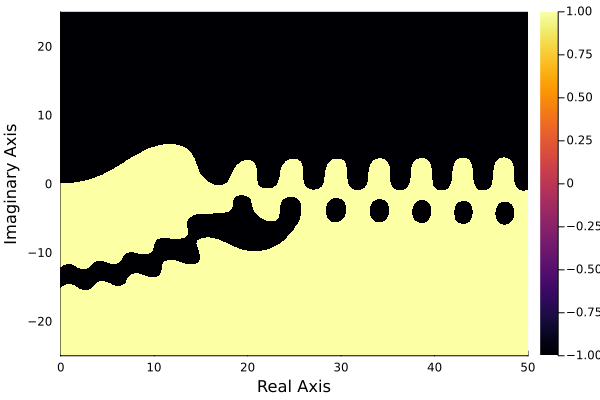}
		\caption{$\sign(\Im(\frac{\partial_k \det_\mathrm{scal}(k)}{\det_\mathrm{scal}(k)}))$}
		\label{fig:signderdetscalratio}
	\end{subfigure}
	\caption{Complex plots of different scalings of $\frac{\partial_k \det(k)}{\det(k)}$ for $\xi=0.5$, $n_1=1.5$, $n_2=1$ and $m=10$.}
	\label{fig:complexplotdetscal}
\end{figure}

An important observation affecting the Newton method using the different scalings can be found in the plots (see Fig.~\ref{fig:complexplotdetscal}) of $\sign(\Im(\frac{\partial_k \det(k)}{\det(k)}))$  related to the phase of $\frac{\partial_k \det(k)}{\det(k)}$ in the complex plane, extending the plot made on the real axis (Figure~\ref{fig:imaginarypart}): one significant advantage of $\det_1$ over the other scalings is that it could be observed in Figure~\ref{fig:signderdet1ratio} that the considered quantity is always negative for $\Im(k)\ge0$, meaning that if a given iterate $k_l$ in the Newton algorithm is such that $\Im(k_l)\ge0$, then the next applications of the Newton algorithm will have a tendency to send it back to the region $\Im(k)<0$ where we are looking for zeros. The other scalings $\det_2$ (Figure~\ref{fig:signderdet2ratio}) and $\det_\mathrm{scal}$ (Figure~\ref{fig:signderdetscalratio}) lack this property, which could explain why for some starting points on the real axis, the next iterates could have their imaginary part that keeps increasing in the region $\Im(k)>0$ (case of $\det_2$ where a big region inside $\Im(k)>0$ has a wrong sign), or may have an oscillatory behavior in some sub-region $0<\Im(k)<\alpha$ where the plotted sign is not constant (case of both $\det_2$ and $\det_\mathrm{scal}$).

Given that we intend to propose a Newton technique that starts on the real axis, we recommend to start the algorithm with a first iteration in the fourth quadrant, so that $\Im(k_1)<0$. The insights gained from the piecewise constant case will be useful to ensure that the problem is scaled properly also in the variable case, where we recommend, by analogy, the use of the scaling~\eqref{eq:databoundarychoice} for the aforementioned benefits. For this purpose, we will need an additional assumption in the variable case:
\begin{assumption}\label{assumption2}
	Let $\left(  \tilde{k},\tilde{u}\right)  \in\mathbb{C}_{<0}\times H_{0}^{1}\left(
	\Omega\right)\setminus\{0\}  $ be a resonance of problem~\eqref{eq:pdenep} and $m\in\mathbb{Z}$.
	Then, there exists a complex neighborhood $\omega_{\tilde{k}}\subset\mathbb{C}_{<0}$
	of $\tilde{k}$ such that $\frac{1}{c_{1,k}}\ne 0$ and $\frac{1}{c_{2,k}}\ne 0$ for all $k\in\omega_{\tilde{k}}$.
\end{assumption}
In this way, it is assumed that the resonances of the non-linear eigenvalue
problem~\eqref{eq:pdenep} are not simultaneously zeros of   $\frac{1}{c_{1,k}}$ or $\frac{1}{c_{2,k}}$.\footnote{Note that these zeros coincide with the critical frequencies for the auxiliary problems in (\ref{finaleqs}) in the piecewise constant case (Remark~\ref{rem:zerosc1kc2k}), thus the additional Assumption~\ref{assumption2} is already included in Assumption~\ref{assumption} in this case, but these assumptions are in general different in the variable case.}

\section{Numerical experiments}\label{sec:pwscase}

\subsection{BVP solver}\label{sec:odesolver}

In the variable case, the Newton method requires to solve the boundary value problems numerically for $f_{1,1,k}$, $f_{2,2,k}$, $\partial_k f_{1,1,k}$ and $\partial_k f_{2,2,k}$ to evaluate the determinant and its derivative. Hence, there is a need for reliable BVP solver of these two-point boundary value problems. 

Since the problems posed in $(0,\xi)$ and $(\xi,1)
$ are Bessel-type equations, their numerical approximation can be challenging given the singularities in the equation and the oscillatory behavior of the solution. To circumvent these problems, we use spectral methods in order to achieve high accuracy approximation. This is based on the asymptotic expansion of the solution on a given basis and we adopt the \texttt{ApproxFun} package in Julia language~\cite{approxfun,ApproxFun.jl-2014,townsend2015automatic,olver2013fast}. 
This appears to be more efficient and reliable for our application than using \texttt{NDSolve} of MATHEMATICA$^{\mbox{\scriptsize{\textregistered}}}$ and it is more suitable for numerical computations.

\subsection{Newton method for variable refractive index}\label{sec:ultimatennewton}

Based on the previous sections, we present here the full Newton algorithm to find (quasi-)resonances in the variable case. We specify in particular the choice of boundary conditions (BC), following \S\ref{sec:scaling}, and starting values, folowing \S\ref{sec:startingvalues}. The Newton algorithm we suggest is the following\footnote{A sample Julia code for this algorithm is available at \url{https://github.com/HankelGretel/WGM}.}%

\begin{algorithm}[H]
	\caption{Newton's algorithm for variable refractive index}
	\begin{algorithmic}
		\STATE $\mathrm{Newton}\_\mathrm{var} (\xi,n_1(r) \text{ on } (0,\xi),n_2(r) \text{ on } (\xi,1),m\in\NN,\varepsilon,l_{\max})$
		\STATE $\bullet$ Take the starting value: $k_0=\frac{m}{\xi n_1(\xi)};$
		\FOR{$l=0\ldots$ until stopping criterion is reached}
		\STATE $\bullet$ Start with $k=k_l;$
		\STATE $\bullet$ Solve two two-point BVPs numerically for $k=k_l$: result: approximations $f_{1,1,k}$ and $f_{2,2,k}$ to~\eqref{finaleqs} with the BC
		\begin{equation}\label{eq:databoundarychoice22}
		\begin{dcases}
		h_{1,1,k}:=k n_1(\xi) (J'_m(k \xi n_1(\xi))-iJ_m(k\xi n_1(\xi))) \quad \text{(assumed to be $\ne 0$)}\\
		h_{2,2,k}:=-k n_2(\xi) (H'_m(k \xi n_2(\xi))+i H_m(k \xi n_2(\xi))) \quad \text{(assumed to be $\ne 0$)};
		\end{dcases}
		\end{equation}
		\STATE $\bullet$ Compute numerically the derivatives with respect to $k$, $\partial_k f_{1,1,k}$ and $\partial_k f_{2,2,k}$, as the solution of two two-point BVPs~\eqref{eq:dkf11k},~\eqref{eq:dkf22k}; 
		\STATE $\bullet$ Compute $\mathrm{det}(k_l)$ by using~\eqref{eq:determinant2} and $\partial_k\mathrm{det}(k_l)$ by~\eqref{eq:derdeterminant2};
		\STATE $\bullet$ Compute
		\begin{equation}\label{eq:newton2}
		k_{l+1}=k_l-\frac{\det(k_l)}{\partial_k \det(k_l)};
		\end{equation}
		\STATE $\bullet$ $l\leftarrow l+1;$
		\STATE The stopping criterion is given by $\frac{\det(k_l)}{\|T(k_l)\|_F}\le \varepsilon$ or when a maximal number of iterations $l_{\max}$ is reached.
		\ENDFOR
	\end{algorithmic}
	\label{algovar}
\end{algorithm}

\subsection{Numerical setups}
\begin{figure}[h!]
	\centering
	\includegraphics[scale=0.65]{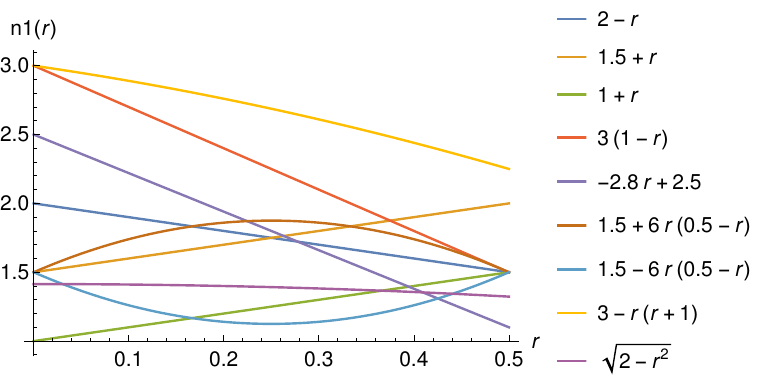}
	\caption{Different test cases in the non-constant case.}
	\label{fig:newtontestsn1}
\end{figure}

We fix without loss of generality $\xi=0.5$ and consider the following numerical setups, including piecewise constant cases, and non-constant cases for $n_1$ of affine or parabolic types. We also consider a special variable $n_1$ with explicit solution ($n_1(r)=\sqrt{2-r^2}$), corresponding to the Luneburg lens~\cite{boriskin2002whispering,lock2008scattering1,lock2008scattering2,lock2008scattering3} (see \ref{appendix:Whittaker}). Finally, we consider a more general case, where $n_1$ and $n_2$ are both variable, making it outside of the scope of the theorems of~\cite{balac2021asymptotics} regarding the asymptotics for resonances. The different setups are summarized in Table~\ref{tab:numericalsetup} and the considered cases for variable $n_1$ are illustrated in Figure~\ref{fig:newtontestsn1}.

\begin{table}[h!]
	\scriptsize
	\centering
	\begin{tabular}{l|cc|ccccc|ccc|ccc}  
		\toprule
		$n_1(r)$ & $\mathbf{1.5}$ & $\mathbf{5}$ & $2-r$  &  $1.5+r$ & $1+r$ & $3(1-r)$  & $-2.8r+2.5$ & $1.5+6r(\xi-r)$ & $1.5-6r(\xi-r)$ & $3-r(r+1)$ & \multicolumn{3}{c}{$\mathbf{\sqrt{2-r^2}}$} \\
		\midrule
		$n_2(r)$ & \multicolumn{2}{c|}{$\mathbf{1}$} & \multicolumn{5}{c|}{$1$} & \multicolumn{3}{c|}{$1$} & $\mathbf{1}$ & $r+0.5$ & $1+(r-0.5)^3$  \\
		\midrule
		Thm. & \multicolumn{2}{c|}{1.A} & \multicolumn{3}{c|}{1.A} & 1.B & 1.C & 1.B  & \multicolumn{2}{|c|}{1.A} & 1.A & \multicolumn{2}{|c}{NA}  \\
		\midrule
		Table & \ref{tab:Newtonconstant1.5} & \ref{tab:Newtonconstant5} & \ref{tab:Newtonnonconstantaffine1} & \ref{tab:Newtonnonconstantaffine2} &  \ref{tab:Newtonnonconstantaffine3} & \ref{tab:Newtonnonconstantaffine4} & \ref{tab:Newtonnonconstantaffine5}\ \ref{tab:Newtonnonconstantaffine5startleadingterm} & \ref{tab:Newtonnonconstantparabolic1} & \ref{tab:Newtonnonconstantparabolic2} & \ref{tab:Newtonnonconstantparabolic3} & \ref{tab:NewtonnonconstantLuneberg} & \ref{tab:Newtonnonconstantvarn21} & \ref{tab:Newtonnonconstantvarn22} \\
		& & & & & & & \ref{tab:Newtonnonconstantaffine5startkth} & & \ref{tab:Newtonnonconstantparabolic2asympt} & & & &  \\
		\bottomrule
	\end{tabular}
	\caption{Summary of the setups considered in the numerical experiments.}
	\label{tab:numericalsetup}
\end{table}

The results of the Newton algorithm are reported in the indicated Tables of the Supplementary Material. We also report on the tables the value $k_\mathrm{asympt}$ corresponding to the first terms of the asymptotics of $\underline{k}_{0}(m)$ up to $O(m^{-1})$ for Thm.~1.A and Thm.~1.C, and up to $O(m^{-\frac12})$  for Thm.~1.B (\ref{appendix:dauge}).

Our algorithm is first validated in cases where fundamental solutions are known explicitly (cases in bold in Table~\ref{tab:numericalsetup}) by comparing the results obtained using the Newton method with \texttt{ApproxFun} and using the exact expression of the determinant in Tables~\ref{tab:comaprison Newton},~\ref{tab:comaprison Newton cont},~\ref{tab:NewtonLaguerre} and~\ref{tab:NewtonLaguerre Cont}. The results are very close (difference of order $10^{-11}-10^{-15}$), which validates our methodology and implementation.

\subsection{Interpretation of results}

Based on the numerical experiments, the following observations can be made regarding the convergence of the Newton method presented in Section~\ref{sec:ultimatennewton}.

\paragraph*{Piecewise constant case:}
\begin{itemize}
	\item  The Newton method converges for $n_1=1.5$ and $n_2=1$ (Table~\ref{tab:Newtonconstant1.5}), for all the considered $m\le 60$ in a number of iterations $l\le 10$.
	\item In high-contrast media, e.g. $n_1=5$ and $n_2=1$ (Table~\ref{tab:Newtonconstant5}), the Newton method also converges for all the considered $m\le 60$ in a number of iterations $l\le 11$. The use of the stopping criterion based on the relative residual $\frac{\det(k)}{\|T(k)\|_F}$ is crucial in this context, as 
	the computed determinant can be very large yet negligible compared to the magnitude of $T$. This allows us to remove the sensitivity of the Newton algorithm with respect to the jump in the refractive index $n^2$ at the interface (low/high-contrast media).
	\item The obtained $k$ is close to the corresponding theoretical asymptotic $k_\mathrm{asympt}$, with a difference of order $10^{-2}$ (low-contrast media) and $10^{-4}$ (high-contrast media) for high values of $m\ge28$.  
\end{itemize}

\paragraph*{Variable case:}
\begin{itemize}
	\item The Newton method converges in all the considered cases using the same formula for the initial guess ($\frac{m}{\xi n_1(\xi)}$), in a number of iterations $l\le 19$.
	\item The Newton method also converges in the case with the highest contrast in the refractive index at the interface (see Figure~\ref{fig:newtontestsn1} and corresponding Table~\ref{tab:Newtonnonconstantparabolic3}), when using the aforementioned stopping criterion, whereas it may fail (maximal number of iterations can be reached) if one uses the simplest stopping criterion $|\det{k}|\le \varepsilon$, for the same reasons as in the piecewise constant case. We note that the imaginary part of the resonance in these cases approaches $0$ (values of order $10^{-12}-10^{-17}$ for high values of $m$).
	\item The considered example satisfying Thm.~1.C converges starting from $k_0=\frac{m}{\xi n_1(\xi)}$ (Table~\ref{tab:Newtonnonconstantaffine5}) which is different from the leading term in the asymptotics in this case (which is $\frac{m}{\xi_0 n_1(\xi_0)}$ with $\xi_0=\frac{25}{56}\approx 0.44$, see \ref{appendix:dauge}). A comparison with the results obtained starting from $k_\mathrm{asympt}$ (Table~\ref{tab:Newtonnonconstantaffine5startkth}) gives comparable results for large $m$ that are also close to $k_\mathrm{asympt}$. However, starting from the leading term $\frac{m}{\xi_0 n_1(\xi_0)}$ (Table~\ref{tab:Newtonnonconstantaffine5startleadingterm}) gives different results for the resonances, that are quite far from $k_\mathrm{asympt}$.  This shows a high sensitivity of the algorithm to the initial guess, especially in this case. In the other cases (Thm.~1.A or Thm.~1.B), the obtained resonances are close to $k_\mathrm{asympt}$ when starting from the leading term of the asymptotics ($k_0=\frac{m}{\xi n_1(\xi)}$ for these theorems), with a difference of orders ranging from $10^{-1}$ to $10^{-3}$ for high values of $m$, except in one case (Tables~\ref{tab:Newtonnonconstantparabolic2} and \ref{tab:Newtonnonconstantparabolic2asympt}). 
	This difference of behavior for Thm.~1.C may be linked to the fact that quasi-resonances are not of whispering gallery type in this case: the quasi-modes are strictly localized inside the cavity (around $r=\xi_0<\xi$)~\cite{balac2021asymptotics}.\\  Nevertheless, this is a good point for our algorithm~\ref{algovar} where we specified the same initial guess independently of the verified theorem. This is particularly interesting when none of the theorems is valid (see next point).
	\item In the case of variable $n_1$ and $n_2$, theoretical asymptotic expansions are not available. In the considered cases (Tables~\ref{tab:Newtonnonconstantvarn21} and~\ref{tab:Newtonnonconstantvarn22}), taking $n_2$ such that $n_2(\xi)=1$ and $n_2'(\xi)=n_2''(\xi)=0$ (Table~\ref{tab:Newtonnonconstantvarn22}) gives results that are asymptotically close to the case $n_2=1$ (Table~\ref{tab:NewtonnonconstantLuneberg}), while a different variable $n_2$ satisfying only $n_2(\xi)=1$ (Table~\ref{tab:Newtonnonconstantvarn21}) yields different results.
\end{itemize}

\subsection{Quasi-resonance solutions and exact modes}\label{quasi}

\begin{figure}[h]
	\centering
	\includegraphics[scale=0.37]{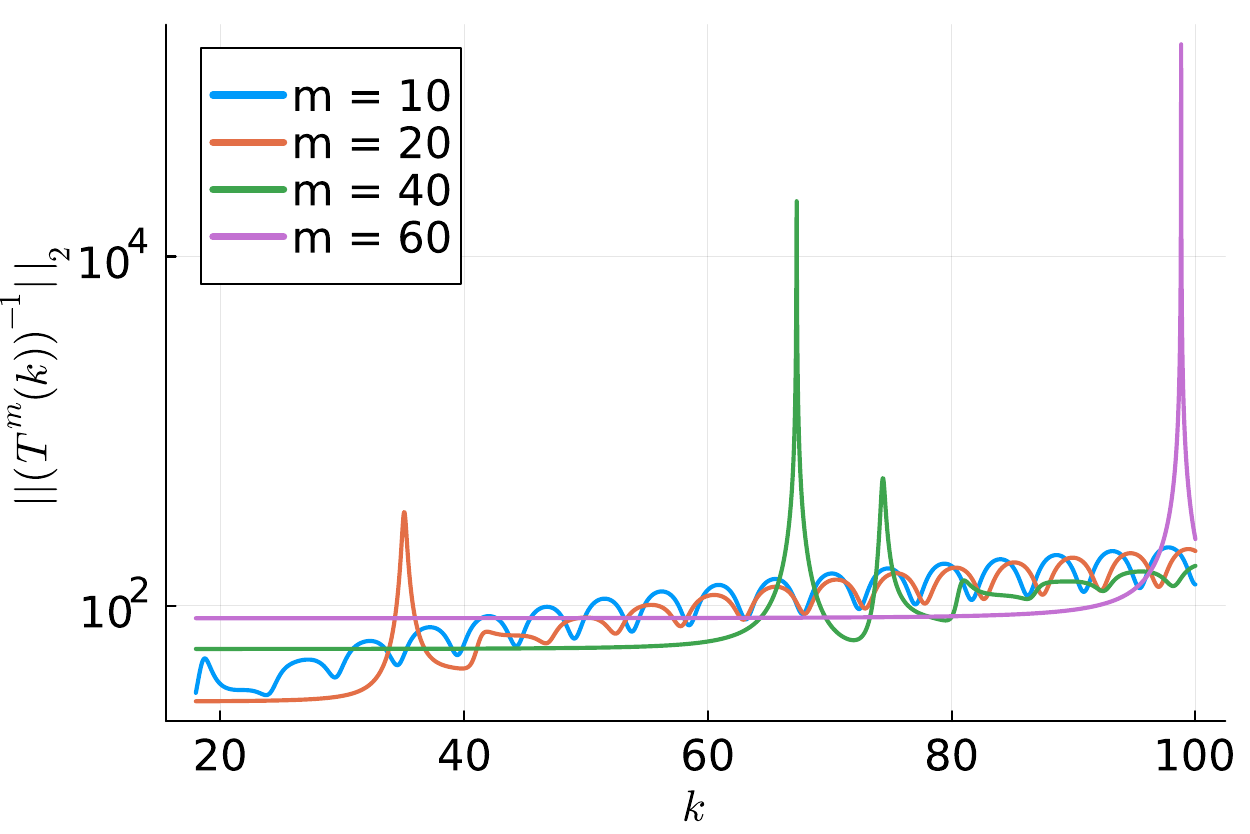}
	\caption{Plot of $\|(T^m(k))^{-1}\|_2$ with respect to $k$ for different values of $m$, in the Luneburg case $n_1(r)=\sqrt{2-r^2}$, $n_2(r)=1$ and $\xi=0.5$. Here $h_{1,1,k}=h_{2,2,k}=1$.}
	\label{fig:opnorm}
\end{figure}

\begin{figure}[H]
	\centering
	\begin{subfigure}[b]{0.38\textwidth}
		\centering
		\includegraphics[width=\textwidth]{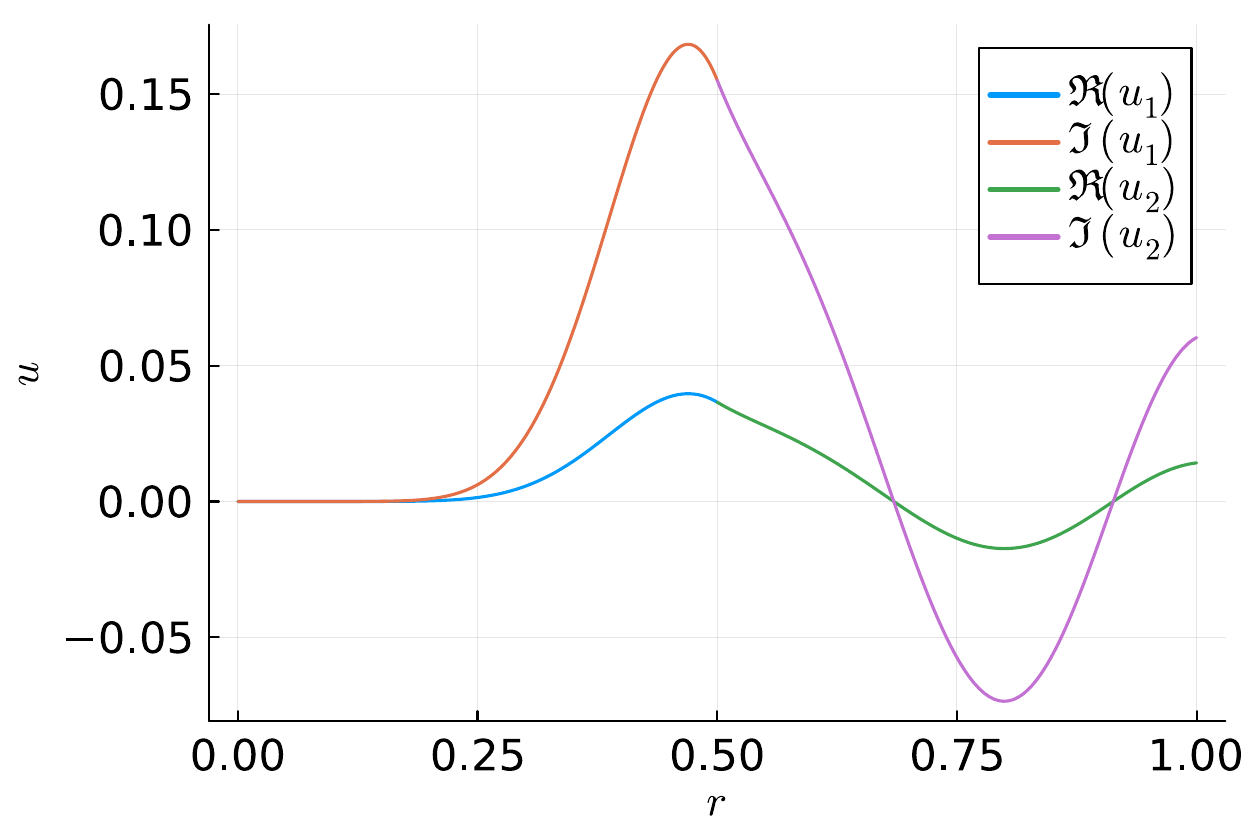}
		\caption{$m=10$, $k=18.588963438926466$}
		\label{fig:Lunebergmequals10}
	\end{subfigure}
	\hfill
	\begin{subfigure}[b]{0.38\textwidth}
		\centering
		\includegraphics[width=\textwidth]{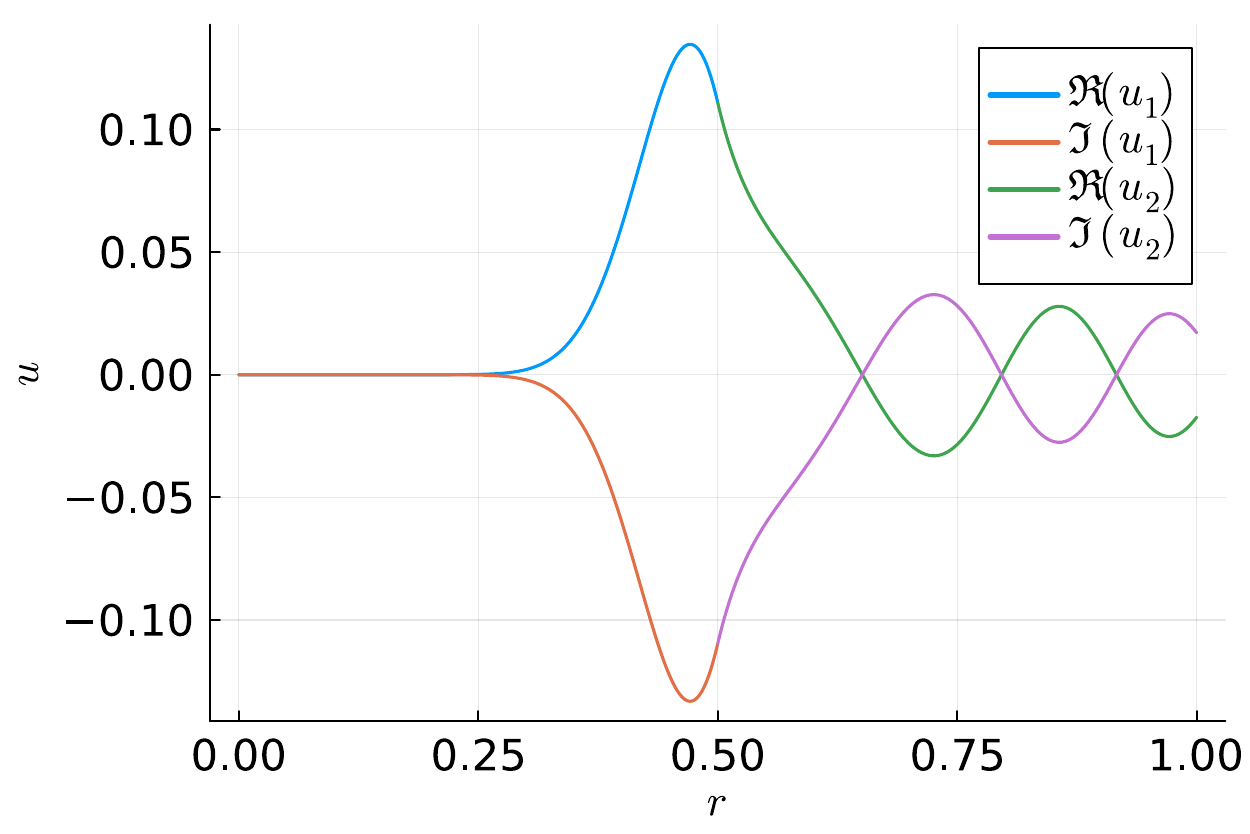}
		\caption{$m=20$, $k=35.09408648067281$}
		\label{fig:Lunebergmequals20}
	\end{subfigure}
	\hfill
	\begin{subfigure}[b]{0.38\textwidth}
		\centering
		\includegraphics[width=\textwidth]{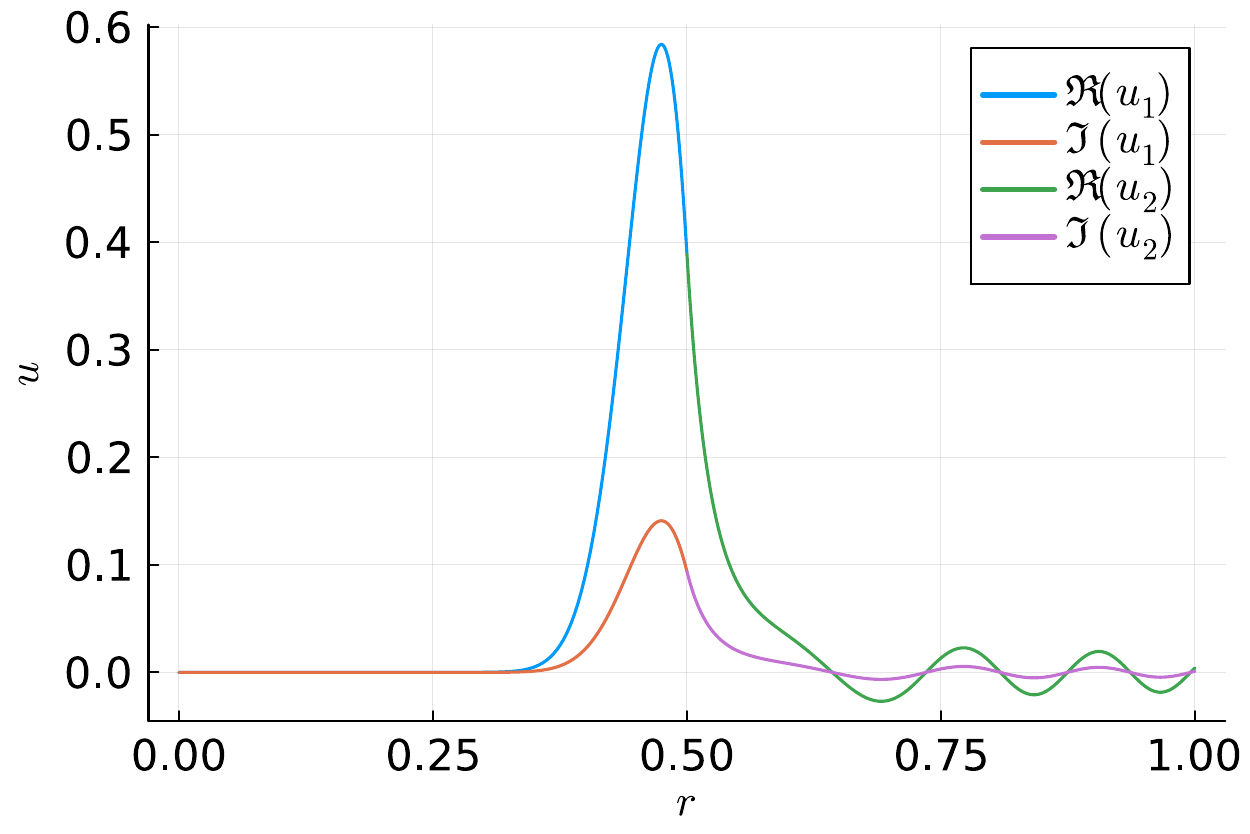}
		\caption{$m=40$, $k=67.28740148972052$}
		\label{fig:Lunebergmequals40}
	\end{subfigure}
	\hfill
	\begin{subfigure}[b]{0.38\textwidth}
		\centering
		\includegraphics[width=\textwidth]{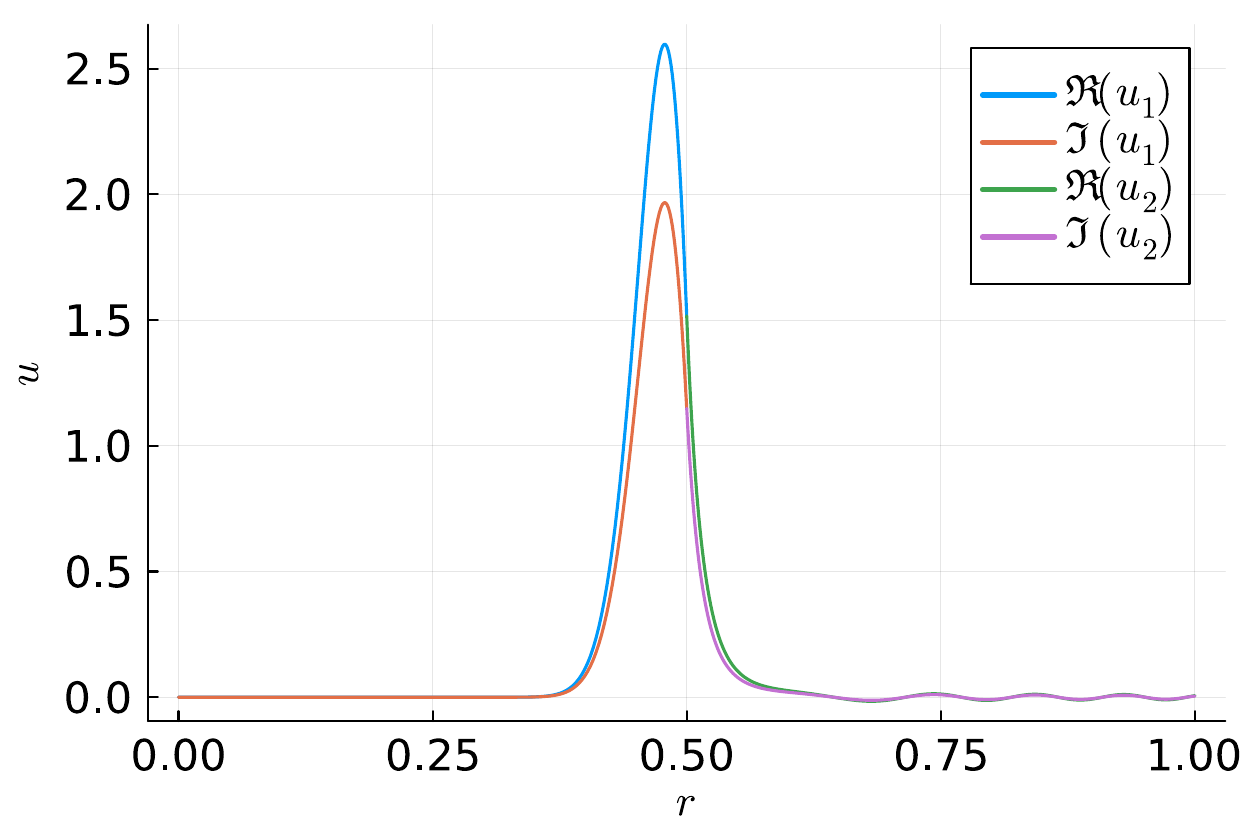}
		\caption{$m=60$, $k=98.82822050605951$}
		\label{fig:Lunebergmequals60}
	\end{subfigure}
	\begin{subfigure}[b]{0.38\textwidth}
		\centering
		\includegraphics[width=\textwidth]{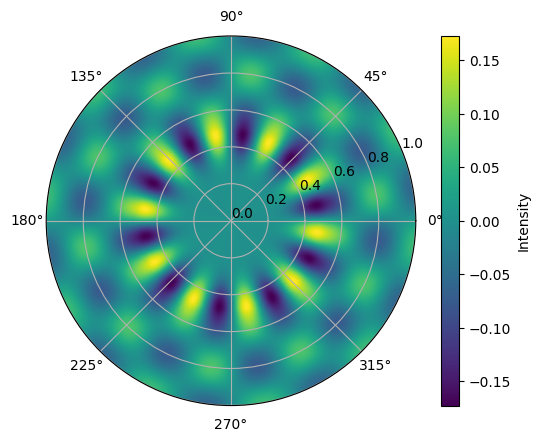}
		\caption{Case (a)}
		\label{fig:PolarLunebergmequals10}
	\end{subfigure}
	\hfill
	\begin{subfigure}[b]{0.38\textwidth}
		\centering
		\includegraphics[width=\textwidth]{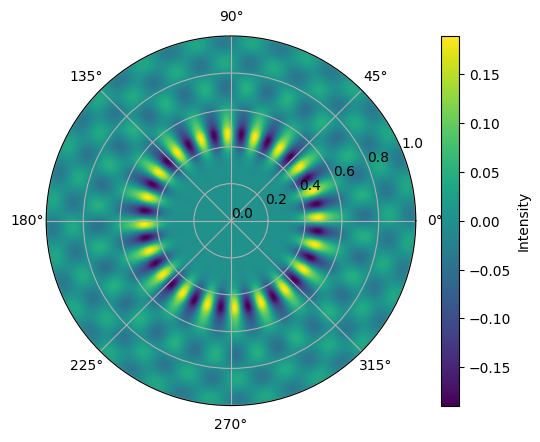}
		\caption{Case (b)}
		\label{fig:PolarLunebergmequals20}
	\end{subfigure}
	\hfill
	\begin{subfigure}[b]{0.38\textwidth}
		\centering
		\includegraphics[width=\textwidth]{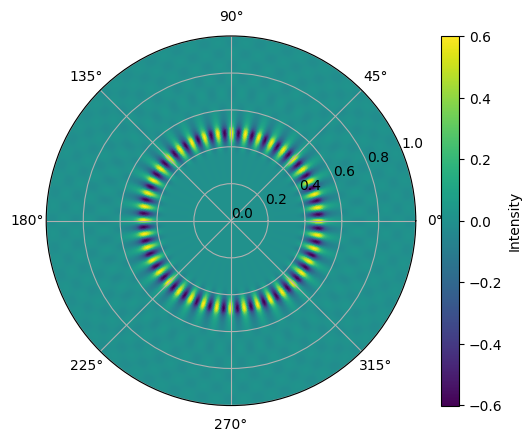}
		\caption{Case (c)}
		\label{fig:PolarLunebergmequals40}
	\end{subfigure}
	\hfill
	\begin{subfigure}[b]{0.38\textwidth}
		\centering
		\includegraphics[width=\textwidth]{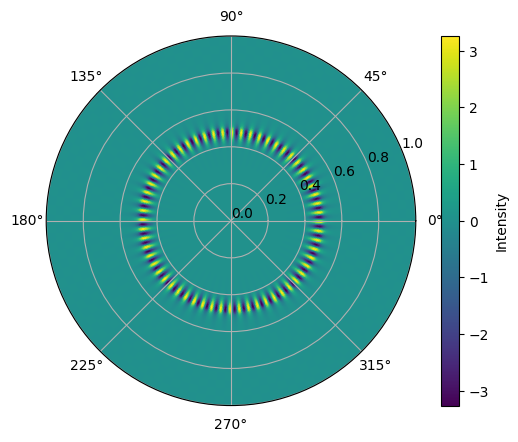}
		\caption{case (d)}
		\label{fig:PolarLunebergmequals60}
	\end{subfigure}
	\caption{(a)-(d): Plots of $\Re(\hat{u}_m(r))$ and $\Im(\hat{u}_m(r))$; (e)-(h): Polar plots of $\Re[\hat{u}_m(r) e^{im\theta}]$, for $\hat{u}_m(r)$  solution of~\eqref{eq:helmholtzpolarcordinates} in the Luneburg case $n_1(r)=\sqrt{2-r^2}$, $n_2(r)=1$ and $\xi=0.5$, where $k$ is a quasi-resonance (close to a resonance), see Table~\ref{tab:NewtonnonconstantLuneberg}. Here $\hat{g}_m=1$.}
	\label{fig:solution}
\end{figure}

After computing a resonance $k^*$, we can solve~\eqref{eq:linearsystemquasi} 
for the quasi-resonance $k=\underline{k}=\Re(k^*)$ and plot the solution $\hat{u}_m(r)$~\eqref{eq:ansatzquasi} and $\Re[\hat{u}_m(r) e^{im\theta}]$ in Figure~\ref{fig:solution}.  
We observe that they concentrate around the interface $r=\xi$ as $m$ gets large ($\Im(k^*)$ gets closer to $0$). 
In Figure~\ref{fig:opnorm}, we represented the norm $\|(T^m(k))^{-1}\|_2$ with respect to $k$ for the different values of $m$ considered: we observe spikes in this norm when $k$ is a quasi-resonance, since the matrix is close to being singular in this case ($\det(T^m(\underline{k}))\approx \det(T^m(k^*))= 0$).

\begin{figure}[H]
	\centering
	\begin{subfigure}[b]{0.38\textwidth}
		\centering
		\includegraphics[width=\textwidth]{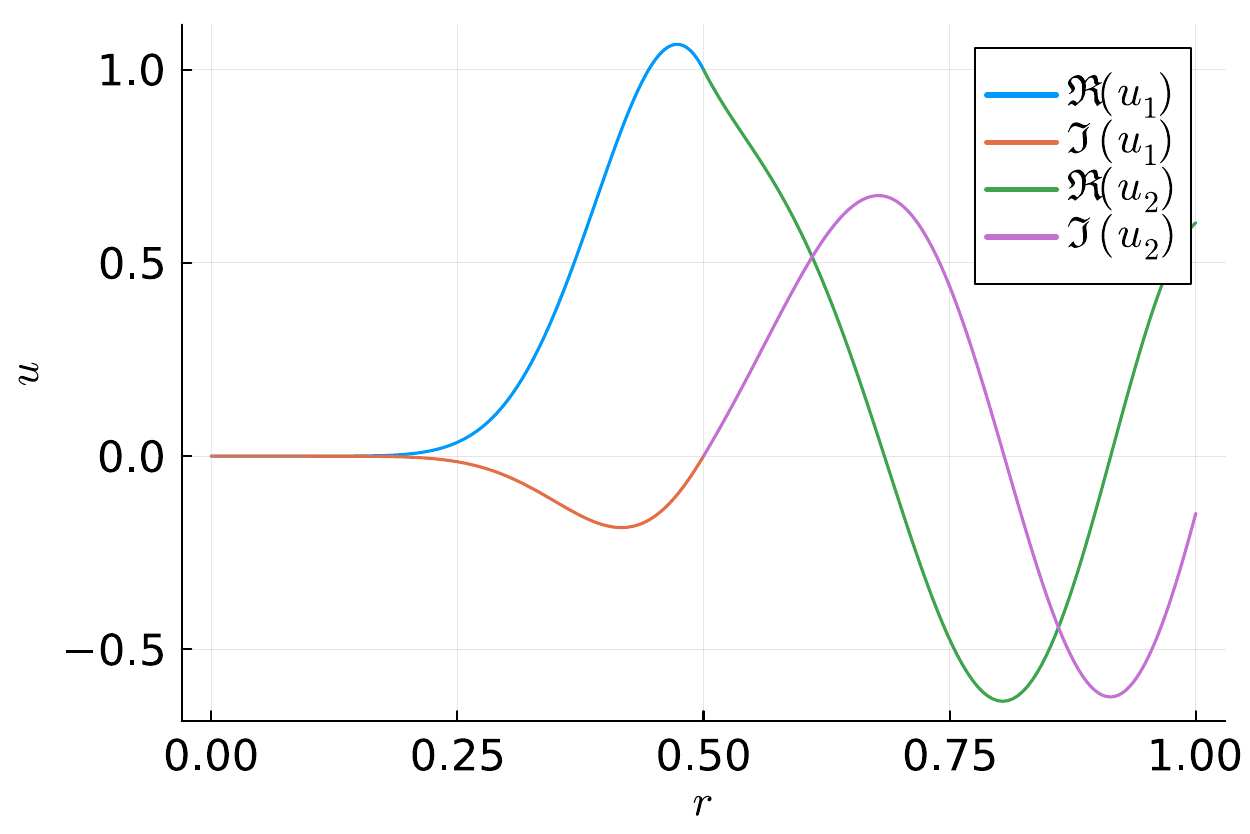}
		\caption{$m=10$, $k=18.588963438926466 - 0.6154425735324377 \,i$}
		\label{fig:modeLunebergmequals10}
	\end{subfigure}
	\hfill
	\begin{subfigure}[b]{0.38\textwidth}
		\centering
		\includegraphics[width=\textwidth]{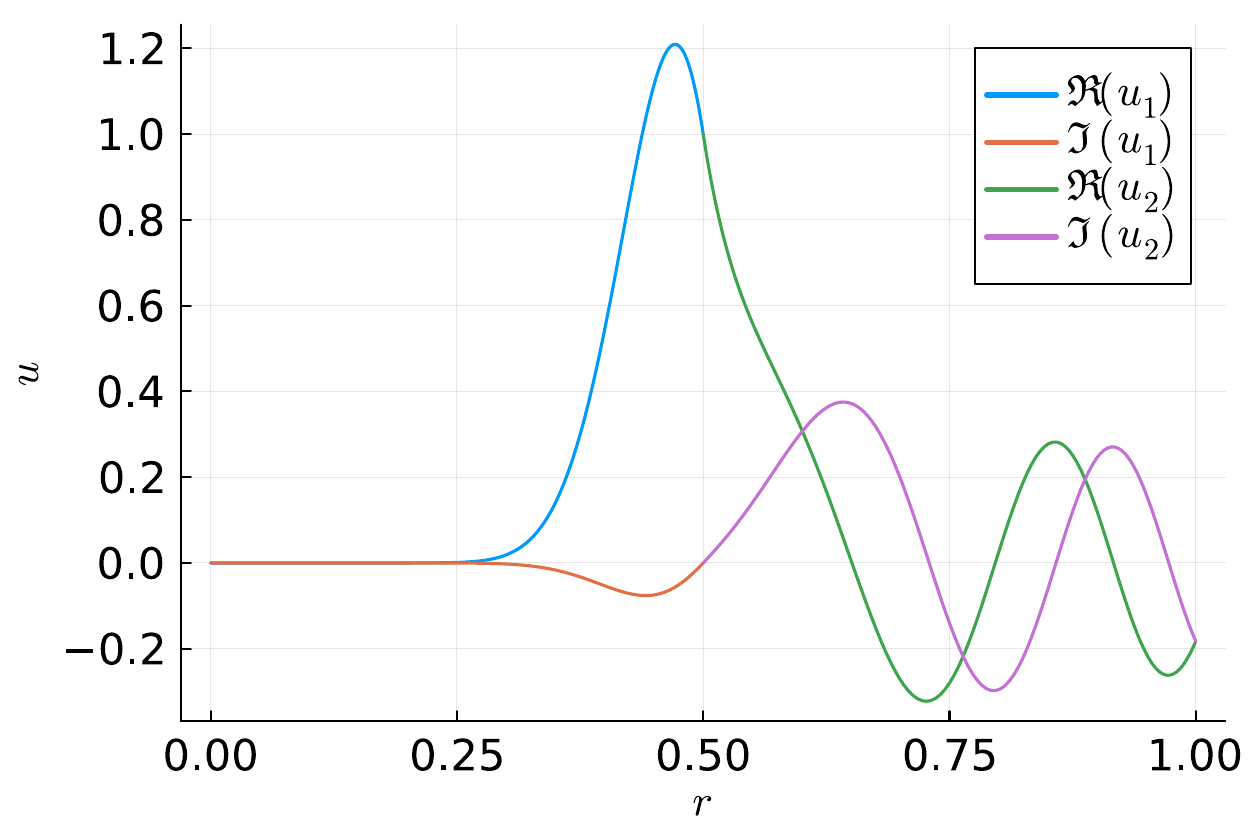}
		\caption{$m=20$, $k=35.09408648067281 - 0.19327141118804717\, i$}
		\label{fig:modeLunebergmequals20}
	\end{subfigure}
	\hfill
	\begin{subfigure}[b]{0.38\textwidth}
		\centering
		\includegraphics[width=\textwidth]{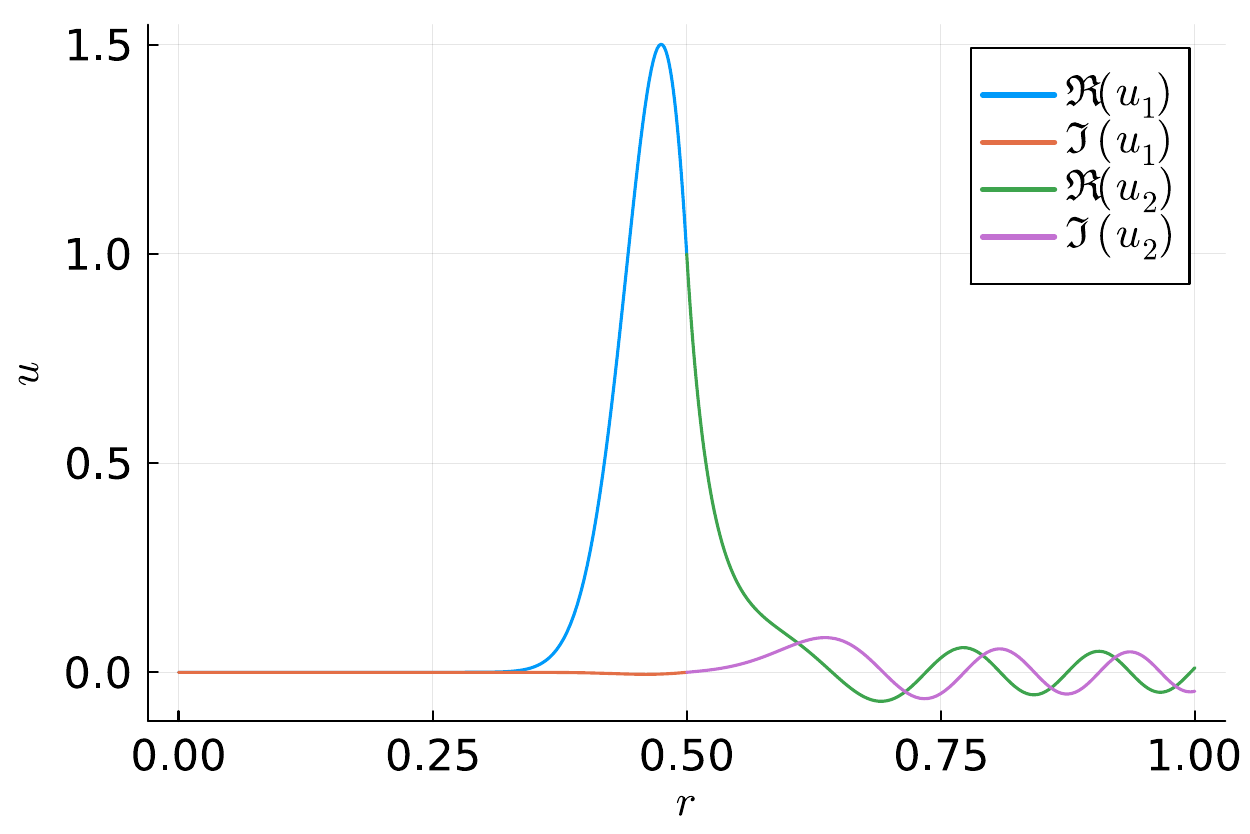}
		\caption{$m=40$, $k=67.28740148972052 - 0.008096455718707863\, i$}
		\label{fig:modeLunebergmequals40}
	\end{subfigure}
	\hfill
	\begin{subfigure}[b]{0.38\textwidth}
		\centering
		\includegraphics[width=\textwidth]{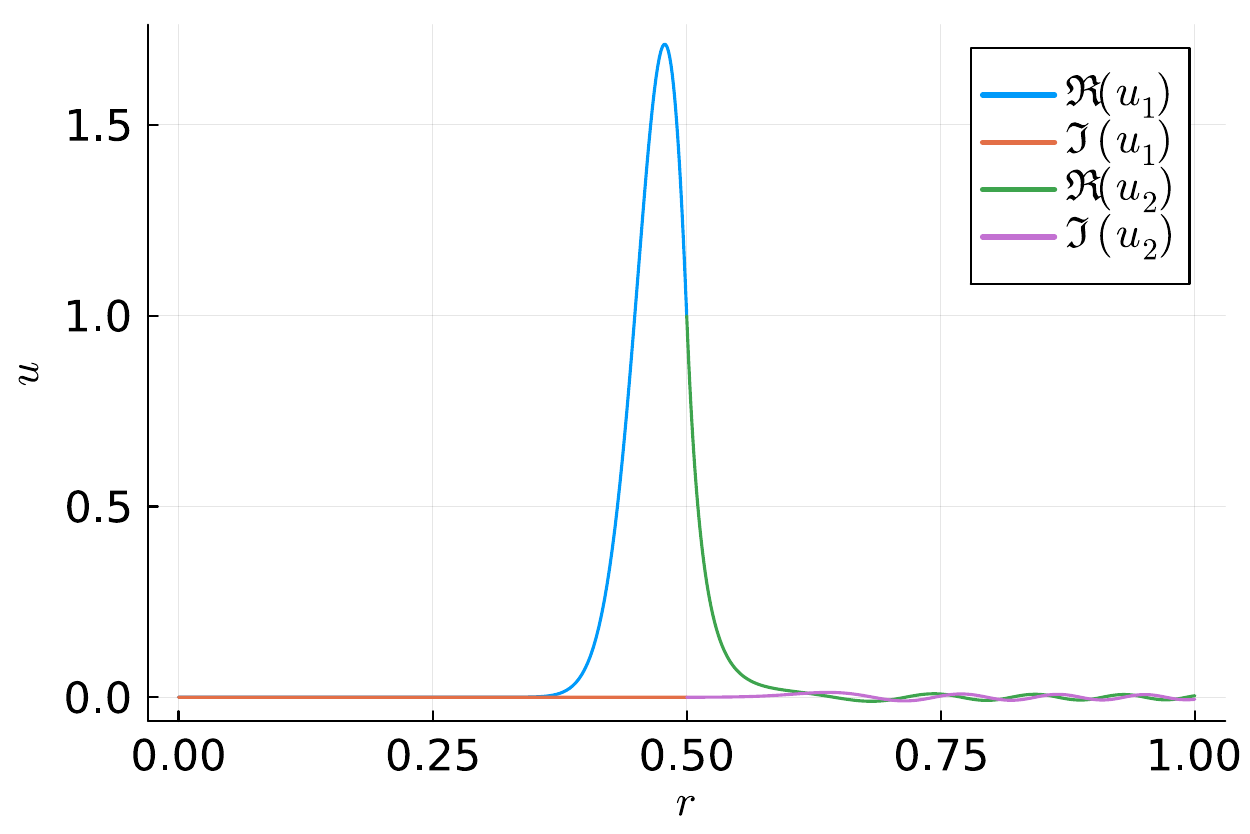}
		\caption{$m=60$, $k=98.82822050605951 - 0.0001666858070041872\, i$}
		\label{fig:modeLunebergmequals60}
	\end{subfigure}
	\begin{subfigure}[b]{0.38\textwidth}
		\centering
		\includegraphics[width=\textwidth]{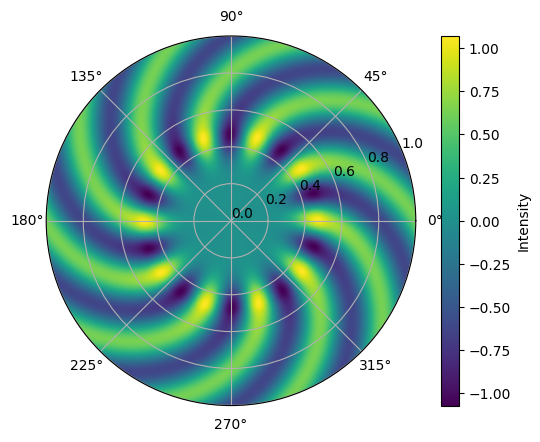}
		\caption{Case (a)}
		\label{fig:PolarmodeLunebergmequals10}
	\end{subfigure}
	\hfill
	\begin{subfigure}[b]{0.38\textwidth}
		\centering
		\includegraphics[width=\textwidth]{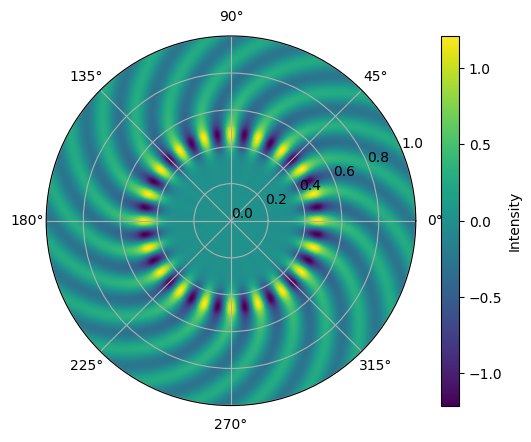}
		\caption{Case (b)}
		\label{fig:PolarmodeLunebergmequals20}
	\end{subfigure}
	\hfill
	\begin{subfigure}[b]{0.38\textwidth}
		\centering
		\includegraphics[width=\textwidth]{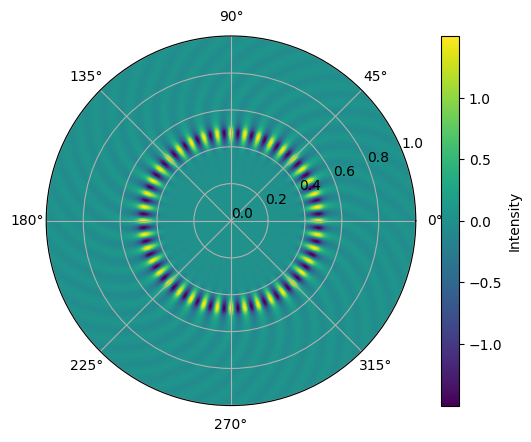}
		\caption{Case (c)}
		\label{fig:PolarmodeLunebergmequals40}
	\end{subfigure}
	\hfill
	\begin{subfigure}[b]{0.38\textwidth}
		\centering
		\includegraphics[width=\textwidth]{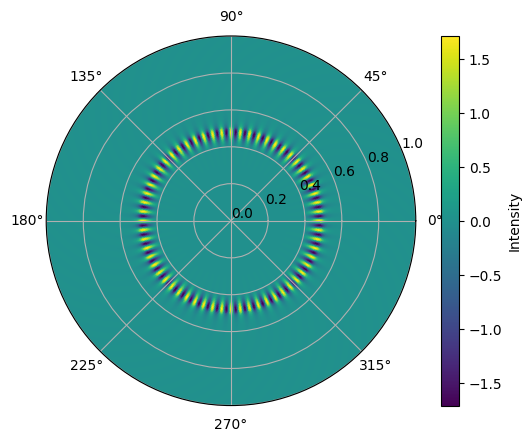}
		\caption{Case (d)}
		\label{fig:PolarmodeLunebergmequals60}
	\end{subfigure}
	\caption{(a)-(d): Plots of $\Re(\hat{u}_m(r))$ and $\Im(\hat{u}_m(r))$; (e)-(h): Polar plots of $\Re[\hat{u}_m(r) e^{im\theta}]$, for the exact mode $\hat{u}_m(r)$~\eqref{eq:exactmode} solution of~\eqref{eq:helmholtzpolarcordinates} for $\hat{g}_m=0$, where $k$ is a resonance, in the Luneburg case $n_1(r)=\sqrt{2-r^2}$, $n_2(r)=1$ and $\xi=0.5$, see Table~\ref{tab:NewtonnonconstantLuneberg}.}
	\label{fig:exactmode}
\end{figure}

We also plot the exact modes for a given resonance $k^*$: they are given, up to a complex multiplicative constant by,
\begin{equation}\label{eq:exactmode}
u(r)=
\begin{dcases}
\frac{f_{1,1,k^*}(r)}{f_{1,1,k^*}(\xi)} & r\in(0,\xi)\\
\frac{f_{2,2,k^*}(r)}{f_{2,2,k^*}(\xi)}    & r\in (\xi,1)
\end{dcases}  
\quad \text{ or }\quad
u(r)=\begin{dcases}
\frac{f_{1,1,k^*}(r)}{f'_{1,1,k^*}(\xi)} & r\in(0,\xi)\\
\frac{f_{2,2,k^*}(r)}{f'_{2,2,k^*}(\xi)}    & r\in (\xi,1),
\end{dcases} 
\end{equation}
see Proposition~\ref{prop:geometricsimplicity}. 
The exact modes are shown in Figure~\ref{fig:exactmode} for the Luneburg case  
(see Table~\ref{tab:NewtonnonconstantLuneberg}), exhibiting the same localization behavior at the interface  for large $m$ as the ``quasi-modes" in Figure~\ref{fig:solution}. 

\section{Conclusion and outlook}

In this paper, we have considered the problem of computing quasi-resonances in spherical symmetric heterogeneous Helmholtz problems with piecewise smooth refractive index. We have developed a general approach based on the splitting of the problem into decoupled problems on each domain and one problem at the interface (skeleton). In the spherical symmetric setting, we have reduced the problem to one-dimensional Sturm-Liouville problems. Using fundamental solutions, the problem of resonances is expressed as a nonlinear eigenvalue problem $T(k)x=0$ involving the values of the fundamental solutions and their derivatives at the interface. We then developed a numerical approach using Newton's method to solve the nonlinear equation $\det(T(k))=0$ where the fundamental solutions and their derivatives are approximated numerically at each iteration.  In the piecewise constant case, we prove the simplicity of the roots, providing a local quadratic convergence of the algorithm. In the variable case, we have specified the initial guess based on known asymptotic expansions, and suggested a proper scaling for the fundamental system in analogy with the piecewise constant case. Various numerical experiments, investigating the convergence and comparing the results to explicit solutions or asymptotic expansions, validate the results and our methodology. Perspectives include comparing the Newton method with other approaches developed for nonlinear eigenvalue problems, such as the contour integral method~\cite{beyn2012integral,asakura2009numerical} or rational approximation~\cite{bruno2024evaluation,SantanaWaves,guttel2024randomized}, and investigating the extension of this method to more general settings following for instance~\cite{alves2006numerical,alves2024wave,alves2020advances}. 
It would also be interesting to consider a piecewise smooth refractive index with $p$ jump points as in~\cite{sauter2021heterogeneous}, for its potential relevance in physical applications (layered or stratified media)~\cite{boriskin2002whispering,lock2008scattering3}.

\paragraph*{Acknowledgements.} 

Part of the research was carried out while the first author was visiting the Institut für Mathematik at the University of Zürich as part of the CIMPA-ICTP Fellowship 2024 ``Research in Pairs''. The first author thanks CIMPA and ICTP for the grant and the Institut für Mathematik for the hospitality. We thank Monique Dauge, Andrea Moiola and Zoïs Moitier for helpful discussions on WGM at the WAVES 2024 conference, and for pointing out the simplicity of eigenvalues for Sturm-Liouville problems.

\bibliographystyle{abbrv}
\bibliography{references}

\pagebreak
\bigskip
\bigskip
\bigskip

\begin{center}
\LARGE Supplementary Material: Computation of whispering gallery modes for spherical symmetric, heterogeneous Helmholtz problems with piecewise smooth refractive index\\
\medskip
\large{Bouchra Bensiali and Stefan Sauter}
\end{center}

\bigskip
\bigskip
\bigskip

\setcounter{equation}{0}
\setcounter{figure}{0}
\setcounter{table}{0}
\setcounter{page}{1}
\setcounter{section}{0}
\makeatletter
\renewcommand{\theequation}{S-\arabic{equation}}
\renewcommand{\thefigure}{S-\arabic{figure}}
\renewcommand{\thetable}{S-\arabic{table}}
\renewcommand{\thepage}{S-\arabic{page}}
\renewcommand{\thesection}{S\arabic{section}}
\renewcommand{\theHsection}{S\arabic{section}}

\addtocontents{toc}{\protect\setcounter{tocdepth}{2}} 

\section*{Note}
\noindent This document provides the detailed results of the numerical experiments performed using Algorithm~\ref{algovar} (Newton's algorithm for variable refractive index). We also include some appendices providing supporting information. We recall that $k_0$ is the initial guess and $k_\mathrm{asympt}$ corresponds to the first terms of the asymptotics of $\underline{k}_{0}(m)$ up to $O(m^{-1})$ for Thm.~1.A and Thm.~1.C, and up to $O(m^{-\frac12})$  for Thm.~1.B (see \ref{appendix:dauge}).

\tableofcontents

\clearpage
\section{Piecewise constant case validation}\label{sec:val1}

\begin{table}[h!]
	\centering
	\begin{tabular}{l|ccccc}  
		\toprule
		$k_0$  & $l$ & $k$  & $|D(k)|$ & $|\partial_k D(k)|$ &  Difference in $k$\\
		\midrule
		\midrule
		1&10&16.92320186058385 - 0.23954559046213897im&8.43e-11&1.19e-01&2.86e-14 \\
		& 10 & 16.923201860583823 - 0.23954559046216217im & 8.43e-11 & 1.19e-01 & \\
		\midrule
		2&9&16.923201860607477 - 0.239545590346429im&7.06e-11&1.19e-01&1.16e-13 \\
		& 9 & 16.923201860607584 - 0.23954559034648804im & 7.06e-11 & 1.19e-01 & \\
		\midrule
		3&8&16.923201958814623 - 0.23954567846387156im&1.58e-08&1.19e-01&2.73e-13 \\
		& 8 & 16.92320195881486 - 0.2395456784639192im & 1.58e-08 & 1.19e-01 & \\
		\midrule
		4&8&16.923201860703298 - 0.2395455900407501im&3.34e-11&1.19e-01&1.26e-13 \\
		& 8 & 16.923201860703177 - 0.23954559004071005im & 3.34e-11 & 1.19e-01 & \\
		\midrule
		5&7&16.92320302561573 - 0.23954583736498444im&1.42e-07&1.19e-01&6.36e-15 \\
		& 7 & 16.923203025615667 - 0.2395458373649745im & 1.42e-07 & 1.19e-01 & \\
		\midrule
		6&7&16.923201876723557 - 0.2395456186602138im&3.93e-09&1.19e-01&1.83e-13 \\
		& 7 & 16.923201876723734 - 0.23954561866023547im & 3.93e-09 & 1.19e-01 & \\
		\midrule
		7&7&16.9232018606048 - 0.23954559035873113im&7.21e-11&1.19e-01&2.07e-14 \\
		& 7 & 16.923201860604816 - 0.23954559035875164im & 7.21e-11 & 1.19e-01 & \\
		\midrule
		8&6&16.92320662364213 - 0.23954499539025922im&5.73e-07&1.19e-01&1.10e-13 \\
		& 6 & 16.923206623642024 - 0.23954499539021604im & 5.73e-07 & 1.19e-01 & \\
		\midrule
		9&6&16.923202147984654 - 0.2395457528343357im&3.94e-08&1.19e-01&4.70e-15 \\
		& 6 & 16.923202147984608 - 0.2395457528343282im & 3.94e-08 & 1.19e-01 & \\
		\midrule
		10&6&16.923201864447833 - 0.2395456022787395im&1.55e-09&1.19e-01&1.79e-13 \\
		& 6 & 16.92320186444801 - 0.2395456022787569im & 1.55e-09 & 1.19e-01 & \\
		\midrule
		11&6&16.923201860726405 - 0.2395455899902375im&2.69e-11&1.19e-01&1.54e-13 \\
		& 6 & 16.923201860726255 - 0.23954558999018935im & 2.69e-11 & 1.19e-01 & \\
		\midrule
		12&5&16.92320391421663 - 0.2395457509419893im&2.46e-07&1.19e-01&4.91e-14 \\
		& 5 & 16.923203914216863 - 0.23954575094203084im & 2.46e-07 & 1.19e-01 & \\
		\midrule
		13&5&16.92320189408275 - 0.23954563529149417im&6.72e-09&1.19e-01&2.25e-14 \\
		& 5 & 16.923201894082876 - 0.2395456352915151im & 6.72e-09 & 1.19e-01 & \\
		\midrule
		14&5&16.92320186069578 - 0.2395455900546274im&3.52e-11&1.19e-01&3.42e-14 \\
		& 5 & 16.92320186069587 - 0.2395455900546616im & 3.52e-11 & 1.19e-01 & \\
		\midrule
		15&4&16.92320227549268 - 0.2395457990936466im&5.54e-08&1.19e-01&2.29e-15 \\
		& 4 & 16.923202275492653 - 0.23954579909366866im & 5.54e-08 & 1.19e-01 & \\
		\midrule
		16&4&16.923201860781877 - 0.2395455898294696im&1.07e-11&1.19e-01&4.43e-14 \\
		& 4 & 16.92320186078177 - 0.2395455898294269im & 1.06e-11 & 1.19e-01 & \\
		\midrule
		17&3&16.923201862834166 - 0.23954559087875893im&2.67e-10&1.19e-01&1.16e-14 \\
		& 3 & 16.923201862834055 - 0.23954559087876123im & 2.67e-10 & 1.19e-01 & \\
		\midrule
		18&4&16.923200862599433 - 0.2395464356945639im&1.56e-07&1.19e-01&1.20e-14 \\
		& 4 & 16.923200862599412 - 0.23954643569455172im & 1.56e-07 & 1.19e-01 & \\
		\midrule
		19&1000&-16.717608272660037+ 0.7591717411433029im&1.24e-01&1.30e-01&1.07e-13 \\
		& 1000 & -16.717608272660062 + 0.7591717411434057im & 1.24e-01 & 1.30e-01 & \\
		\midrule
		20&6&22.119804063752024 - 0.7063456283412871im&3.82e-09&6.03e-02&6.21e-15 \\
		& 6 & 22.119804063752017 - 0.7063456283412931im & 3.82e-09 & 6.03e-02 & \\
		\bottomrule
	\end{tabular}
	\caption{Comparison between the Newton method using ApproxFun (down) in the piecewise constant case and the Newton method using the exact expression~\eqref{eq:determinantpwconstant} of $D(k)$ %
		(up) for $\xi=0.5$, $m=10$, $n_1=1.5$ and $n_2=1$. We took here $\varepsilon=1e{-6}$ and $l_{\max} =1000$.}
	\label{tab:comaprison Newton}
\end{table}

\begin{table}[h!]
	\centering
	\begin{tabular}{l|ccccc}  
		\toprule
		$k_0$  & $l$ & $k$  & $|D(k)|$ & $|\partial_k D(k)|$ &  Difference in $k$\\
		\midrule
		\midrule
		21&4&22.119804057435196 - 0.7063456946543558im&3.00e-10&6.03e-02&1.26e-14 \\
		& 4 & 22.1198040574352 - 0.706345694654344im & 3.00e-10 & 6.03e-02 & \\
		\midrule
		22&3&22.11980251129046 - 0.7063450854527952im&1.00e-07&6.03e-02&7.10e-15 \\
		& 3 & 22.119802511290455 - 0.7063450854527893im & 1.00e-07 & 6.03e-02 & \\
		\midrule
		23&4&22.119804096508354 - 0.706345754755474im&4.35e-09&6.03e-02&5.12e-15 \\
		& 4 & 22.119804096508357 - 0.7063457547554691im & 4.35e-09 & 6.03e-02 & \\
		\midrule
		24&5&16.92320279370931 - 0.23954483086628545im&1.44e-07&1.19e-01&1.80e-13 \\
		& 5 & 16.923202793709482 - 0.2395448308663334im & 1.44e-07 & 1.19e-01 & \\
		\midrule
		25&5&31.730346846554976 - 0.9553141417588785im&7.90e-08&4.30e-02&2.70e-15 \\
		& 5 & 31.730346846554987 - 0.9553141417588782im & 7.90e-08 & 4.30e-02 & \\
		\midrule
		26&4&27.04249139796374 - 0.8848480038591856im&3.98e-07&4.93e-02&1.71e-14 \\
		& 4 & 27.042491397963737 - 0.8848480038592026im & 3.98e-07 & 4.93e-02 & \\
		\midrule
		27&4&27.042488355965958 - 0.8848405153357658im&1.31e-10&4.93e-02&1.26e-14 \\
		& 4 & 27.04248835596597 - 0.8848405153357701im & 1.31e-10 & 4.93e-02 & \\
		\midrule
		28&5&27.042488356712457 - 0.8848405174101011im&3.00e-11&4.93e-02&1.83e-14 \\
		& 5 & 27.042488356712454 - 0.8848405174101192im & 3.00e-11 & 4.93e-02 & \\
		\midrule
		29&7&16.92320259545549 - 0.23954508315910544im&1.06e-07&1.19e-01&1.96e-14 \\
		& 7 & 16.923202595455503 - 0.23954508315912526im & 1.06e-07 & 1.19e-01 & \\
		\midrule
		30&7&36.2794673729414 - 0.9907691069700367im&3.96e-10&3.83e-02&2.52e-14 \\
		& 7 & 36.279467372941376 - 0.9907691069700371im & 3.96e-10 & 3.83e-02 & \\
		\midrule
		31&4&31.73034262206993 - 0.9553087127872721im&2.33e-07&4.30e-02&2.08e-14 \\
		& 4 & 31.73034262206995 - 0.9553087127872671im & 2.33e-07 & 4.30e-02 & \\
		\midrule
		32&4&31.730344889489842 - 0.9553135682753766im&9.74e-09&4.30e-02&7.80e-15 \\
		& 4 & 31.730344889489842 - 0.9553135682753844im & 9.74e-09 & 4.30e-02 & \\
		\midrule
		33&6&31.73034327288489 - 0.9553121097293579im&1.00e-07&4.30e-02&2.11e-14 \\
		& 6 & 31.73034327288491 - 0.9553121097293374im & 1.00e-07 & 4.30e-02 & \\
		\midrule
		34&6&92.23783269919007 - 1.06289908087784im&4.07e-11&1.60e-02&1.13e-14 \\
		& 6 & 92.23783269919007 - 1.0628990808778287im & 4.07e-11 & 1.60e-02 & \\
		\midrule
		35&6&36.27946761966935 - 0.9907688668035165im&1.36e-08&3.83e-02&6.90e-15 \\
		& 6 & 36.27946761966934 - 0.9907688668035097im & 1.36e-08 & 3.83e-02 & \\
		\midrule
		36&4&36.27946660785861 - 0.9907693889504177im&3.09e-08&3.83e-02&7.10e-15 \\
		& 4 & 36.27946660785862 - 0.9907693889504107im & 3.09e-08 & 3.83e-02 & \\
		\midrule
		37&5&36.279467367290934 - 0.9907691165798495im&9.81e-11&3.83e-02&2.86e-14 \\
		& 5 & 36.279467367290906 - 0.9907691165798471im & 9.81e-11 & 3.83e-02 & \\
		\midrule
		38&5&31.730342433122015 - 0.9553140640209747im&1.18e-07&4.30e-02&5.70e-15 \\
		& 5 & 31.73034243312201 - 0.9553140640209804im & 1.18e-07 & 4.30e-02 & \\
		\midrule
		39&5&45.14774543822885 - 1.0249856558400918im&1.86e-10&3.16e-02&4.96e-14 \\
		& 5 & 45.147745438228895 - 1.0249856558400874im & 1.86e-10 & 3.16e-02 & \\
		\midrule
		40&5&40.74237122193321 - 1.0115554177751054im&1.58e-10&3.46e-02&2.50e-14 \\
		& 5 & 40.74237122193319 - 1.0115554177751027im & 1.58e-10 & 3.46e-02 & \\
		\bottomrule
	\end{tabular}
	\caption{(Cont'd) Comparison between the Newton method using ApproxFun (down) in the piecewise constant case and the Newton method using the exact expression~\eqref{eq:determinantpwconstant} of $D(k)$ 
		(up) for $\xi=0.5$, $m=10$, $n_1=1.5$ and $n_2=1$. We took here $\varepsilon=1e{-6}$ and $l_{\max} =1000$.
		.}
	\label{tab:comaprison Newton cont}
\end{table}

\clearpage

\begin{remark}
	For the validation (\S\ref{sec:val1}-\ref{sec:val2}), the Newton algorithm was applied on $D(k)=\frac{\det(k)}k$ and stopping criterion $|D(k)|\le \varepsilon$.
\end{remark}

\begin{remark}
	Investigating the iterations for the case where $k_0=19$, one finds that the process terminates by oscillating between two values $-16.71760827266006 + 0.7591717411434071 i$ and $-17.13411129251172 - 0.09502475320676629 i$. It is worth noting that, if one works with $\det$ instead of $D$, this case converges.
\end{remark}

\begin{remark}
	For information, the behavior of $\frac{\partial_k D(k)}{D(k)}$ is shown in Figure~\ref{fig:complexplotD} to be compared with Figure~\ref{fig:complexplotdetscal}.
	
	\begin{figure}[h]
		\centering
		\centering
		\begin{subfigure}[b]{0.35\textwidth}
			\centering
			\includegraphics[width=\textwidth]{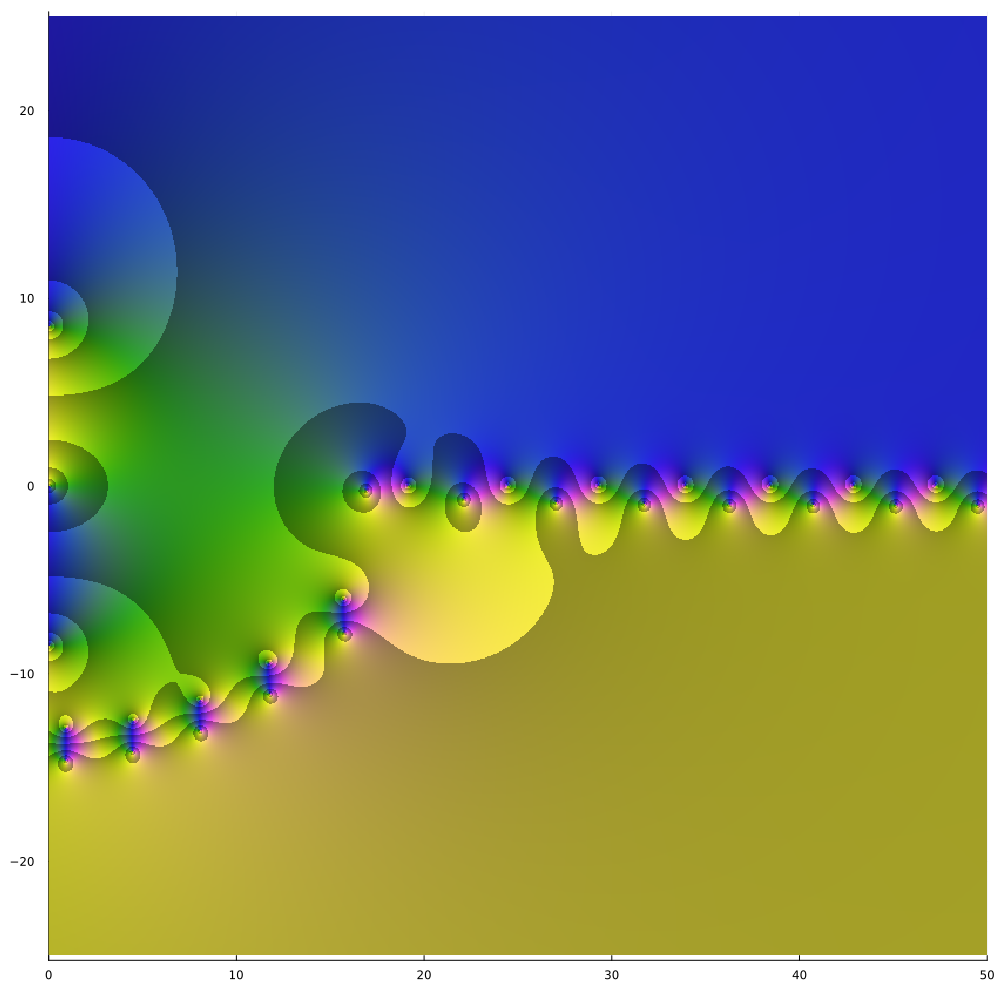}
			\caption{$\frac{\partial_k D(k)}{D(k)}$}
			\label{fig:derD}
		\end{subfigure}\qquad \qquad
		\\
		\begin{subfigure}[b]{0.4\textwidth}
			\centering
			\includegraphics[width=\textwidth]{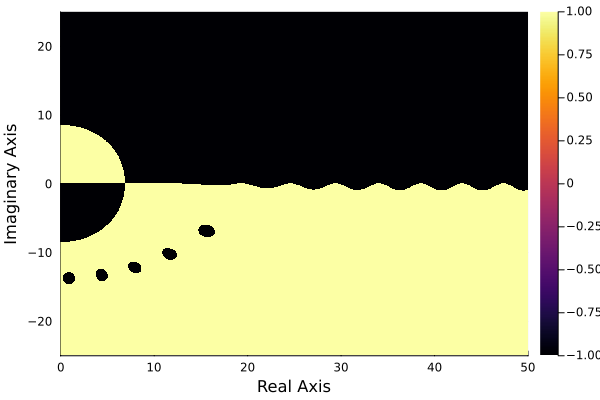}
			\caption{$\sign(\Im(\frac{\partial_k D(k)}{D(k)}))$}
			\label{fig:signderDratio}
		\end{subfigure}
		\caption{Complex plots of $\frac{\partial_k D(k)}{D(k)}$ where $D(k)=\frac{\det_1(k)}{k}$~\eqref{eq:determinantpwconstant} for $\xi=0.5$, $n_1=1.5$, $n_2=1$ and $m=10$.}
		\label{fig:complexplotD}
	\end{figure}
\end{remark}

\clearpage
\section{Variable case validation: A special variable case with explicit solution (Luneburg lens)}\label{sec:val2}

\begin{table}[h!]
	\centering
\begin{tabular}{l|ccccc}  
		\toprule
		$k_0$  & $l$ & $k$  & $|D(k)|$ & $|\partial_k D(k)|$ &  Difference in $k$\\
		\midrule
		\midrule
		1&9&18.58896340188007 - 0.6154425719496848im&2.74e-09&6.84e-02&1.75e-14 \\
		& 9 & 18.588963401880036 - 0.615442571949669im & 2.74e-09 & 6.84e-02 & \\
		\midrule
		2&8&18.58896340488882 - 0.6154425690916594im&2.51e-09&6.84e-02&3.46e-15 \\
		& 8 & 18.588963404888805 - 0.6154425690916552im & 2.51e-09 & 6.84e-02 & \\
		\midrule
		3&7&18.58896695779626 - 0.6154417792268551im&2.47e-07&6.84e-02&1.14e-14 \\
		& 7 & 18.588966957796284 - 0.6154417792268659im & 2.47e-07 & 6.84e-02 & \\
		\midrule
		4&7&18.58896341614563 - 0.6154425612670845im&1.74e-09&6.84e-02&7.09e-15 \\
\bottomrule
	\end{tabular}
	\caption{Comparison between the Newton method using ApproxFun (down) and the Newton method using the exact expression of $D(k)$ (up) in a special variable case $n_1$ given by $n_1(r)=\sqrt{2-r^2}$ 
		for $n_2=1$, $\xi=0.5$ and $m=10$. We took here $\varepsilon=1e{-6}$ and $l_{\max} =1000$.}
	\label{tab:NewtonLaguerre}
\end{table}

\begin{table}[h!]
	\centering
	\begin{tabular}{l|ccccc}  
		\toprule
		$k_0$  & $l$ & $k$  & $|D(k)|$ & $|\partial_k D(k)|$ &  Difference in $k$\\
		\midrule
		\midrule
		21&8&24.31472259294407 - 1.132533380322222im&3.20e-09&4.37e-02&9.86e-14 \\
		& 8 & 24.31472259294416 - 1.1325333803221875im & 3.20e-09 & 4.37e-02 &  \\
		\midrule
		22&5&24.31472055171216 - 1.1325342728376im&9.49e-08&4.37e-02&2.00e-13 \\
		& 5 & 24.314720551712146 - 1.1325342728377994im & 9.49e-08 & 4.37e-02 &  \\
		\midrule
		23&4&24.31472272287542 - 1.132533490562084im&1.05e-08&4.37e-02&5.38e-13 \\
		& 4 & 24.314722722875874 - 1.132533490562372im & 1.05e-08 & 4.37e-02 &  \\
		\midrule
		24&4&24.31472252029847 - 1.132533369059882im&4.79e-10&4.37e-02&6.79e-13 \\
\bottomrule
	\end{tabular}
	\caption{(Cont'd) Comparison between the Newton method using ApproxFun (down) and the Newton method using the exact expression of $D(k)$ (up) in a special variable case $n_1$ given by $n_1(r)=\sqrt{2-r^2}$  
		for $n_2=1$, $\xi=0.5$ and $m=10$. We took here $\varepsilon=1e{-6}$ and $l_{\max} =1000$.}
	\label{tab:NewtonLaguerre Cont}
\end{table}

\clearpage
\section{Piecewise constant cases}

\subsection{Low-contrast media}

\begin{table}[h!]
	\begin{equation*}
		n_1(r)=1.5
	\end{equation*}
	\footnotesize
	\centering
	\begin{tabular}{ll|ccccccc}  
		\toprule
		$m$  & $k_0$ & $l$ & $k$  & $k_\mathrm{asympt}$ & $|\det(k)|$ & $|\partial_k \det(k)|$  & $\frac{|\det(k)|}{\|T(k)\|_F}$ & $|k-k_\mathrm{asympt}|$  \\
		\midrule
		\midrule
		1 & 1.33 & 10 & 2.9412214517988846 - 1.0213193155493827im & 1.404547952561922 & 1.07e-14 & 1.40e+00 & 3.15e-15 & 1.85e+00 \\ 
		2 & 2.67 & 6 & 4.567191225356701 - 0.9043901168618548im & 3.8338779784797286 & 3.35e-11 & 1.32e+00 & 1.05e-11 & 1.16e+00 \\ 
		3 & 4 & 5 & 6.153857666206802 - 0.7885737875639456im & 5.777332794430966 & 2.85e-10 & 1.32e+00 & 9.29e-11 & 8.74e-01 \\ 
		4 & 5.33 & 5 & 7.720812756597912 - 0.6808760720924771im & 7.549595817571761 & 3.49e-14 & 1.35e+00 & 1.13e-14 & 7.02e-01 \\ 
		5 & 6.67 & 4 & 9.275045821677303 - 0.5831461096495909im & 9.233276944106649 & 1.64e-08 & 1.41e+00 & 5.07e-09 & 5.85e-01 \\ 
		6 & 8 & 4 & 10.819644524989346 - 0.49566872221472635im & 10.862185118573025 & 6.08e-09 & 1.49e+00 & 1.72e-09 & 4.97e-01 \\ 
		7 & 9.33 & 4 & 12.356132077506965 - 0.41815032448175377im & 12.45353017076617 & 2.37e-08 & 1.59e+00 & 5.93e-09 & 4.29e-01 \\
\bottomrule
	\end{tabular}
	\caption{Results obtained by the Newton method for $n_1=1.5$, $n_2=1$, $\xi=0.5$ starting from $k_0=\frac{m}{\xi n_0}$. We took here $\varepsilon=1e{-8}$ and $l_{\max} =2000$.}
	\label{tab:Newtonconstant1.5}
\end{table}

\clearpage
\subsection{High-contrast media}

\begin{table}[h!]
	\begin{equation*}
		n_1(r)=5
	\end{equation*}
	\footnotesize
	\centering
	\begin{tabular}{ll|ccccccc}  
		\toprule
		$m$  & $k_0$ & $l$ & $k$  & $k_\mathrm{asympt}$ & $|\det(k)|$ & $|\partial_k \det(k)|$  & $\frac{|\det(k)|}{\|T(k)\|_F}$ & $|k-k_\mathrm{asympt}|$  \\
		\midrule
		\midrule
		1 & 0.4 & 5 & 0.9289987680826443 - 0.04523834315971972im & 0.8842557552463222 & 1.44e-15 & 9.62e+00 & 4.39e-16 & 6.36e-02 \\ 
		2 & 0.8 & 4 & 1.4998099150133781 - 0.010215157356845911im & 1.4892827214976112 & 1.19e-11 & 1.99e+01 & 1.25e-12 & 1.47e-02 \\ 
		3 & 1.2 & 5 & 2.0304915035306306 - 0.0016737245701525074im & 2.0224106828300403 & 1.90e-11 & 4.78e+01 & 6.16e-13 & 8.25e-03 \\ 
		4 & 1.6 & 5 & 2.5332600608053313 - 0.0002250224324238292im & 2.5260265729729126 & 2.37e-07 & 1.27e+02 & 2.32e-09 & 7.24e-03 \\ 
		5 & 2 & 6 & 3.0190057464666977 - 2.6964505481422285e-5im & 3.0127828023877194 & 1.39e-11 & 3.60e+02 & 4.06e-14 & 6.22e-03 \\ 
		6 & 2.4 & 6 & 3.4936803250165007 - 3.004752173541452e-6im & 3.4883629888379457 & 1.34e-08 & 1.06e+03 & 1.17e-11 & 5.32e-03 \\ 
		7 & 2.8 & 6 & 3.960459416225663 - 3.183017415564268e-7im & 3.9558733069976633 & 2.31e-06 & 3.21e+03 & 5.96e-10 & 4.59e-03 \\
\bottomrule
	\end{tabular}
	\caption{Results obtained by the Newton method for $n_1=5$, $n_2=1$, $\xi=0.5$ starting from $k_0=\frac{m}{\xi n_0}$. We took here $\varepsilon=1e{-8}$ and $l_{\max} =2000$.
	}
	\label{tab:Newtonconstant5}
\end{table}

\clearpage
\section{Variable cases}

\subsection{$n_1$ variable and $n_2=1$}

\begin{table}[h!]
	\begin{equation*}
		n_1(r)=2-r
	\end{equation*}
	\footnotesize
	\centering
	\begin{tabular}{ll|ccccccc}  
		\toprule
		$m$  & $k_0$ & $l$ & $k$  & $k_\mathrm{asympt}$ & $|\det(k)|$ & $|\partial_k \det(k)|$  & $\frac{|\det(k)|}{\|T(k)\|_F}$ & $|k-k_\mathrm{asympt}|$  \\
		\midrule
		\midrule
		1 & 1.33 & 13 & 6.319821637534675 - 0.8852775094525012im & 1.4051097872513703 & 3.67e-08 & 1.40e+00 & 7.33e-09 & 4.99e+00 \\ 
		2 & 2.67 & 7 & 4.311909988078464 - 0.6755610763676703im & 3.708252723915477 & 4.23e-14 & 1.69e+00 & 1.41e-14 & 9.06e-01 \\ 
		3 & 4 & 6 & 5.847774226712011 - 0.5742451453806543im & 5.553028501041968 & 3.01e-15 & 1.72e+00 & 1.00e-15 & 6.45e-01 \\ 
		4 & 5.33 & 5 & 7.369474014406212 - 0.4824173452925744im & 7.243859029864192 & 5.15e-11 & 1.81e+00 & 1.63e-11 & 4.99e-01 \\ 
		5 & 6.67 & 5 & 8.881738598735172 - 0.4012506982673458im & 8.857569745641534 & 3.92e-14 & 1.94e+00 & 1.12e-14 & 4.02e-01 \\ 
		6 & 8 & 4 & 10.386432418688315 - 0.3305887598822313im & 10.424704364080949 & 2.80e-08 & 2.10e+00 & 7.01e-09 & 3.33e-01 \\ 
		7 & 9.33 & 4 & 11.884338704967615 - 0.26981120308410683im & 11.96046139230838 & 6.11e-09 & 2.31e+00 & 1.30e-09 & 2.80e-01 \\
\bottomrule
	\end{tabular}
	\caption{Results obtained by the Newton method for $n_1(r)=2-r$, $n_2=1$, $\xi=0.5$ starting from $k_0=\frac{m}{\xi n_0}$. We took here $\varepsilon=1e{-8}$ and $l_{\max} =2000$.}
	\label{tab:Newtonnonconstantaffine1}
\end{table}

\begin{table}[h!]
	\begin{equation*}
		n_1(r)=1.5+r
	\end{equation*}
	\footnotesize
	\centering
	\begin{tabular}{ll|ccccccc}  
		\toprule
		$m$  & $k_0$ & $l$ & $k$  & $k_\mathrm{asympt}$ & $|\det(k)|$ & $|\partial_k \det(k)|$  & $\frac{|\det(k)|}{\|T(k)\|_F}$ & $|k-k_\mathrm{asympt}|$  \\
		\midrule
		\midrule
		1 & 1 & 6 & 2.3573676407577873 - 0.5820419051788126im & 1.8200470263362465 & 1.57e-15 & 1.71e+00 & 6.10e-16 & 7.92e-01 \\
		2 & 2 & 4 & 3.676563387273686 - 0.43451946674734865im & 3.542548925725047 & 2.25e-08 & 1.81e+00 & 8.74e-09 & 4.55e-01 \\
		3 & 3 & 4 & 4.9645216831943495 - 0.3132395128737777im & 4.979499565390129 & 2.58e-10 & 2.03e+00 & 8.57e-11 & 3.14e-01 \\
		4 & 4 & 5 & 6.2318202354950385 - 0.21855409808217305im & 6.308316746647625 & 1.42e-14 & 2.35e+00 & 3.59e-15 & 2.32e-01 \\
		5 & 5 & 5 & 7.480957938826707 - 0.1474128807045029im & 7.578289435972096 & 1.40e-10 & 2.78e+00 & 2.53e-11 & 1.77e-01 \\
		6 & 6 & 5 & 8.712536676603357 - 0.09598428277781516im & 8.810501058485475 & 5.20e-08 & 3.36e+00 & 6.64e-09 & 1.37e-01 \\
		7 & 7 & 6 & 9.926896921138045 - 0.06030041686723931im & 10.016070612813092 & 1.25e-12 & 4.12e+00 & 1.11e-13 & 1.08e-01 \\
\bottomrule
	\end{tabular}
	\caption{Results obtained by the Newton method for $n_1(r)=1.5+r$, $n_2=1$, $\xi=0.5$ starting from $k_0=\frac{m}{\xi n_0}$. We took here $\varepsilon=1e{-8}$ and $l_{\max} =2000$.}
	\label{tab:Newtonnonconstantaffine2}
\end{table}

\begin{table}[h!]
	\begin{equation*}
		n_1(r)=1+r
	\end{equation*}
	\footnotesize
	\centering
	\begin{tabular}{ll|ccccccc}  
		\toprule
		$m$  & $k_0$ & $l$ & $k$  & $k_\mathrm{asympt}$ & $|\det(k)|$ & $|\partial_k \det(k)|$  & $\frac{|\det(k)|}{\|T(k)\|_F}$ & $|k-k_\mathrm{asympt}|$  \\
		\midrule
		\midrule
		1 & 1.33 & 9 & 13.602746351567303 - 1.2941451159878599im & 0.97222155979372 & 1.50e-11 & 1.09e+00 & 1.73e-12 & 1.27e+01 \\ 
		2 & 2.67 & 6 & 4.800842869405789 - 1.183070967313181im & 3.739093521053414 & 3.75e-11 & 1.05e+00 & 1.08e-11 & 1.59e+00 \\ 
		3 & 4 & 5 & 6.435809549914782 - 1.0457273636322453im & 5.852012676075551 & 1.46e-10 & 1.02e+00 & 4.53e-11 & 1.20e+00 \\ 
		4 & 5.33 & 5 & 8.044440569723108 - 0.918011144092654im & 7.73893601066758 & 6.66e-14 & 1.02e+00 & 2.15e-14 & 9.68e-01 \\ 
		5 & 6.67 & 5 & 9.636426835812768 - 0.8011294949004661im & 9.510506071020615 & 1.11e-14 & 1.04e+00 & 3.55e-15 & 8.11e-01 \\ 
		6 & 8 & 5 & 11.216277798909314 - 0.6951459775408366im & 11.211474721375804 & 2.64e-14 & 1.07e+00 & 8.12e-15 & 6.95e-01 \\ 
		7 & 9.33 & 5 & 12.786368558462707 - 0.5996993138799153im & 12.864406898379844 & 1.56e-12 & 1.12e+00 & 4.44e-13 & 6.05e-01 \\
\bottomrule
	\end{tabular}
	\caption{Results obtained by the Newton method for $n_1(r)=1+r$, $n_2=1$, $\xi=0.5$ starting from $k_0=\frac{m}{\xi n_0}$. We took here $\varepsilon=1e{-8}$ and $l_{\max} =2000$.}
	\label{tab:Newtonnonconstantaffine3}
\end{table}

\begin{table}[h!]
	\begin{equation*}
		n_1(r)=3(1-r)
	\end{equation*}
	\footnotesize
	\centering
	\begin{tabular}{ll|ccccccc}  
		\toprule
		$m$  & $k_0$ & $l$ & $k$  & $k_\mathrm{asympt}$ & $|\det(k)|$ & $|\partial_k \det(k)|$  & $\frac{|\det(k)|}{\|T(k)\|_F}$ & $|k-k_\mathrm{asympt}|$  \\
		\midrule
		\midrule
		1 & 1.33 & 9 & 5.0550018091851285 - 0.6248886757558838im & 4.161760458079524 & 8.68e-11 & 2.43e+00 & 1.91e-11 & 1.09e+00 \\ 
		2 & 2.67 & 9 & 3.793068555579194 - 0.3612720296024689im & 5.495093791412857 & 1.43e-10 & 2.73e+00 & 5.06e-11 & 1.74e+00 \\ 
		3 & 4 & 7 & 5.211463134482563 - 0.2775990463435686im & 6.82842712474619 & 3.45e-11 & 2.96e+00 & 1.11e-11 & 1.64e+00 \\ 
		4 & 5.33 & 6 & 6.624138531594931 - 0.20983062483433862im & 8.161760458079524 & 1.13e-08 & 3.33e+00 & 3.09e-09 & 1.55e+00 \\ 
		5 & 6.67 & 6 & 8.031924435987753 - 0.15618584407704272im & 9.495093791412858 & 1.32e-10 & 3.82e+00 & 2.94e-11 & 1.47e+00 \\ 
		6 & 8 & 6 & 9.43467283612071 - 0.11452520508310852im & 10.82842712474619 & 7.23e-12 & 4.46e+00 & 1.28e-12 & 1.40e+00 \\ 
		7 & 9.33 & 6 & 10.832115071451929 - 0.08276743632092393im & 12.161760458079524 & 1.03e-12 & 5.27e+00 & 1.43e-13 & 1.33e+00 \\
\bottomrule
	\end{tabular}
	\caption{Results obtained by the Newton method for $n_1(r)=3(1-r)$, $n_2=1$, $\xi=0.5$ starting from $k_0=\frac{m}{\xi n_0}$. We took here $\varepsilon=1e{-8}$ and $l_{\max} =2000$.}
	\label{tab:Newtonnonconstantaffine4}
\end{table}

\begin{table}[h!]
	\begin{equation*}
		n_1(r)=-2.8 r+2.5
	\end{equation*}
	\footnotesize
	\centering
	\begin{tabular}{ll|ccccccc}  
		\toprule
		$m$  & $k_0$ & $l$ & $k$  & $k_\mathrm{asympt}$ & $|\det(k)|$ & $|\partial_k \det(k)|$  & $\frac{|\det(k)|}{\|T(k)\|_F}$ & $|k-k_\mathrm{asympt}|$  \\
		\midrule
		\midrule
		1 & 1.82 & 17 & 7.924919667036129 - 2.7112347145700677im & 3.059135351886293 & 5.18e-08 & 7.76e+00 & 4.36e-09 & 5.57e+00 \\ 
		2 & 3.64 & 7 & 8.37697526088352 - 1.1970651616679695im & 4.851135351886293 & 9.84e-09 & 2.13e+00 & 1.11e-09 & 3.72e+00 \\ 
		3 & 5.45 & 7 & 10.164172357599108 - 1.1715053777822293im & 6.643135351886293 & 3.11e-14 & 1.85e+00 & 3.69e-15 & 3.71e+00 \\ 
		4 & 7.27 & 7 & 11.945784263931984 - 1.1385426036460573im & 8.435135351886293 & 5.35e-14 & 1.71e+00 & 6.59e-15 & 3.69e+00 \\ 
		5 & 9.09 & 7 & 13.724183225760296 - 1.1022390897528125im & 10.227135351886291 & 8.15e-14 & 1.64e+00 & 1.03e-14 & 3.67e+00 \\ 
		6 & 10.91 & 7 & 15.500792465870116 - 1.0645281109812703im & 12.019135351886291 & 2.05e-12 & 1.61e+00 & 2.65e-13 & 3.64e+00 \\ 
		7 & 12.73 & 7 & 17.276486339139304 - 1.0264260048504945im & 13.811135351886294 & 1.31e-10 & 1.60e+00 & 1.73e-11 & 3.61e+00 \\
\bottomrule
	\end{tabular}
	\caption{Results obtained by the Newton method for $n_1(r)=-2.8r+2.5$, $n_2=1$, $\xi=0.5$ starting from $k_0=\frac{m}{\xi n_0}$. We took here $\varepsilon=1e{-8}$ and $l_{\max} =2000$.}
	\label{tab:Newtonnonconstantaffine5}
\end{table}

\begin{table}[h!]
	\begin{equation*}
		n_1(r)=-2.8 r+2.5
	\end{equation*}
	\footnotesize
	\centering
	\begin{tabular}{ll|ccccccc}  
		\toprule
		$m$  & $k_0$ & $l$ & $k$  & $k_\mathrm{asympt}$ & $|\det(k)|$ & $|\partial_k \det(k)|$  & $\frac{|\det(k)|}{\|T(k)\|_F}$ & $|k-k_\mathrm{asympt}|$  \\
		\midrule
		\midrule
		1 & 1.79 & 19 & 7.924919673377272 - 2.711234712507853im & 3.059135351886293 & 1.15e-12 & 7.76e+00 & 9.73e-14 & 5.57e+00 \\ 
		2 & 3.58 & 8 & 8.376975264130033 - 1.1970651649564188im & 4.851135351886293 & 4.52e-12 & 2.13e+00 & 5.08e-13 & 3.72e+00 \\ 
		3 & 5.38 & 6 & 10.164172393807487 - 1.1715053821567618im & 6.643135351886293 & 6.75e-08 & 1.85e+00 & 8.01e-09 & 3.71e+00 \\ 
		4 & 7.17 & 6 & 11.94578427594684 - 1.1385426064549284im & 8.435135351886293 & 2.11e-08 & 1.71e+00 & 2.60e-09 & 3.69e+00 \\ 
		5 & 8.96 & 7 & 13.724183225760338 - 1.1022390897528427im & 10.227135351886291 & 5.62e-15 & 1.64e+00 & 7.10e-16 & 3.67e+00 \\ 
		6 & 10.75 & 7 & 15.500792465869143 - 1.0645281109805322im & 12.019135351886291 & 9.42e-14 & 1.61e+00 & 1.22e-14 & 3.64e+00 \\ 
		7 & 12.54 & 7 & 17.276486339104792 - 1.0264260047759763im & 13.811135351886294 & 1.93e-12 & 1.60e+00 & 2.55e-13 & 3.61e+00 \\
\bottomrule
	\end{tabular}
	\caption{Results obtained by the Newton method for $n_1(r)=-2.8r+2.5$, $n_2=1$, $\xi=0.5$ starting from $k_0=\frac{m}{\xi_0 n(\xi_0)}$ (leading term in the asymptotics). We took here $\varepsilon=1e{-8}$ and $l_{\max} =2000$.}
	\label{tab:Newtonnonconstantaffine5startleadingterm}
\end{table}

\begin{table}[h!]
	\begin{equation*}
		n_1(r)=-2.8 r+2.5
	\end{equation*}
	\footnotesize
	\centering
	\begin{tabular}{ll|ccccccc}  
		\toprule
		$m$  & $k_0$ & $l$ & $k$  & $k_\mathrm{asympt}$ & $|\det(k)|$ & $|\partial_k \det(k)|$  & $\frac{|\det(k)|}{\|T(k)\|_F}$ & $|k-k_\mathrm{asympt}|$  \\
		\midrule
		\midrule
		1 & 3.06 & 8 & 3.126030006549902 - 0.8490620306257869im & 3.059135351886293 & 9.45e-09 & 2.32e+00 & 2.24e-09 & 8.52e-01 \\ 
		2 & 4.85 & 8 & 4.960265572425262 - 0.7930888853914064im & 4.851135351886293 & 2.58e-15 & 2.21e+00 & 6.01e-16 & 8.01e-01 \\ 
		3 & 6.64 & 7 & 6.770901080871078 - 0.7367987701752973im & 6.643135351886293 & 1.14e-14 & 2.05e+00 & 2.72e-15 & 7.48e-01 \\ 
		4 & 8.44 & 6 & 8.572070703674939 - 0.6839781193220764im & 8.435135351886293 & 1.98e-11 & 1.97e+00 & 4.79e-12 & 6.98e-01 \\ 
		5 & 10.23 & 6 & 10.36876216651674 - 0.6350860184782311im & 10.227135351886291 & 1.55e-14 & 1.92e+00 & 3.78e-15 & 6.51e-01 \\ 
		6 & 12.02 & 5 & 12.16318428414382 - 0.5899319796825272im & 12.019135351886291 & 1.49e-08 & 1.91e+00 & 3.63e-09 & 6.07e-01 \\ 
		7 & 13.81 & 5 & 13.956444792753562 - 0.5481968118098185im & 13.811135351886294 & 1.27e-09 & 1.91e+00 & 3.07e-10 & 5.67e-01 \\
\bottomrule
	\end{tabular}
	\caption{Results obtained by the Newton method for $n_1(r)=-2.8r+2.5$, $n_2=1$, $\xi=0.5$ starting from $k_\mathrm{asympt}$. We took here $\varepsilon=1e{-8}$ and $l_{\max} =2000$.}
	\label{tab:Newtonnonconstantaffine5startkth}
\end{table}

\begin{table}[h!]
	\begin{equation*}
		n_1(r)=1.5+6r(\xi-r)
	\end{equation*}
	\footnotesize
	\centering
	\begin{tabular}{ll|ccccccc}  
		\toprule
		$m$  & $k_0$ & $l$ & $k$  & $k_\mathrm{asympt}$ & $|\det(k)|$ & $|\partial_k \det(k)|$  & $\frac{|\det(k)|}{\|T(k)\|_F}$ & $|k-k_\mathrm{asympt}|$  \\
		\midrule
		\midrule
		1 & 1.33 & 15 & 6.2127404455770465 - 0.8387720532805358im & 5.333333333333333 & 4.74e-08 & 1.54e+00 & 9.63e-09 & 1.22e+00 \\ 
		2 & 2.67 & 7 & 4.047468166574287 - 0.5246463375796181im & 6.666666666666666 & 2.97e-09 & 2.06e+00 & 1.03e-09 & 2.67e+00 \\ 
		3 & 4 & 6 & 5.499773898985259 - 0.41551950976663776im & 8.0 & 4.90e-11 & 2.20e+00 & 1.63e-11 & 2.53e+00 \\ 
		4 & 5.33 & 6 & 6.942042250192962 - 0.3249938617322503im & 9.333333333333332 & 1.06e-14 & 2.42e+00 & 3.15e-15 & 2.41e+00 \\ 
		5 & 6.67 & 5 & 8.377083767868388 - 0.25104186154205443im & 10.666666666666668 & 1.58e-09 & 2.71e+00 & 3.99e-10 & 2.30e+00 \\ 
		6 & 8 & 5 & 9.805651396895824 - 0.1914390997765093im & 12.0 & 6.90e-11 & 3.10e+00 & 1.42e-11 & 2.20e+00 \\ 
		7 & 9.33 & 5 & 11.227871705396403 - 0.14407559550030327im & 13.333333333333334 & 5.52e-12 & 3.57e+00 & 9.16e-13 & 2.11e+00 \\
\bottomrule
	\end{tabular}
	\caption{Results obtained by the Newton method for $n_1(r)=1.5+6r(\xi-r)$, $n_2=1$, $\xi=0.5$ starting from $k_0=\frac{m}{\xi n_0}$. We took here $\varepsilon=1e{-8}$ and $l_{\max} =2000$.}
	\label{tab:Newtonnonconstantparabolic1}
\end{table}

\begin{table}[h!]
	\begin{equation*}
		n_1(r)=1.5-6r(\xi-r)
	\end{equation*}
	\footnotesize
	\centering
	\begin{tabular}{ll|ccccccc}  
		\toprule
		$m$  & $k_0$ & $l$ & $k$  & $k_\mathrm{asympt}$ & $|\det(k)|$ & $|\partial_k \det(k)|$  & $\frac{|\det(k)|}{\|T(k)\|_F}$ & $|k-k_\mathrm{asympt}|$  \\
		\midrule
		\midrule
		1 & 1.33 & 8 & 5.81160775798412 - 1.9882387793184597im & 0.8013387620671644 & 1.80e-08 & 1.77e+00 & 2.61e-09 & 5.39e+00 \\
		2 & 2.67 & 7 & 5.062619436026797 - 1.553700758832008im & 3.797740929268689 & 7.67e-15 & 8.14e-01 & 1.93e-15 & 2.00e+00 \\
		3 & 4 & 5 & 6.77380298797884 - 1.4179046129451796im & 6.073077254855459 & 6.90e-09 & 7.44e-01 & 1.92e-09 & 1.58e+00 \\
		4 & 5.33 & 5 & 8.45020498175942 - 1.2821303643219im & 8.090067850271312 & 5.62e-11 & 7.03e-01 & 1.69e-11 & 1.33e+00 \\
		5 & 6.67 & 5 & 10.104099467260863 - 1.1519395713087264im & 9.971935246328307 & 1.46e-10 & 6.78e-01 & 4.63e-11 & 1.16e+00 \\
		6 & 8 & 5 & 11.741692047461065 - 1.0294931116728006im & 11.769619297419684 & 3.13e-09 & 6.64e-01 & 1.01e-09 & 1.03e+00 \\
		7 & 9.33 & 6 & 13.366517237596524 - 0.9155941154546037im & 13.509248610554106 & 8.05e-15 & 6.58e-01 & 2.55e-15 & 9.27e-01 \\
\bottomrule
	\end{tabular}
	\caption{Results obtained by the Newton method for $n_1(r)=1.5-6r(\xi-r)$, $n_2=1$, $\xi=0.5$ starting from $k_0=\frac{m}{\xi n_0}$. We took here $\varepsilon=1e{-8}$ and $l_{\max} =2000$.
	}
	\label{tab:Newtonnonconstantparabolic2}
\end{table}

\begin{table}[h!]
	\begin{equation*}
		n_1(r)=1.5-6r(\xi-r)
	\end{equation*}
	\footnotesize
	\centering
	\begin{tabular}{ll|ccccccc}  
		\toprule
		$m$  & $k_0$ & $l$ & $k$  & $k_\mathrm{asympt}$ & $|\det(k)|$ & $|\partial_k \det(k)|$  & $\frac{|\det(k)|}{\|T(k)\|_F}$ & $|k-k_\mathrm{asympt}|$  \\
		\midrule
		\midrule
		1 & 0.80 & 7 & 5.811607768328178 - 1.988238794791164im & 0.8013387620671644 & 1.71e-08 & 1.77e+00 & 2.48e-09 & 5.39e+00 \\ 
		2 & 3.80 & 5 & 5.0626194360283705 - 1.553700758846383im & 3.797740929268689 & 1.18e-11 & 8.14e-01 & 2.96e-12 & 2.00e+00 \\ 
		3 & 6.07 & 5 & 6.773802979521228 - 1.4179046167704703im & 6.073077254855459 & 2.43e-14 & 7.44e-01 & 6.74e-15 & 1.58e+00 \\ 
		4 & 8.09 & 5 & 8.450204981711662 - 1.2821303643861277im & 8.090067850271312 & 1.52e-14 & 7.03e-01 & 4.58e-15 & 1.33e+00 \\ 
		5 & 9.97 & 5 & 10.104099467097347 - 1.1519395711676987im & 9.971935246328307 & 3.14e-15 & 6.78e-01 & 9.95e-16 & 1.16e+00 \\ 
		6 & 11.77 & 5 & 11.741692042775735 - 1.0294931112001742im & 11.769619297419684 & 8.57e-15 & 6.64e-01 & 2.76e-15 & 1.03e+00 \\ 
		7 & 13.51 & 5 & 13.366517237596518 - 0.9155941154545839im & 13.509248610554106 & 6.86e-15 & 6.58e-01 & 2.17e-15 & 9.27e-01 \\
\bottomrule
	\end{tabular}
	\caption{Results obtained by the Newton method for $n_1(r)=1.5-6r(\xi-r)$, $n_2=1$, $\xi=0.5$ starting from $k_0=k_\mathrm{asympt}$. We took here $\varepsilon=1e{-8}$ and $l_{\max} =2000$.
	}
	\label{tab:Newtonnonconstantparabolic2asympt}
\end{table}

\begin{table}[h!]
	\begin{equation*}
		n_1(r)=3-r(r+1)
	\end{equation*}
	\footnotesize
	\centering
	\begin{tabular}{ll|ccccccc}  
		\toprule
		$m$  & $k_0$ & $l$ & $k$  & $k_\mathrm{asympt}$ & $|\det(k)|$ & $|\partial_k \det(k)|$  & $\frac{|\det(k)|}{\|T(k)\|_F}$ & $|k-k_\mathrm{asympt}|$  \\
		\midrule
		\midrule
		1 & 0.89 & 6 & 1.8143997879339349 - 0.2594411374549328im & 1.3379979374850162 & 2.17e-10 & 3.21e+00 & 8.99e-11 & 5.42e-01 \\
		2 & 1.78 & 5 & 2.91191502137055 - 0.15714010829805797im & 2.7985740610105108 & 8.42e-12 & 3.95e+00 & 2.74e-12 & 1.94e-01 \\ 
		3 & 2.67 & 4 & 3.989187559993246 - 0.08804724456565523im & 3.9793015281921194 & 2.68e-08 & 5.21e+00 & 5.66e-09 & 8.86e-02 \\ 
		4 & 3.56 & 4 & 5.046457936436096 - 0.04590716693828803im & 5.069595700420508 & 8.50e-09 & 7.19e+00 & 1.09e-09 & 5.14e-02 \\ 
		5 & 4.44 & 5 & 6.083657425358343 - 0.022395492175469326im & 6.1155946033931174 & 7.50e-15 & 1.03e+01 & 5.74e-16 & 3.90e-02 \\ 
		6 & 5.33 & 5 & 7.10244321017077 - 0.01031427758860502im & 7.13522154082962 & 2.96e-12 & 1.52e+01 & 1.34e-13 & 3.44e-02 \\ 
		7 & 6.22 & 5 & 8.105641723913287 - 0.00453303804189889im & 8.137248806935698 & 2.46e-10 & 2.30e+01 & 6.63e-12 & 3.19e-02 \\
\bottomrule
	\end{tabular}
	\caption{Results obtained by the Newton method for $n_1(r)=3-r(r+1)$, $n_2=1$, $\xi=0.5$ starting from $k_0=\frac{m}{\xi n_0}$. We took here $\varepsilon=1e{-8}$ and $l_{\max} =2000$.}
	\label{tab:Newtonnonconstantparabolic3}
\end{table}

\begin{table}[h!]
	\begin{equation*}
		n_1(r)=\sqrt{2-r^2}
	\end{equation*}
	\footnotesize
	\centering
	\begin{tabular}{ll|ccccccc}  
		\toprule
		$m$  & $k_0$ & $l$ & $k$  & $k_\mathrm{asympt}$ & $|\det(k)|$ & $|\partial_k \det(k)|$  & $\frac{|\det(k)|}{\|T(k)\|_F}$ & $|k-k_\mathrm{asympt}|$  \\
		\midrule
		\midrule
		1 & 1.51 & 6 & 7.8876676918624335 - 1.4052428170656204im & 0.2779979601087338 & 2.81e-10 & 1.26e+00 & 3.85e-11 & 7.74e+00 \\ 
		2 & 3.02 & 14 & 5.043983011022969 - 1.255099671399734im & 3.3992834444609796 & 5.01e-15 & 1.21e+00 & 1.26e-15 & 2.07e+00 \\ 
		3 & 4.54 & 7 & 6.78117621695404 - 1.1657519246857233im & 5.725929985613505 & 6.55e-14 & 1.15e+00 & 1.72e-14 & 1.57e+00 \\ 
		4 & 6.05 & 6 & 8.495101178211378 - 1.07544263327515im & 7.79147286644759 & 4.59e-13 & 1.13e+00 & 1.26e-13 & 1.29e+00 \\ 
		5 & 7.56 & 5 & 10.195114750093254 - 0.9876143271520355im & 9.728823554860602 & 2.67e-09 & 1.12e+00 & 7.47e-10 & 1.09e+00 \\ 
		6 & 9.07 & 5 & 11.885660825533378 - 0.9037283263467333im & 11.590106722069358 & 9.57e-12 & 1.13e+00 & 2.69e-12 & 9.51e-01 \\ 
		7 & 10.58 & 5 & 13.569145340854345 - 0.8244171609148544im & 13.400921462365059 & 5.13e-14 & 1.16e+00 & 1.42e-14 & 8.41e-01 \\
\bottomrule
	\end{tabular}
	\caption{Results obtained by the Newton method for $n_1(r)=\sqrt{2-r^2}$, $n_2=1$, $\xi=0.5$ starting from $k_0=\frac{m}{\xi n_0}$. We took here $\varepsilon=1e{-8}$ and $l_{\max} =2000$.}
	\label{tab:NewtonnonconstantLuneberg}
\end{table}

\clearpage
\subsection{$n_1$ and  $n_2$ variable}

\begin{table}[h!]
	\begin{equation*}
		n_1(r)=\sqrt{2-r^2} \qquad n_2(r)=r+ 0.5
	\end{equation*}
	\footnotesize
	\centering
	\begin{tabular}{ll|ccccc}  
		\toprule
		$m$  & $k_0$ & $l$ & $k$ & $|\det(k)|$ & $|\partial_k \det(k)|$  & $\frac{|\det(k)|}{\|T(k)\|_F}$  \\
		\midrule
		\midrule
		1 & 1.51 & 17 & 0.40053408329828327 - 1.2449374864234752im & 8.48e-09 & 1.30e+00 & 5.31e-09 \\ 
		2 & 3.02 & 15 & 11.44346347372721 - 2.286602295636471im & 3.34e-13 & 6.05e+00 & 1.80e-14 \\ 
		3 & 4.54 & 15 & -7.6810822010348785 - 1.774902301665455im & 1.10e-10 & 3.66e+00 & 1.48e-11 \\ 
		4 & 6.05 & 6 & 8.89369366419431 - 1.0592866230998579im & 3.32e-11 & 1.39e+00 & 7.37e-12 \\ 
		5 & 7.56 & 5 & 10.40252005138865 - 0.9746801306001142im & 6.80e-09 & 1.36e+00 & 1.69e-09 \\ 
		6 & 9.07 & 5 & 11.971622435387683 - 1.0228695592136103im & 1.53e-12 & 1.16e+00 & 3.84e-13 \\ 
		7 & 10.58 & 4 & 13.757647634726336 - 1.0978647143490887im & 3.83e-11 & 8.58e-01 & 8.39e-12 \\
\bottomrule
	\end{tabular}
	\caption{Results obtained by the Newton method for $n_1(r)=\sqrt{2-r^2}$, $n_2=r+ 0.5$, $\xi=0.5$ starting from $k_0=\frac{m}{\xi n_0}$. We took here $\varepsilon=1e{-8}$ and $l_{\max} =2000$.}
	\label{tab:Newtonnonconstantvarn21}
\end{table}

\begin{table}[h!]
	\begin{equation*}
		n_1(r)=\sqrt{2-r^2} \qquad n_2(r)=1+(r- 0.5)^3
	\end{equation*}
	\footnotesize
	\centering
	\begin{tabular}{ll|ccccc}  
		\toprule
		$m$  & $k_0$ & $l$ & $k$ & $|\det(k)|$ & $|\partial_k \det(k)|$  & $\frac{|\det(k)|}{\|T(k)\|_F}$  \\
		\midrule
		\midrule
		1 & 1.51 & 19 & 8.33731377407904 - 1.5145282232893382im & 8.63e-15 & 1.55e+00 & 1.00e-15 \\ 
		2 & 3.02 & 10 & 10.018989709571363 - 1.2044185189436711im & 9.92e-15 & 1.56e+00 & 1.46e-15 \\ 
		3 & 4.54 & 8 & 6.560084964097422 - 0.9823851559418392im & 1.05e-12 & 1.47e+00 & 3.47e-13 \\ 
		4 & 6.05 & 7 & 8.20434829959099 - 1.1276540761930103im & 2.47e-12 & 1.14e+00 & 7.43e-13 \\ 
		5 & 7.56 & 6 & 10.30091455556008 - 1.2279993570070407im & 7.50e-10 & 7.14e-01 & 1.73e-10 \\ 
		6 & 9.07 & 5 & 12.006156469238457 - 0.9155868162462955im & 6.06e-10 & 1.14e+00 & 1.57e-10 \\ 
		7 & 10.58 & 5 & 13.62772269891376 - 0.7681906372433964im & 1.53e-12 & 1.26e+00 & 4.19e-13 \\
\bottomrule
	\end{tabular}
	\caption{Results obtained by the Newton method for $n_1(r)=\sqrt{2-r^2}$, $n_2=1+(r- 0.5)^3$, $\xi=0.5$ starting from $k_0=\frac{m}{\xi n_0}$. We took here $\varepsilon=1e{-8}$ and $l_{\max} =2000$.}
	\label{tab:Newtonnonconstantvarn22}
\end{table}

\appendix

\renewcommand{\thesection}{S-Appendix \Alph{section}}

\clearpage
\section[\hspace{5em} Additional figures]{Additional figures}\label{appendix:figures}

\begin{table}[H]
	\footnotesize
	\centering
	
	\begin{subtable}{\textwidth}
		\centering
		\begin{tabular}{lccc}
			\toprule
			$l$ & $k_l$ & $|\det_1(k_l)|$ &$|\partial_k\det_1(k_l)|$ \\
			\midrule
			1& 15.2027551239541 - 0.0014194570359291362im& 4.67e+00& 3.35e+00\\
			2& 16.595070718259585 - 0.09584456745332563im& 7.67e-01& 2.26e+00  \\
			3& 16.910865265369715 - 0.21866765298920093im& 4.90e-02& 2.03e+00\\
			4& 16.92329111475096 - 0.2394557194458im& 2.56e-04& 2.02e+00 \\
			5& 16.923201859904534 - 0.239545593149898im& 7.01e-09& 2.02e+00\\
			\bottomrule
		\end{tabular}
		\caption{$\det_1$ starting from $k_0=11$.}
	\end{subtable}
	
	\begin{subtable}{\textwidth}
		\centering
		\begin{tabular}{lccc}
			\toprule
			$l$ & $k_l$ & $|\det_2(k_l)|$ &$|\partial_k\det_2(k_l)|$ \\
			\midrule
			1& 13.396229610406712 + 4.6308732874575025im& 3.28e-02& 3.63e-03\\
			2& 21.37329855682114 + 8.869886753253347im& 1.81e-02& 7.99e-04    \\
			3& 42.399466900838604 + 17.295737648975493im& 9.11e-03& 1.99e-04\\
			4& 84.49957023255998 + 34.942960837304014im& 4.56e-03& 4.99e-05   \\
			5& 168.80515382942275 + 70.15689112416922im& 2.28e-03& 1.25e-05 \\
			6& 337.50591933464125 + 140.53062554154792im& 1.14e-03& 3.12e-06\\
			7& 674.9581092588107 + 281.25208333979583im& 5.70e-04& 7.79e-07 \\ 
			8& 1349.8889995007594 + 562.6826361056818im& 2.85e-04& NaN\\
			\bottomrule
		\end{tabular}
		\caption{$\det_2$ starting from $k_0=11$.}
	\end{subtable}
	\begin{subtable}{\textwidth}
		\centering
		\begin{tabular}{lccc}
			\toprule
			$l$ & $k_l$ & $|\det_1(k_l)|$ &$|\partial_k\det_1(k_l)|$ \\
			\midrule
			1& 16.933172576735238 - 0.24982842615297063im& 2.89e-02& 2.01e+00\\
			2& 16.92321608496663 - 0.23950336537694566im& 8.99e-05& 2.02e+00   \\
			3& 16.92320186112372 - 0.23954559016258034im& 8.67e-10& 2.02e+00\\
			\bottomrule
		\end{tabular}
		\caption{$\det_1$ starting from $k_0=17$.}
	\end{subtable}
	\\
	\begin{subtable}{\textwidth}
		\centering
		\begin{tabular}{lccc}
			\toprule
			$l$ & $k_l$ & $|\det_2(k_l)|$ &$|\partial_k\det_2(k_l)|$ \\
			\midrule
			1& 16.903562544030315 - 0.2594084783861693im& 5.85e-04& 2.11e-02\\
			2& 16.922896121039905 - 0.2394985971301382im& 6.44e-06& 2.08e-02\\
			3& 16.92320187831679 - 0.2395455557149322im& 7.97e-10& 2.08e-02\\
			\bottomrule
		\end{tabular}
		\caption{$\det_2$ starting from $k_0=17$.}
	\end{subtable}
	\caption{Newton iterations  using $\det_1$ or $\det_2$ for $\xi=0.5$, $n_1=1.5$, $n_2=1$ and $m=10$, starting from $k_0=11$ or $k_0=17$.}
	\label{tab:newtontrajectories}
\end{table}

\begin{figure}[H]
	\centering
	\centering
	\begin{subfigure}[b]{0.32\textwidth}
		\centering
		\includegraphics[width=\textwidth]{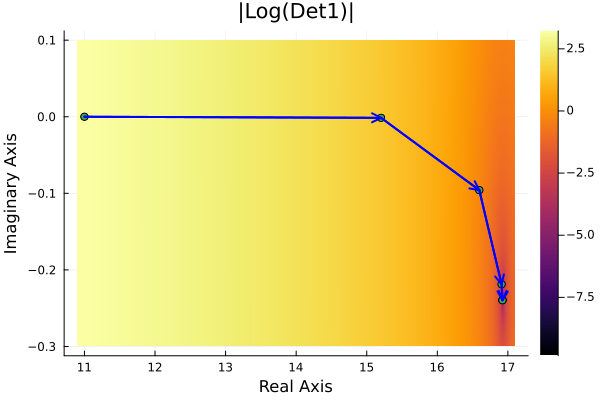}
		\caption{}
		\label{fig:k011det1}
	\end{subfigure}\quad
	\begin{subfigure}[b]{0.32\textwidth}
		\centering
		\includegraphics[width=\textwidth]{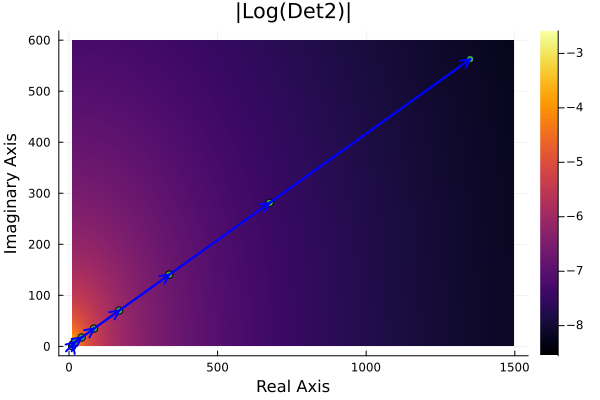}
		\caption{}
		\label{fig:k011det2}
	\end{subfigure}\\
	\vspace{0.2cm}
	\begin{subfigure}[b]{0.32\textwidth}
		\centering
		\includegraphics[width=\textwidth]{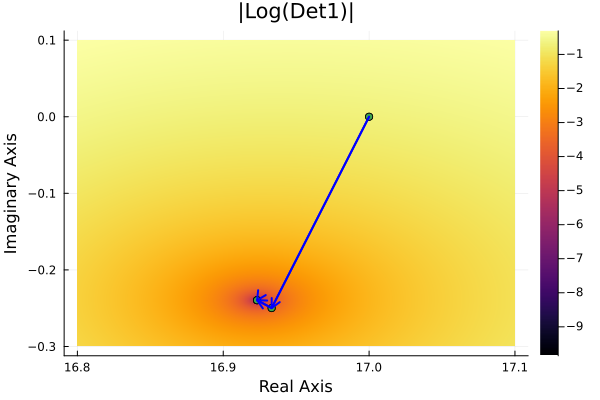}
		\caption{}
		\label{fig:k017det1}
	\end{subfigure}\quad
	\begin{subfigure}[b]{0.32\textwidth}
		\centering
		\includegraphics[width=\textwidth]{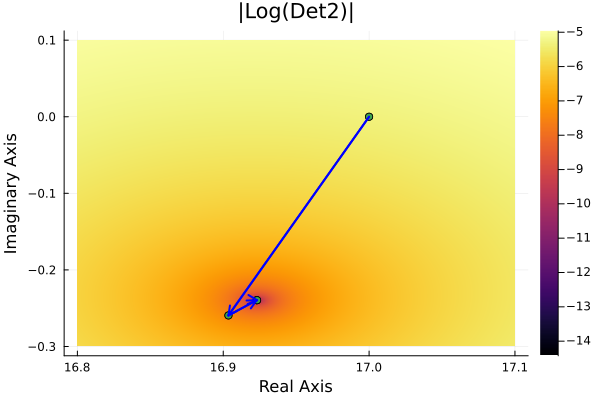}
		\caption{}
		\label{fig:k017det2}
	\end{subfigure}
	\caption{Newton iterations in the complex plane using $\det_1$ or $\det_2$ for $\xi=0.5$, $n_1=1.5$, $n_2=1$ and $m=10$, starting from $k_0=11$ (top) or $k_0=17$ (bottom).}
	\label{fig:newtontrajectories}
\end{figure}

\section[\hspace{5em} Basic properties of nonlinear eigenvalue problems~\cite{voss2013nonlinear}]{Basic properties of nonlinear eigenvalue problems~\cite{voss2013nonlinear}}\label{appendix:nep}

We gather here some definitions and basic properties of nonlinear eigenvalue problems.

We consider the problem of finding $\lambda$ such that the linear system
\begin{equation}
T(\lambda) x=0
\end{equation}
has a nontrivial solution $x$, where $T(.)\colon D\to \CC^{n\times n}$ is a family of matrices depending on a complex parameter $\lambda \in D$.

\bigskip
\noindent\textbf{Definitions:}

$\hat\lambda\in D$ is called an eigenvalue of $T(\cdot)$ if $T(\hat\lambda)x=0$ has a nontrivial solution $\hat x\ne 0$. $\hat x$ is called a corresponding eigenvector or right eigenvector, and $(\hat\lambda,\hat x)$ is called eigenpair of $T(\cdot)$.

Any nontrivial solution $\hat y\ne 0$ of the adjoint equation $T(\hat\lambda)^*y=0$ is called left eigenvector of $T(\cdot)$ and the vector-scalar-vector triplet $(\hat y,\hat\lambda,\hat x)$ is called an eigentriplet of $T(\cdot)$.

An eigenvalue $\hat\lambda$ of $T(\cdot)$ has algebraic multiplicity $k$ if $\frac{d^\ell}{d\lambda^\ell}\det{(T(\lambda))}|_{\lambda=\hat\lambda}=0$ for $\ell=0,\ldots,k-1$ and $\frac{d^k}{d\lambda^k}\det{(T(\lambda))}|_{\lambda=\hat\lambda}\ne 0$.

An eigenvalue $\hat\lambda$ is simple if its algebraic multiplicity is one.

The geometric multiplicity of an eigenvalue $\hat\lambda $ is the dimension of the kernel $\ker(T(\hat\lambda))$ of $T(\hat\lambda)$.

An eigenvalue is called semi-simple if its algebraic and geometric multiplicity coincide.

\bigskip
\noindent\textbf{Properties:}

If $\hat\lambda$ is an algebraically simple eigenvalue of $T(\cdot)$, then $\hat\lambda$ is geometrically simple.~\cite{schreiber2008nonlinear}

Let $(\hat y,\hat\lambda,\hat x)$ be an eigentriplet of $T(\cdot)$. Then $\hat\lambda$ is algebraically simple if and only if $\hat{\lambda}$ is geometrically simple and $\hat y^* T'(\hat \lambda) \hat x\ne 0$.~\cite{neumaier1985residual,schreiber2008nonlinear}

\section[\hspace{5em} Asymptotic expansions of quasi-resonances in dimension 2~\cite{balac2021asymptotics}]{Asymptotic expansions of quasi-resonances in dimension 2~\cite{balac2021asymptotics}}\label{appendix:dauge}

We gather here some of the results in~\cite{balac2021asymptotics} regarding the construction of quasimodes $(\underline{k}(m),\underline{u}(m))$ as $m\to\infty$, i.e. approximate solutions of~\eqref{eq:helmholtzpolarcordinates} 
with $\hat{g}_m=0$, in three different configurations. 
Moreover, it is shown in~\cite{balac2021asymptotics} that the constructed quasi-resonances
$\underline{k}_j(m)$ are real and positive and close to true resonances $k_j(m)$ modulo
a super-algebraic error $O(m^{-\infty})$.

\begin{customassump}{1.1}
	The radial function $n\colon r\mapsto n(r)$ satisfies the following properties:
	\begin{enumerate}
		\item $n(r)=1$ if $r>\xi$ ;
		\item The function $n$ belongs to $C^\infty([0,\xi])$ and $n(r)>1$ for all $r\le \xi$.
	\end{enumerate}
\end{customassump}

\begin{customnot}{1.2}
	\begin{equation}
	n_0=\lim\limits_{r\nearrow\xi} n(r), \quad n_I=\lim\limits_{r\nearrow\xi} n'(r), \quad n_{II}=\lim\limits_{r\nearrow\xi} n''(r).
	\end{equation}
	
	Effective adimensional curvature:
	\begin{equation}
	\check\kappa:=\xi\Bigl(\frac1\xi+\frac{n_I}{n_0}\Bigr).
	\end{equation}
	
	Adimensional Hessian:
	\begin{equation}
	\check\mu:=\xi^2\Bigl(\frac2{\xi^2}-\frac{n_{II}}{n_0}\Bigr).
	\end{equation}
	
\end{customnot}

\begin{customthm}{1.A}\label{th:1.A}
	Assume the radial function $n$ satisfies Assumption~1.1 and
	\begin{equation}
	\check\kappa>0.
	\end{equation}
	Then, for any integer $j\ge 0$, there exists a quasi-pair $(\underline{k}_j(m),\underline{u}_j(m))$ such that the quasi-resonance $\underline{k}_j(m)$ has an expansion in integer powers of $m^{-1/3}$ starting as
	\begin{align}
		\underline{k}_j(m)=\frac{m}{\xi n_0}\Bigl[&1+\frac{a_j}2\Bigl(\frac{2\check\kappa}{m}\Bigr)^{\frac23}-\frac{n_0}{2\sqrt{n_0^2-1}}\Bigl(\frac{2\check\kappa}{m}\Bigr)+\frac{a_j^2}{15}\Bigl(\frac{17}8-\frac{3}{\check\kappa}+\frac{\check\mu}{\check\kappa^2}\Bigr)\Bigl(\frac{2\check\kappa}{m}\Bigr)^{\frac43}\notag\\ 
		&\quad -\frac{a_j n_0}{12\sqrt{n_0^2-1}}\Bigl(\frac{n_0^2}{n_0^2-1}+2-\frac{6}{\check\kappa}+\frac{2\check\mu}{\check\kappa^2}\Bigr)\Bigl(\frac{2\check\kappa}{m}\Bigr)^{\frac53}+O(m^{-2})\Bigr]
	\end{align}
	with the numbers $a_j$ being the successive roots of the flipped Airy function  $z \in\CC\mapsto \mathrm{Ai}(-z)$ where $\mathrm{Ai}$ denotes the Airy function~\cite{abramowitz1968handbook,fabijonas2004computation}.
\end{customthm}

\begin{customthm}{1.B}\label{th:1.B}
	Assume the radial function $n$ satisfies Assumption~1.1 and
	\begin{equation}
	\check\kappa=0 \quad \text{with} \quad \check\mu>0.
	\end{equation}
	Then, for any integer $j\ge 0$, there exists a quasi-pair $(\underline{k}_j(m),\underline{u}_j(m))$ such that the quasi-resonance $\underline{k}_j(m)$ has an expansion in integer powers of $m^{-1/2}$ starting as
	\begin{equation}
	\underline{k}_j(m)=\frac{m}{\xi n_0}\Bigl[1+\frac{4j+3}2\Bigl(\frac{\sqrt{\check\mu}}{m}\Bigr)+O(m^{-\frac32})\Bigr],
	\end{equation}
	the coefficient of degree $1$ being zero.
\end{customthm}

\begin{customthm}{1.C}\label{th:1.C}
	Assume the radial function $n$ satisfies Assumption~1.1 and that 
	\begin{equation}
	\check\kappa<0.
	\end{equation}
	Let $\xi_0\in(0,\xi)$ such that $1+\frac{\xi_0n'(\xi_0)}{n(\xi_0)}=0$ and assume further that
	\begin{equation}
	\check\mu_0 := \xi_0^2\Bigl(\frac2{\xi_0^2}-\frac{n''(\xi_0)}{n(\xi_0)}\Bigr)>0.
	\end{equation}
	Then, for any integer $j\ge 0$, there exists a quasi-pair $(\underline{k}_j(m),\underline{u}_j(m))$ such that the quasi-resonance $\underline{k}_j(m)$ has an expansion in integer powers of $m^{-1/2}$ starting as
	\begin{equation}\label{eq:th1C}
	\underline{k}_j(m)=\frac{m}{\xi_0 n(\xi_0)}\Bigl[1+\frac{2j+1}2\Bigl(\frac{\sqrt{\check\mu_0}}{m}\Bigr)+O(m^{-2})\Bigr],
	\end{equation}
	the coefficients of degree $1$ and $3$ being zero.
	
\end{customthm}

\section[\hspace{5em} Special variable $n$ with explicit solution (Luneburg lens)]{Special variable $n$ with explicit solution (Luneburg lens)}\label{appendix:Whittaker}

The fundamental solutions of the second-order differential equation
\begin{equation}\label{eq:luneburg}
-y''(x)-\frac{y'(x)}{x}+\Bigl(\frac{m^2}{x^2}-k^2 (2-x^2)\Bigr)y(x)=0
\end{equation}
(where $n(x)=\sqrt{2-x^2}$) are given by
\begin{itemize}
	\item  $ \displaystyle f_1(x)=\frac1x M_{\frac{k}2,\frac{m}2}(k x^2)$ (vanishes at $x=0$)
	\item   $\displaystyle f_2(x)=\frac1x W_{\frac{k}2,\frac{m}2}(k x^2)$ (singular at $x=0$)
\end{itemize}
where $M_{k,m}(z)$ and $W_{k,m}(z)$ are the Whittaker functions~\cite{abramowitz1968handbook,whittaker1920course,olver2016nist}, solutions of the Whittaker differential equation
\begin{equation}
w''+\Big(\frac{1/4-m^2}{z^2}+\frac{k}z-\frac14\Bigr)w=0.
\end{equation}
Indeed, if one looks for a solution of~\eqref{eq:luneburg} in the form $y(x)=\frac1x w(kx^2)$ then using
\begin{align}
	y'(x)&=2k w'(kx^2)-\frac1{x^2}w(kx^2)\\
	y''(x)&=4k^2 x w''(kx^2)-\frac{2k}x w'(kx^2)+\frac2{x^3} w(kx^2),
\end{align}
the differential equation~\eqref{eq:luneburg} becomes
\begin{align}
	-4k^2 x w''(kx^2)-\Bigl(\frac{1-m^2}{x^3}+\frac{k^2(2-x^2)}x\Bigr)
	w(kx^2)&=0\\
	w''(kx^2)+\Bigl(\frac{1-m^2}{4k^2x^4}+\frac{1}{2x^2}-\frac{1}{4}\Bigr)w(kx^2)&=0,
\end{align}
hence the announced result.

We have used the Whittaker representation in our implementation for its convenience and reduced form.

This special case in linked to the Luneburg lens which is a spherically symmetric gradient-index lens, whose typical refractive index has the same form considered here~\cite{boriskin2002whispering,lock2008scattering1,lock2008scattering2,lock2008scattering3}.

\end{document}